\documentclass[12pt]{amsart}
\usepackage[all, knot]{xy}
\usepackage{epsfig}
\usepackage{xspace}
\usepackage{layout}
\usepackage{latexsym}
\usepackage{amsmath}
\usepackage{amsthm}
\usepackage{amssymb}
\usepackage{amsfonts}
\usepackage{amsxtra}     
\usepackage{url}
\usepackage{verbatim}
\usepackage{longtable}

\newcommand{\ie}{\textit{i.e., }}
\newcommand{\eg}{\textit{e.g., }}

\newcommand{\vcomp}[2]{\genfrac{[}{]}{0pt}{}{#1}{#2}}
\newcommand{\bbA}{\mathbb{A}}
\newcommand{\bbB}{\mathbb{B}}
\newcommand{\bbC}{\mathbb{C}}
\newcommand{\bbD}{\mathbb{D}}
\newcommand{\bbE}{\mathbb{E}}
\newcommand{\bbF}{\mathbb{F}}

\newcommand{\bbH}{\mathbb{H}}
\newcommand{\bbP}{\mathbb{P}}
\newcommand{\bbQ}{\mathbb{Q}}
\newcommand{\bbR}{\mathbb{R}}
\newcommand{\bbS}{\mathbb{S}}
\newcommand{\bbV}{\mathbb{V}}
\newcommand{\bbX}{\mathbb{X}}

\newcommand{\bfA}{\mathbf{A}}
\newcommand{\bfB}{\mathbf{B}}
\newcommand{\bfC}{\mathbf{C}}
\newcommand{\bfD}{\mathbf{D}}
\newcommand{\bfE}{\mathbf{E}}
\newcommand{\bfF}{\mathbf{F}}

\newcommand{\bfH}{\mathbf{H}}
\newcommand{\bfI}{\mathbf{I}}
\newcommand{\bfM}{\mathbf{M}}
\newcommand{\bfP}{\mathbf{P}}
\newcommand{\bfQ}{\mathbf{Q}}
\newcommand{\bfR}{\mathbf{R}}
\newcommand{\bfU}{\mathbf{U}}
\newcommand{\bfV}{\mathbf{V}}
\newcommand{\bfX}{\mathbf{X}}

\newcommand{\Obj}{\text{\it Obj}\;}
\newcommand{\Mor}{\text{\it Mor}\;}
\newcommand{\Hor}{\text{\it Hor}\;}
\newcommand{\Ver}{\text{\it Ver}\;}

\newcommand{\Sq}{\text{\it Sq}\;}
\newcommand{\colim}{\text{\rm colim}\;}
\newcommand{\im}{\text{\it Im}}

\newcommand{\dom}{\text{\it dom}\;}
\newcommand{\fib}{\text{\sf fib}}
\newcommand{\cof}{\text{\sf cof}}
\newcommand{\we}{\text{\sf we}}
\newcommand{\Cat}{\text{\sf Cat}}
\newcommand{\iiso}{\text{\sf iso}}
\newcommand{\Grpd}{\text{\sf Grpd}}
\newcommand{\diag}{\text{\rm diag}}
\newcommand{\Sd}{\text{\rm Sd}}
\newcommand{\Ex}{\text{\rm Ex}}

\newcommand{\co}{\colon\thinspace}

\newcommand{\tb}[1]{\phantom{\sum^\Sigma_\Sigma} #1 \phantom{\sum^\Sigma_\Sigma}}
\newcommand{\lr}[1]{\hspace{.5mm}#1\hspace{.5mm}}

\begin{document}

\newtheorem{thm}{Theorem}[section]
\newtheorem{conj}[thm]{Conjecture}
\newtheorem{lem}[thm]{Lemma}
\newtheorem{cor}[thm]{Corollary}
\newtheorem{prop}[thm]{Proposition}
\newtheorem{rem}[thm]{Remark}

\theoremstyle{definition}
\newtheorem{defn}[thm]{Definition}
\newtheorem{examp}[thm]{Example}
\newtheorem{notation}[thm]{Notation}
\newtheorem{rmk}[thm]{Remark}

\theoremstyle{remark}

\makeatletter
\renewcommand{\maketag@@@}[1]{\hbox{\m@th\normalsize\normalfont#1}}%
\makeatother
\renewcommand{\theenumi}{\roman{enumi}}

\def\qed {{
   \parfillskip=0pt        
   \widowpenalty=10000     
   \displaywidowpenalty=10000  
   \finalhyphendemerits=0  
  %
   \leavevmode             
   \unskip                 
   \nobreak                
   \hfil                   
   \penalty50              
   \hskip.2em              
   \null                   
   \hfill                  
   $\square$
  %
   \par}}                  

\newenvironment{pf}{{\it Proof:}\quad}{\qed \vskip 12pt}


\title[Model Structures on {\bf DblCat}]{Model Structures on the Category of Small Double Categories}
\author[Fiore, Paoli, Pronk]{Thomas M. Fiore, Simona Paoli, Dorette A. Pronk}
\address{Thomas M. Fiore \\ Department of Mathematics \\ University of Chicago \\ Chicago, IL 60637 \\ USA \\
and Departament de Matem\`{a}tiques \\ Universitat Aut\`{o}noma de Barcelona \\
08193 Bellaterra (Barcelona)  \\ Spain}
\email{fiore@math.uchicago.edu}
\address{Simona Paoli \\ Department of Mathematics \\ Macquarie University \\ NSW 2109 \\ Australia}
\email{simonap@maths.mq.edu.au}
\address{Dorette Pronk \\ Department of Mathematics and Statistics \\ Dalhousie
University \\ Halifax, NS B3H 3J5 \\ Canada}
\email{pronk@mathstat.dal.ca}



\begin{abstract}
In this paper we obtain several model structures on {\bf DblCat},
the category of small double categories. Our model structures have
three sources. We first transfer across a categorification-nerve
adjunction. Secondly, we view double categories as internal
categories in {\bf Cat} and take as our weak equivalences various
internal equivalences defined via Grothendieck topologies. Thirdly,
{\bf DblCat} inherits a model structure as a category of algebras
over a 2-monad. Some of these model structures coincide and the
different points of view give us further results about cofibrant
replacements and cofibrant objects. As part of this program we give
explicit descriptions and discuss properties of free double
categories, quotient double categories, colimits of double
categories, several nerves, and horizontal categorification.
\end{abstract}

\maketitle

\tableofcontents



\section{Introduction} \label{section:Introduction}
The theory of categories enriched in {\bf Cat}, called 2-categories,
has been highly developed over the past 40 years and has found
numerous applications. Beginning with B\'enabou's bicategories (weak
2-categories) in \cite{benabou}, through Kelly's monograph
\cite{kelly3} on enriched categories, and including the more recent
article \cite{lackcompanion}, as well as many others, we have seen
the $n=2$ case for higher category theory become very well
understood. Limits in 2-categories \cite{kelly2}, 2-monads on
2-categories \cite{blackwellkellypower}, and Kan extensions for
2-functors \cite{dubuc} are now widely known. Model structures on
{\bf 2-Cat} have also been studied recently in \cite{lack2Cat},
\cite{lackBiCat}, and \cite{worytkiewicz2Cat}. Model structures,
more generally, have been used in the study of
$(\infty,1)$-categories as a means of comparison by Bergner,
Joyal-Tierney, Rezk, and To{\"e}n \cite{bergnerthreemodels},
\cite{bergnersimplicialcategories}, \cite{joyaltierneyquasisegal},
\cite{rezkhomotopytheory}, and \cite{toenaxiomatization}.

Recent examples, however, show that 2-categories are not enough, and
that one must invoke Ehresmann's earlier notion of {\it double}
category \cite{ehresmann}, \cite{ehresmann2}. In many mathematical
situations one is interested in two types of morphisms, which may or
may not interact. Between rings, for example, there are ring
homomorphisms as well as bimodules. Between manifolds there are
diffeomorphisms and cobordisms, which are both used in field theory.
Between categories there are functors as well as adjunctions. The
notion of 2-category does not capture both types of morphisms, but
the notion of (pseudo) double category certainly does.

Concisely, a {\it small double category} is an internal category in
{\bf Cat}. A small double category consists of a set of objects, a
set of horizontal morphisms, a set of vertical morphisms, and a set
of squares, equipped with various associative and unital
compositions satisfying the interchange law. In addition to the
early work of Bastiani-Ehresmann \cite{ehresmannone},
A.~Ehresmann-C.~Ehresmann \cite{ehresmanntwo},
\cite{ehresmannthree}, \cite{ehresmannfour}, C.~Ehresmann
\cite{ehresmann}, \cite{ehresmann2}, and Brown-Spencer
\cite{brownspencer76}, recent work on double categories has been
completed by Brown and collaborators, Dawson, Fiore, Garner,
Grandis, Kock, Par\'e, Pronk, Shulman, and others
\cite{brownicen2003}, \cite{brownmackenzie}, \cite{brownmosa99},
\cite{dawsonparetiling}, \cite{dawsonparepronkspan},
\cite{dawsonpareassociativity}, \cite{dawsonparefreedouble},
\cite{fiore2}, \cite{garner2005}, \cite{grandisdouble3},
\cite{grandisdouble4}, \cite{grandisdouble1}, \cite{grandisdouble2},
\cite{kock}, \cite{shulmanonquillenfunctors}, and
\cite{shulmanframed}.

Double categories are the $n=2$ case for $n$-fold categories, which
have been studied and applied for some time now. In the same way
that {\it higher} categories may be defined by iterated enrichment,
one may define {\it wider} categories or {\it $n$-fold categories}
via iterated internalization. The edge symmetric\footnote{{\it Edge
symmetric} means that the $n$-morphisms in all $n+1$ directions are
the same.} case has been studied by Brown under the name of cubical
$\omega$-categories. Further, $n$-fold categories internal to the
category of groups have been used to model connected homotopy
$(n+1)$-types in \cite{lodayfinitelymany} as summarized in the
survey paper \cite{paoliinternalstructures}. Recent work includes
\cite{lackpaolioperadic} and \cite{paoliontamsamani}. Applications
of versions of the $n=2$ case of internalized categories include
\cite{dawsonparepronkpaths}, \cite{fiore1}, \cite{fiore2},
\cite{kerlerlyubashenko}, \cite{maysigurdsson}, \cite{mortondouble},
\cite{shulmanonquillenfunctors}, \cite{shulmanframed}. Thus, there
has been a general trend towards $n$-fold categories, especially the
$n=\omega$ and $n=2$ cases.

In this article we introduce model categories into the theory of
double categories, anticipating a utility in the theory of {\it
wider categories} analogous to that of model structures in the
theory of {\it higher categories}. Already in the $n=2$ case we see
that $n$-fold categories and $n$-categories diverge: even though the
homotopy theory of 2-categories resembles that of categories, the
homotopy theory of double categories is much richer. This results
from the numerous ways to view a double category: as an internal
category in {\bf Cat}, as a categorical structure with two
directions, as certain simplicial objects in {\bf Cat}, as certain
bisimplicial sets, or as algebras over a 2-monad. Each point of view
suggests different notions of weak equivalence and fibration. The
new types of pasting diagrams available in a double category also
create new phenomena. We take these various points of view into
consideration when constructing the model structures.

Thus, our model structures have three sources. First, we transfer
the categorical diagram structure and Thomason diagram structure on
the category of simplicial objects in {\bf Cat} to {\bf DblCat} via
a horizontal categorification-horizontal nerve adjunction. In the
Thomason structure on {\bf Cat} in \cite{thomasonCat}, a functor is
a weak equivalence if and only if its nerve is a weak homotopy
equivalence of simplicial sets. In the categorical structure on {\bf
Cat} of \cite{joyaltierney} and \cite{rezkcat}, a functor is a weak
equivalence if and only if it is an equivalence of categories. Both
the Thomason structure and the categorical structure on
$\mathbf{Cat}$ are cofibrantly generated, and thus induce
cofibrantly generated model structures on simplicial objects in {\bf
Cat} where weak equivalences and fibrations are defined levelwise.
We apply Kan's Lemma on Transfer of cofibrantly generated model
structures (Theorem \ref{Kan}) to transfer both of these diagram
structures to {\bf DblCat} (Theorems \ref{diagramtransfer} and
\ref{categoricaltransfer}). However, the application is not
straightforward, and we must make several double categorical
preparations, including horizontal categorification and a pushout
formula in {\bf DblCat}. We also prove one negative result in
Theorem \ref{NoReedy}: it is impossible to transfer the Reedy
categorical structure on $\mathbf{Cat}^{\Delta^{op}}$ to {\bf
DblCat}. The transfer from bisimplicial sets will be the subject of
a later article.

We arrive at a second source for model structures on {\bf DblCat}
when we view double categories as internal categories in {\bf Cat}.
In this way we obtain double categorical versions of the categorical
structure on {\bf Cat}, where a functor is a weak equivalence if and
only if it is fully faithful and essentially surjective. Although
the notion of fully faithfulness makes sense internally, essential
surjectivity does not, and therefore equivalences of internal
categories need further explanation. Model structures on categories
internal to a good category $\bfC$ have already been developed in
\cite{everaertinternal}, and we apply their results to the case
$\bfC=\mathbf{Cat}$. They define essential surjectivity (and hence
also weak equivalences) with respect to a Grothendieck topology
$\mathcal{T}$ on $\bfC$. We take simplicially surjective functors
and categorically surjective functors as bases for Grothendieck
topologies on {\bf Cat}, and obtain two distinct model structures in
Sections \ref{subsection:simpliciallysurjective} and
\ref{subsection:categoricallysurjective}. Additionally, we show in
Theorem \ref{trivialtopologyproducestrivialstructure} that the model
structure induced by the trivial topology coincides with the trivial
model structure from the 2-category of internal categories.

Third, {\bf DblCat} inherits a model structure as a category of
algebras over a 2-monad as in \cite{lack2monads}. The underlying
1-category of a 2-category with finite limits and finite colimits
always admits the so-called trivial model structure, whose weak
equivalences are equivalences and fibrations are isofibrations. If
$\mathcal{K}$ is a locally finitely presentable 2-category equipped
with a 2-monad $T$ with rank, then the category of (strict)
$T$-algebras is a model category: a morphism of $T$-algebras is a
weak equivalence or fibration if and only if its underlying morphism
is a weak equivalence or fibration in the trivial model structure on
$\mathcal{K}$. In our application of \cite{lack2monads} in Section
\ref{2monadstructure}, $\mathcal{K}$ is the 2-category $\Cat({\bf
Graph})$ of internal categories in small non-reflexive graphs, and
$T$ is the 2-monad $\overline{M}$ induced by the Cartesian monad $M$
on $\mathbf{Graph}$ whose algebras are categories.

Depending on the reader's experience, certain model structures will
be of more interest than others. Simplicially-minded readers will no
doubt find the transfers from $\mathbf{Cat}^{\Delta^{op}}$ most
interesting, while categorically minded readers may find the model
structures on $\Cat(\mathbf{Cat})$ arising from
\cite{everaertinternal} more interesting. Universal algebraists may
find the third point of view most appealing, namely double
categories as algebras for a 2-monad. Nevertheless, certain model
structures can be defined from two points of view, and will thus be
of interest to readers working in different fields.

In other words, we prove that some of these model structures
coincide. The model structure obtained by transferring the
categorical diagram structure across the vertical
categorification-vertical nerve adjunction is the same as the model
structure associated to the simplicially surjective topology on {\bf
Cat} (Corollary \ref{simpliciallysurjectivecoincidence}). The
algebra structure is the same as the model structure associated to
the categorically surjective topology on {\bf Cat} (Theorem
\ref{2monadstructure=tau'structure}).

These two different constructions of the same model structures yield
more refined information about cofibrant replacements and cofibrant
objects. For example, the cofibrant objects in the algebra structure
are known to be precisely the flexible algebras, but from the
categorically-surjective-topology structure we see that the flexible
double categories are precisely those with object category free on a
graph (Corollary \ref{flexibledoublecategories} and Remark
\ref{alternativetocorollary}). Such a description allows us to
conclude that the flexible 2-categories of \cite{lack2Cat} are
indeed flexible algebras for a 2-monad. Lack's Theorem 4.8 (iv) in
\cite{lack2Cat}, which characterizes flexible 2-categories as those
2-categories with underlying 1-category a free category on a graph,
now extends to double categories.

We also compare our model structures on $\mathbf{DblCat}$ with the
analogous ones for $\mathbf{Cat}$ in Propositions
\ref{ThomasonExtension}, \ref{CategoricalExtension}, and
\ref{TrivialExtension}. The vertical embedding of $\mathbf{Cat}$
into $\mathbf{DblCat}$ preserves and reflects weak equivalences,
fibrations, and cofibrations from the Thomason structure into the
transferred diagram Thomason structure, as well as from the
categorical structure into the transferred diagram categorical
structure. The horizontal embedding of $\mathbf{Cat}$ into
$\mathbf{DblCat}$ preserves and reflects weak equivalences,
fibrations, and cofibrations from the categorical structure into the
trivial structure. However, the vertical inclusion of {\bf 2-Cat}
into $\mathbf{DblCat}$ preserves neither the weak equivalences nor
the cofibrations of the categorical structure into the 2-monad
structure, as shown at the end of Section \ref{2monadstructure}.
Nevertheless, a 2-category is cofibrant in {\bf 2-Cat} if and only
if its vertical embedding into $\mathbf{DblCat}$ is cofibrant.

In order to build our model structures we prove various general
results about double categories, so far not available in the
literature. These results are also of independent interest for the
theory of double categories in its own right. We develop free double
categories, their quotients, and colimits of double categories using
a double categorical version of Street's 2-categorical notion of
{\it derivation scheme} \cite{streetstructures}. In particular we
obtain an explicit formula for two pushouts of double categories in
Theorem \ref{pushoutinDblCat}, which is essential for our
application of Kan's Lemma on Transfer in Theorems
\ref{diagramtransfer} and \ref{categoricaltransfer}. We also prove
that the 2-categories $\mathbf{DblCat_v}$ and $\mathbf{DblCat_h}$
are 2-cocomplete in Theorem \ref{2cocomplete}.

Free double categories on reflexive double graphs have been studied
in \cite{dawsonparefreedouble}. By {\it reflexive double graph} we
mean a collection of objects, vertical edges, horizontal edges, and
squares equipped with source and target maps, identity edges, and
identity squares. In this paper, we will instead use double graphs
with 1-identities. A {\it double graph with 1-identities} is like a
reflexive double graph, except there are no identity squares.
Between double graphs with 1-identities and double categories, there
is the intermediate notion of double derivation scheme.  A {\it
double derivation scheme} is a double graph with 1-identities in
which the horizontal and vertical reflexive 1-graphs are categories.
In the free double category on a double derivation scheme, the
vertical and horizontal 1-categories are preserved, but nontrivial
squares consist of allowable compatible arrangements. Since we are
considering compatible arrangements of squares in a double
derivation scheme rather than in a double reflexive graph, our {\it
allowable} compatible arrangements are different than the {\it
composable} compatible arrangements of
\cite{dawsonpareassociativity}.

Free double categories on double derivation schemes and their
quotients allow us to construct colimits of double categories. First
one takes the colimits of the vertical and horizontal 1-categories.
These, together with the colimit of the {\it sets} of squares, form
a double derivation scheme. Finally, we mod out the free double
category on this double derivation scheme by the smallest congruence
which guarantees that the natural maps are double functors, and the
result is the colimit in {\bf DblCat}. This colimit formula is the
basis of Theorem \ref{pushoutinDblCat} which gives an explicit
description of the pushouts of a double functor along two inclusions
of external products. This theorem is crucial for our application of
Kan's Lemma on Transfer. These two pushouts are special cases of a
more general theorem on pushouts along inclusions of external
products, which will appear in a separate paper with a comparison to
\cite{dawsonparepronkextensions}.

Free double categories on double derivation schemes and their
quotients find further application in the construction of
fundamental double categories of simplicial objects in {\bf Cat},
\ie in our construction of a left adjoint to the horizontal nerve.
We obtain an important example of our explicit constructions of
fundamental double categories in a second way as well, namely via
weighted colimits (see Example \ref{horcatcatssetfromformula} and
Proposition \ref{horcatcatsset}).

We begin in Section \ref{section:DoubleCategories} with a review of
double categories, including horizontal 2-categories, vertical
2-categories, double functors, horizontal and vertical natural
transformations, the external product of 2-categories, and Cartesian
closedness of the category $\mathbf{DblCat}$, as well as the
2-categories $\mathbf{DblCat_v}$ and $\mathbf{DblCat_h}$. Free
double categories on double derivation schemes are introduced in
Section \ref{section:FreeDouble} and are used in Section
\ref{section:limitscolimits} to describe colimits in {\bf DblCat}.
Horizontal and double nerves are discussed in Section
\ref{section:Nerves} along with their representable definitions in
terms of external products of finite ordinals. In Section
\ref{section:categorification}, free double categories on double
derivation schemes and their quotients are applied to construct the
left adjoint to the horizontal nerve. Section
\ref{section:simplicialcategories} focuses on transferring model
structures across the horizontal categorification-horizontal nerve
adjunction, and recalls model structures on {\bf Cat}, smallness
issues, and Kan's Lemma on Transfer. Section
\ref{section:topologies} begins with an exposition of the methods of
\cite{everaertinternal}, and then applies them to obtain model
structures on $\Cat({\bf Cat})={\bf DblCat}$ induced by three
Grothendieck topologies on {\bf Cat}: the simplicially surjective
topology, the categorically surjective topology, and the trivial
topology. The model structure induced by the simplicially surjective
topology coincides with the transfer of the diagram categorical
structure across the adjunction $c_v \dashv N_v$. In Section
\ref{2monadstructure} we prove that the 2-monad structure on {\bf
DblCat} coincides with the model structure induced by the
categorically surjective topology. In Section
\ref{section:appendix}, the Appendix, we obtain an explicit
description of certain pushouts in {\bf DblCat}, namely Theorem
\ref{pushoutinDblCat}. We use this to characterize the behavior of
the horizontal nerve on such pushouts in Theorem \ref{pushoutiso}.
The essential application is to the generating acyclic cofibrations
in the transfer in Section \ref{section:simplicialcategories}.

\vspace{.5in} {\bf Acknowledgements.} Thomas M. Fiore was supported
by National Science Foundation Grant DMS 0501208 at the University
of Chicago. At the Universitat Aut\`{o}noma de Barcelona he was
supported by Grant SB2006-0085 of the Spanish Ministerio de
Educaci\'{o}n y Ciencia under the Programa Nacional de ayudas para
la movilidad de profesores de universidad e investigadores
espa$\tilde{\text{n}}$oles y extranjeros. Simona Paoli was supported
by an Australian Research Council Postdoctoral Fellowship (project
number DP0558598) and by a Macquarie University New Staff Grant
Scheme. Dorette Pronk was supported by an NSERC Discovery Grant, and
she also thanks Macquarie University and the University of Chicago
for their hospitality and financial support, as well as Calvin
College and Utrecht University for their hospitality. All three
authors gratefully acknowledge the support and hospitality of the
Fields Institute during the Thematic Program on Geometric
Applications of Homotopy Theory, at which a significant portion of
this work was completed. The authors thank Steve Lack for suggesting
the comparison of the model structure induced by the categorically
surjective topology and the model structure induced by a 2-monad.
They also thank Michael Shulman for the simplified proof of
Corollary \ref{likeinternalhom}. The authors express their gratitude
to Peter May, Robert Par{\'e}, and Robert Dawson for some discussion
of this work. Additionally, they thank the anonymous referee for
many helpful suggestions.

\section{Double Categories} \label{section:DoubleCategories}

We first recall the elementary notions of double category theory. In
many mathematical contexts there are two interesting types of
morphisms; double categories organize them into one structure. For
example, between rings there are morphisms of rings as well as
bimodules, between objects of any 2-category there are morphisms as
well as adjunctions, and so on. Sometimes one would like to
distinguish a family of squares, such as the pullback squares among
the commutative squares, and double categories are also of use here.
The notion of double category is not new, and goes back to Ehresmann
in \cite{ehresmann} and \cite{ehresmann2}.

\begin{defn} \label{defnofdoublecategory}
A {\it small double category}
$\mathbb{D}=(\mathbb{D}_0,\mathbb{D}_1)$ is a category object in the
category of small categories. This means that $\mathbb{D}_0$ and
$\mathbb{D}_1$ are categories equipped with functors
$$\xymatrix@C=3pc{\mathbb{D}_1 \times_{\mathbb{D}_0} \mathbb{D}_1
\ar[r]^-m & \mathbb{D}_1 \ar@/^1pc/[r]^s \ar@/_1pc/[r]_t &
\ar[l]|{\lr{u}} \mathbb{D}_0 }$$ that satisfy the usual axioms of a
category. We call the objects and morphisms of $\mathbb{D}_0$
respectively the {\it objects} and {\it vertical morphisms} of
$\mathbb{D}$, and we call the objects and morphisms of
$\mathbb{D}_1$ respectively the {\it horizontal morphisms} and {\it
squares} of $\mathbb{D}$.
\end{defn}

When one expands this definition, one sees that a small double
category consists of a set of {\it objects}, a set of {\it
horizontal morphisms}, a set of {\it vertical morphisms}, and a set
of {\it squares} equipped with various sources, targets, and
associative and unital compositions. Since we only deal with small
double categories, we will usually leave off the adjective small.
Sources and targets are indicated as follows.
\begin{equation} \label{sourcetarget}
\xymatrix{A \ar[r]^f &  B & A \ar[d]_j
& A \ar[r]^f \ar[d]_j \ar@{}[dr]|\alpha & B \ar[d]^k \\
& & C & C \ar[r]_g & D}
\end{equation}
We denote the
set of squares with the boundary
$$\xymatrix{A \ar[r]^f \ar[d]_j & B \ar[d]^k \\ C
\ar[r]_g & D}$$ by $\bbD\begin{pmatrix} & f &  \\ j & & k
\\ & g & \end{pmatrix}$.
Then one has the categories
$$(\Obj \bbD, \Hor \bbD)\mbox{ and }(\Ver \bbD, \Sq \bbD)$$ under
horizontal composition and the categories
$$(\Obj \bbD, \Ver \bbD)\mbox{ and }(\Hor \bbD, \Sq \bbD)$$
under vertical composition. We will write $[f\  g]$ for the
horizontal composition of horizontal morphisms $f$ and $g$, and
similarly $[\alpha\ \beta]$ for the horizontal composition of
squares $\alpha$ and $\beta$. We will write $\vcomp{v}{w}$ for the
vertical composition of vertical morphisms $v$ and $w$, and
similarly $\vcomp{\alpha}{\beta}$ for the vertical composition of
squares $\alpha$ and $\beta$. Composition of squares in $\bbD$
satisfies the usual interchange law.

There are many examples of double categories. The commutative
squares in a given 1-category form a double category. More
generally, for a 2-category $\bfC$, Ehresmann defined the double
category $\mathbb{Q}\mathbf{C}$ of {\it quintets of $\mathbf{C}$}.
Its objects are the objects of $\mathbf{C}$, horizontal and vertical
morphisms are the morphisms of $\mathbf{C}$, and the squares
$\alpha$ as in (\ref{sourcetarget}) are the 2-cells
$\xymatrix@1{\alpha\co k \circ f \ar@{=>}[r] & g \circ j }$. In many
situations, one has examples of a slightly more general notion
called {\it pseudo double category}, defined in
\cite{grandisdouble1}. This is like a double category, except one
direction is a bicategory (weak 2-category) rather than a
2-category. For example, the double category of rings, bimodules,
ring homomorphism, and twisted maps of bimodules is weak in one
direction. Another example is given by finite sets, Riemann surfaces
with labelled analytically parametrized boundary components,
bijections of finite sets, and holomorphic maps preserving the given
structure. In these two examples we choose the horizontal direction
to be weak, so that bimodules respectively Riemann surfaces are the
horizontal morphisms. In this paper we work only with strict double
categories, though pseudo double categories can also fit into our
framework.

The notion of double category contains many familiar structures. If
we view a category as an internal category in ${\bf Cat}$ with
object and morphism categories discrete, it is equivalent to viewing
an ordinary category as a double category with trivial vertical
morphisms and trivial squares. Every 2-category $\bfC$ can be
considered a double category in at least four ways: as a double
category $\bbH \bfC$ with trivial vertical morphisms, as a double
category $\bbV \bfC$ with trivial horizontal morphisms, as
Ehresmann's quintets $\bbQ\bfC$, or as the transpose of Ehresmann's
quintets $(\bbQ\bfC)^t$. Any double category $\bbD$ has an {\it
underlying horizontal 2-category} $\bfH \bbD$ and an {\it underlying
vertical 2-category} $\bfV \bbD$: we obtain these substructures as
the full sub-double categories with only trivial vertical morphisms
or trivial horizontal morphisms respectively. We denote the
underlying 1-categories of $\bfH \bbD$ and $\bfV \bbD$ by $(\bfH
\bbD)_0$ and $(\bfV \bbD)_0$ respectively. The subscript 0 here
means underlying 1-category of a 2-category, and is unrelated to the
subscript 0 in Definition \ref{defnofdoublecategory}. Though the
formula $(\bfV \bbD)_0=\bbD_0$ holds, $(\bfH \bbD)_0$ is {\it not}
the same as $\bbD_0$.

A {\em double functor} $\xymatrix@1{F\co\bbD \ar[r] & \bbE}$ is an
internal functor in {\bf Cat}. Such a functor consists of functions
$$\xymatrix{\Obj \hspace{1mm} \mathbb{D} \ar[r] & \Obj \hspace{1mm}  \mathbb{E}}$$
$$\xymatrix{\Hor \hspace{1mm} \mathbb{D} \ar[r] & \Hor \hspace{1mm} \mathbb{E}}$$
$$\xymatrix{\Ver \hspace{1mm} \mathbb{D} \ar[r] & \Ver \hspace{1mm} \mathbb{E}}$$
$$\xymatrix{\Sq \hspace{1mm} \mathbb{D} \ar[r] & \Sq \hspace{1mm} \mathbb{E}}$$
which preserve all sources, targets, compositions, and identities.

Internal natural transformations in $\mathbf{Cat}$ are also called
{\it horizontal natural transformations}.

\begin{defn} \label{horizontalnaturaltransformation}
If $\xymatrix@1{F,G\co\mathbb{D} \ar[r] & \mathbb{E}}$ are double
functors, then a {\it horizontal natural transformation}
$\xymatrix@1{\theta\co F \ar@{=>}[r] & G}$ as in
\cite{grandisdouble1} assigns to each object $A$ a horizontal
morhism $\xymatrix@1{\theta A\co FA \ar[r] & GA}$ and assigns to
each vertical morphism $j$ a square
$$\xymatrix@R=3pc@C=3pc{FA \ar[r]^{\theta A} \ar[d]_{Fj} \ar@{}[dr]|{\theta j} & GA \ar[d]^{Gj} \\
FC \ar[r]_{\theta C} & GC}$$ such that:
\begin{enumerate}
\item
For all $A \in \mathbb{D}$, we have $\theta 1^v_A =i^v_{\theta A }$,
\item
For composable vertical morphisms $j$ and $k$,
$$
\begin{array}{c}
\xymatrix@R=3pc@C=3pc{FA \ar[r]^{\theta A } \ar[d]_{F\vcomp{j}{k}}
\ar@{}[dr]|{\theta \vcomp{j}{k} } & GA \ar[d]^{F\vcomp{j}{k}} \\
FE \ar[r]_{\theta E} & GE}
\end{array}
=
\begin{array}{c}
\xymatrix@R=3pc@C=3pc{FA \ar[r]^{\theta A } \ar[d]_{Fj} \ar@{}[dr]|{\theta j } & GA \ar[d]^{Gj} \\
FC \ar[r]|{\lr{\theta C }} \ar[d]_{Fk} \ar@{}[dr]|{\theta k } & GC \ar[d]^{Gk} \\
FE \ar[r]_{\theta E} & GE,}
\end{array}
$$
\item
For all $\alpha$ as in (\ref{sourcetarget}),
$$
\begin{array}{c}
\xymatrix@R=3pc@C=3pc{FA \ar[r]^{Ff} \ar[d]_{Fj}
\ar@{}[dr]|{F\alpha} & FB \ar[r]^{\theta B} \ar[d]|{\tb{Fk}}
\ar@{}[dr]|{\theta k} & GB \ar[d]^{Gk} \\ FC \ar[r]_{Fg} & FC
\ar[r]_{\theta C} & GD}
\end{array}
=
\begin{array}{c}
\xymatrix@R=3pc@C=3pc{FA \ar[r]^{\theta A} \ar[d]_{Fj}
\ar@{}[dr]|{\theta j} & GA \ar[r]^{Gf} \ar[d]|{\tb{Gj}}
\ar@{}[dr]|{G\alpha} & GB \ar[d]^{Gk} \\ FC \ar[r]_{\theta C} & GC
\ar[r]_{Gg} & GD.}
\end{array}
$$
\end{enumerate}
\end{defn}

We also need the analogous notion of vertical natural
transformation.

\begin{defn} \label{verticalnaturaltransformation}
If $\xymatrix@1{F,G\co \mathbb{D} \ar[r] & \mathbb{E}}$ are double
functors, then a {\it vertical natural transformation}
$\xymatrix@1{\sigma\co F \ar@{=>}[r] & G}$ as in
\cite{grandisdouble1} assigns to each object $A$ a vertical morphism
$\xymatrix@1{\sigma A\co FA \ar[r] & GA}$ and assigns to each
horizontal morphism $f$ a square
$$\xymatrix@R=3pc@C=3pc{FA \ar[d]_{\sigma A} \ar[r]^{Ff} \ar@{}[dr]|{\sigma f} & FB  \ar[d]^{\sigma B} \\
GA  \ar[r]_{Gf} & GB}$$ such that:
\begin{enumerate}
\item
For all objects $A \in \mathbb{D}$, we have $\sigma
1^h_A=i^h_{\sigma A}$,
\item
For all composable horizontal morphisms $f$ and $g$, $$\sigma [ f \
g  ]=[\sigma f \ \ \sigma g  ],$$
\item
For all $\alpha$ as in (\ref{sourcetarget}),
$$\vcomp{F \alpha}{\sigma g}=\vcomp{\sigma f}{G\alpha}.$$
\end{enumerate}
\end{defn}

Thus, double categories form a 2-category in two different ways, depending on the choice of
2-cell. Further, there is a useful adjunction with {\bf 2-Cat}.

\begin{prop} \label{HHVV}
Let $\mathbf{DlbCat_h}$ respectively $\mathbf{DblCat_v}$ denote the
2-categories of small double categories, double functors, and
horizontal natural transformations respectively vertical natural
transformations. Let {\bf 2-Cat} denoted the 2-category of small
2-categories, 2-functors, and 2-natural transformations.\footnote{In
this article we follow the convention that 2-functors and 2-natural
transformations are {\it strict} 2-functors and {\it strict}
2-natural transformations.} Then the inclusion 2-functors
$$\xymatrix{\bbH\co \text{\bf 2-Cat} \ar[r] & \mathbf{DblCat_h}}$$
$$\xymatrix{\bbV\co \text{\bf 2-Cat} \ar[r] & \mathbf{DblCat_v}}$$
have as right 2-adjoints the 2-functors
$$\xymatrix{\bfH\co \mathbf{DblCat_h} \ar[r] & \text{\bf 2-Cat}}$$
$$\xymatrix{\bfV\co \mathbf{DblCat_v} \ar[r] & \text{\bf 2-Cat}}$$
respectively. Moreover, the inclusion 2-functors
$$\xymatrix{\bbH\co\mathbf{Cat} \ar[r] & \mathbf{DblCat_h}}$$
$$\xymatrix{\bbV\co\mathbf{Cat} \ar[r] & \mathbf{DblCat_v}}$$
have as right 2-adjoints the 2-functors
$$\xymatrix{(\bfH\text{-})_0\co \mathbf{DblCat_h} \ar[r] & \text{\bf Cat}}$$
$$\xymatrix{(\bfV\text{-})_0\co \mathbf{DblCat_v} \ar[r] & \text{\bf Cat}}$$
respectively.
\end{prop}

\begin{defn} \label{defnofexternalproduct}
If $\mathbf{C}$ and $\mathbf{D}$ are 2-categories, then their {\it
external product} $\mathbf{C}\boxtimes \mathbf{D}$ is the double
category with objects $\Obj \mathbf{C} \times \Obj \mathbf{D}$,
vertical morphisms $$\xymatrix@1{(f,D)\co (C,D) \ar[r] & (C',D)},$$
horizontal morphisms
$$\xymatrix@1{(C,g)\co (C,D) \ar[r] & (C,D')},$$ and squares
$$\xymatrix{(C,D) \ar[r]^{(C,g_1)} \ar[d]_{(f_1,D)} \ar@{}[dr]|{\alpha} & (C,D')
\ar[d]^{(f_2,D')} \\ (C',D) \ar[r]_{(C',g_2)} & (C',D') }$$ given by
pairs $\alpha=(\gamma,\delta)$ of 2-cells $\xymatrix@1{\gamma \co
f_1 \ar@{=>}[r] & f_2}$ and $\xymatrix@1{\delta\co g_1 \ar@{=>}[r] &
g_2}$ in $\mathbf{C}$ and $\mathbf{D}$ respectively.
\end{defn}

We may simplify the foregoing definitions using the operation of
transposition, which interchanges the roles of horizontal and
vertical.

\begin{defn}
The {\it transpose} of a double category $\bbD$ is the double
category $\bbD^t$ with
$$\aligned
\Obj \bbD^t &= \Obj \bbD \\
\Hor \bbD^t &= \Ver \bbD \\
\Ver \bbD^t &= \Hor \bbD \\
\bbD^t\begin{pmatrix} & f &  \\ j & & k
\\ & g & \end{pmatrix} &= \bbD\begin{pmatrix} &  j &  \\ f & & g
\\ & k & \end{pmatrix}
\endaligned$$
and the expected compositions and units. Transposition defines
2-functors
$$\xymatrix{(\text{-})^t\co\mathbf{DblCat_h} \ar[r] &
\mathbf{DblCat_v}}$$
$$\xymatrix{(\text{-})^t\co\mathbf{DblCat_v} \ar[r] &
\mathbf{DblCat_h}}$$ that are mutually inverse.
\end{defn}

\begin{rmk}
A vertical natural transformation $\xymatrix@1{\sigma\co F
\ar@{=>}[r] & G}$ is a horizontal natural transformation
$\xymatrix@1{\sigma^t\co F^t \ar@{=>}[r] & G^t}$. The transpose of
the 2-adjunction $\bbH \dashv \bfH$ is the 2-adjunction $\bbV \dashv
\bfV$. The external product of 2-categories $\bfC$ and $\bfD$ is
$\mathbf{C} \boxtimes
\mathbf{D}=\bbV\bfC\times\bbH\bfD=(\bbH\mathbf{C})^{t} \times
\bbH\mathbf{D}$. More generally, the {\it external product} of
double categories $\bbC$ and $\bbD$ is $\bbC \boxtimes \bbD:=\bbC^t
\times \bbD$.
\end{rmk}

\begin{lem}
The external product of 2-categories is a functor
$$\xymatrix@1{\boxtimes\co \text{{\bf 2-Cat}} \times \text{{\bf 2-Cat}} \ar[r] & \text{{\bf DblCat}}}.$$
\end{lem}
\begin{pf}
Transpose is functorial.
\end{pf}

\begin{examp}
Let $[m]$ denote the partially ordered set $\{0,1,2,\dots,m\}$. Then the double category
$[m] \boxtimes [n]$ has the shape
$$
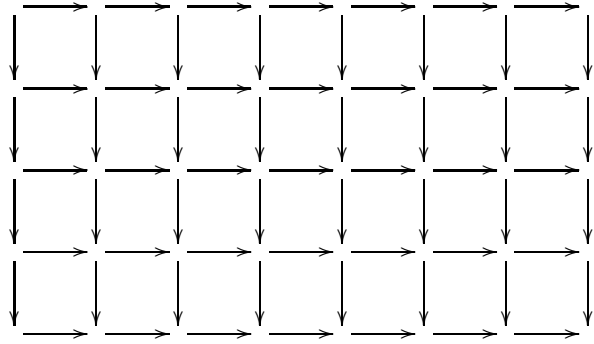
\xymatrix{\ar[r] \ar[d] & \ar[r] \ar[d] & \ar[r] \ar[d]
& \ar[r]\ar[d] & \ar[d] \ar[r] &\ar[d] \ar[r] & \ar[r] \ar[d] &
\ar[d] \\ \ar[r] \ar[d] & \ar[r] \ar[d] & \ar[r] \ar[d] &
\ar[r]\ar[d] & \ar[d] \ar[r] &\ar[d] \ar[r] & \ar[r] \ar[d] &
\ar[d]\\ \ar[r] \ar[d] & \ar[r] \ar[d] & \ar[r] \ar[d] &
\ar[r]\ar[d] & \ar[d] \ar[r] &\ar[d] \ar[r] & \ar[r] \ar[d] &
\ar[d]\\ \ar[r] \ar[d] & \ar[r] \ar[d] & \ar[r] \ar[d] &
\ar[r]\ar[d] & \ar[d] \ar[r] &\ar[d] \ar[r] & \ar[r] \ar[d] & \ar[d]
\\ \ar[r] & \ar[r] & \ar[r] & \ar[r]& \ar[r] & \ar[r] & \ar[r] &}$$
with $m$ rows and $n$ columns of squares.
\end{examp}

We round off this section with a discussion of Cartesian closedness
for $\mathbf{DblCat}$, $\mathbf{DblCat_v}$, and $\mathbf{DblCat_h}$.

\begin{prop}[$n=2$ case of \cite{ehresmannthree}] \label{DblCatCartesianclosednessstatement}
The category $\mathbf{DblCat}$ is Cartesian closed. In other words
for each $\bbD$ there is an endofunctor $(\text{-})^\bbD$ of
$\mathbf{DblCat}$ and a bijection of sets
\begin{equation} \label{DblCatCartesianclosedness}
\mathbf{DblCat}(\bbC \times \bbD, \bbE) \cong
\mathbf{DblCat}\big(\bbC, \bbE^\bbD\big)
\end{equation}
natural in $\bbC$ and $\bbE$.
\end{prop}
\begin{cor}[1.6 of \cite{grandisdouble1}] \label{explicitDblCatCartesianclosedness}
The objects of $\bbE^\bbD$ are double functors $\xymatrix@1{\bbD
\ar[r] & \bbE}$, horizontal morphisms are horizontal natural
transformations, vertical morphisms are vertical natural
transformations, and squares are modifications.
\end{cor}
\begin{pf}
In equation (\ref{DblCatCartesianclosedness}), we take $\bbC$ to be
the terminal double category, $\bbH [1]$, $\bbV [1]$, or $[1]
\boxtimes [1]$. See 1.6 of \cite{grandisdouble1} for the definition
of modification.
\end{pf}

\begin{prop} \label{DblCatvCartesianclosednessstatement}
The 2-category $\mathbf{DblCat_v}$ is Cartesian closed. More
precisely, the functor $(\text{-})^\bbE$ of Proposition
\ref{DblCatCartesianclosednessstatement} and Corollary
\ref{explicitDblCatCartesianclosedness} extends to an endo-2-functor
of $\mathbf{DblCat_v}$, and there is an isomorphism of categories
\begin{equation} \label{DblCatvCartesianclosedness}
\mathbf{DblCat_v}(\bbC \times \bbD, \bbE) \cong
\mathbf{DblCat_v}\big(\bbC, \bbE^\bbD\big)
\end{equation}
2-natural in $\bbC$ and $\bbE$. Similarly, the 2-category
$\mathbf{DblCat_h}$ is Cartesian closed.
\end{prop}
\begin{pf}
Equation (\ref{DblCatCartesianclosedness}) is the object part of the
isomorphism in equation (\ref{DblCatvCartesianclosedness}). For the
bijection of morphism sets we have
$$\aligned
\Mor \mathbf{DblCat_v}(\bbC \times \bbD, \bbE) & \cong
\mathbf{DblCat}(\bbV[1] \times \bbC \times \bbD, \bbE) \\
&\cong \mathbf{DblCat}\big(\bbV[1] \times \bbC, \bbE^\bbD\big) \\
&\cong \Mor \mathbf{DblCat_v}\big(\bbC,\bbE^\bbD\big).
\endaligned$$
The proof that this isomorphism of graphs is a 2-natural functor is
similar to the analogous proof of the Cartesian closedness of
$\mathbf{Cat}$.
\end{pf}

\begin{cor} \label{likeinternalhom}
For a small category $\bfC$ and small double categories $\bbD$ and
$\bbE$, we have an isomorphism of categories
$$
\mathbf{DblCat_v}(\bbV\bfC \times \bbD, \bbE) \cong
\mathbf{Cat}(\bfC, \mathbf{DblCat_v}(\bbD,\bbE))
$$
2-natural in $\bfC$ and $\bbE$.
\end{cor}
\begin{pf}
From equation (\ref{DblCatvCartesianclosedness}) and the
2-adjunction $\bbV \dashv (\bfV\text{-})_0$ of Proposition
\ref{HHVV} we have a 2-natural isomorphism of categories
$$\aligned
\mathbf{DblCat_v}(\bbV\bfC \times \bbD, \bbE) & \cong
\mathbf{DblCat_v}\big(\bbV\bfC, \bbE^\bbD \big) \\
& \cong \mathbf{Cat}\Big(\bfC, \big(\bfV(\bbE^\bbD)\big)_0 \Big) \\
& \cong \mathbf{Cat}(\bfC,\mathbf{DblCat_v}(\bbD,\bbE)).
\endaligned$$
\end{pf}

\section{Free Double Categories and Quotients}
\label{section:FreeDouble} As expected, there is a notion of free
double category and quotient double category. However, the situation
is richer than for ordinary categories, as there is an intermediate
step between double categories and double graphs, which we
call {\it double derivation schemes}. Double derivation schemes and
quotients are crucial in the explicit description of colimits in
Section \ref{section:limitscolimits}, the construction of a left
adjoint to horizontal nerve in Section
\ref{section:categorification}, and the computation of pushouts in
Theorems \ref{pushoutinDblCat} and \ref{pushoutiso}.

In this section we introduce double analogues to some of the
concepts in \cite{streetstructures}. The special kind of double
graphs we will work with have 1-identities but are not equipped with
identity squares. This is important because nontrivial squares in a
double category may very well have one or more trivial edges. Recall
that a {\it reflexive} graph is a graph equipped with a
distinguished {\it identity} edge $\xymatrix@1{1_A \co A \ar[r] &
A}$ for each vertex $A$. All graphs in this paper are directed and
small, so we often leave off these adjectives.

\begin{defn}
A {\it double graph} $\bbA$ is an internal graph in the category of
small graphs. This consists of a set of vertices (objects) $\Obj
\bbA$, a set of horizontal edges $\Hor \bbA$, a set of vertical
edges $\Ver \bbA$, and a set of squares $\Sq \bbA$ equipped with
source and target maps as in (\ref{sourcetarget}).  A {\it morphism
of double graphs} is a morphism of internal graphs in the category
of small graphs, or equivalently, a map which preserves the sources
and targets of (\ref{sourcetarget}). We denote the horizontal and
vertical 2-graphs of a double graph $\bbA$ by $\bfH \bbA$ and $\bfV
\bbA$.
\end{defn}

\begin{defn}
A {\it double graph with 1-identities} is a double graph in which
the horizontal and vertical 1-graphs are reflexive graphs. This
means for each object $A$, there is a distinguished horizontal edge
$\xymatrix@1{1_A^h \co A \ar[r] & A}$ as well as a distinguished
vertical edge $\xymatrix@1{1_A^v \co A \ar[r] & A}$. There are no
distinguished squares. A {\it morphism of double graphs with
1-identities} is a morphism of double graphs which preserves the
distinguished edges. Double graphs with 1-identities form a category
which we denote by {\bf DblGr1-Id}.
\end{defn}

A double graph with 1-identities is a double category without any of the compositions and without identity squares.
The intermediate structure between double graphs with 1-identities and double categories
is analogous to Street's notion of {\it derivation scheme} in \cite{streetstructures}.

\begin{defn}
A {\it double derivation scheme} is a double graph with 1-identities whose vertical
reflexive 1-graph and horizontal reflexive 1-graph are categories. A {\it morphism of
double derivation schemes} is a morphism of double graphs with 1-identities
which is a functor on both the horizontal and vertical 1-categories.
Double derivation schemes form a category which we denote by $\mathbf{DblDerSch}$. We denote
the horizontal and vertical derivation schemes of a double
derivation scheme $\bbS$ by $\bfH \bbS$ and $\bfV \bbS$, and their
underlying categories by $(\bfH \bbS)_0$ and $(\bfV \bbS)_0$.
\end{defn}

To take a free category on a reflexive graph, one merely takes paths
of composable edges and identifies paths which differ only by
insertion or deletion of identity edges. However, the 2-dimensional
situation is more subtle, as evidenced by
\cite{dawsonparefreedouble}, \cite{johnsonpasting},
\cite{power2pasting}, and \cite{powernpasting}. Thus, in the
construction of a free double category we need a careful definition
of {\it allowable compatible arrangement}. We use the notion of
compatible arrangement from \cite{dawsonpareassociativity}, and
develop it further for our purposes.

\begin{defn}
In a double derivation scheme $\bbS$, a {\em compatible arrangement}
consists of a subdivision of a rectangle into smaller rectangles and
a function which assigns to each vertex an object, to each
horizontal line segment a horizontal morphism, to each vertical line
segment a vertical morphism, and to each constituent rectangle a
square in $\bbS$, which are compatible in the sense that
    \begin{enumerate}
     \item
    for each horizontal edge in the subdivision, the
    domain and codomain respectively of the morphism assigned to it are the objects
    assigned to the left and right vertices respectively;
    \item
    for each vertical edge in the subdivision, the
    domain and codomain respectively of the morphism assigned to it are the objects
    assigned to the top and bottom vertices respectively;
    \item
    for each constituent rectangle the composition of the morphisms
    assigned to the edges on
        \begin{enumerate}
        \item
        the left side is the horizontal domain of the square assigned to it;
        \item
        the right side is the horizontal codomain of the square assigned to it;
        \item
        the top is the vertical domain of the square assigned to it;
        \item
        the bottom is the vertical codomain of the square assigned to it.
        \end{enumerate}
    \end{enumerate}
\end{defn}

In the free double category on a double derivation scheme, a square
is a compatible arrangement for which the image under any morphism
of double derivation schemes into a double category becomes
composable to a single square by a sequence of horizontal and
vertical compositions. We will call such compatible arrangements
{\it allowable}. However, an image of a compatible arrangement is
just a compatible arrangement in the target double category with the
same underlying subdivision of the rectangle. So whether a
compatible arrangement is allowable in the free double category
depends only on its shape, {\it i.e.}, the underlying subdivision of
the rectangle.

A {\em horizontal} (resp. {\em vertical}) {\em cut} in a compatible
arrangement is a horizontal (resp. vertical) line segment which
consists of edges of the underlying subdivision of the rectangle.  A
horizontal (respectively vertical) cut is {\em full length} if it
stretches from the left (respectively top) edge of the arrangement
to the right (respectively bottom) edge of the arrangement. We can
use this to characterize when a compatible arrangement is allowable.

\begin{defn}\label{allowable}
A subdivision of a rectangle is {\em allowable} if it is either the
trivial subdivision, consisting of just the rectangle itself, or
contains a full length horizontal or vertical cut which divides it
into two allowable subdivisions. A compatible arrangement is {\it
allowable} if its underlying subdivision of the rectangle is
allowable.
\end{defn}
As an illustration, consider the following two examples of
subdivisions of a rectangle.
$$
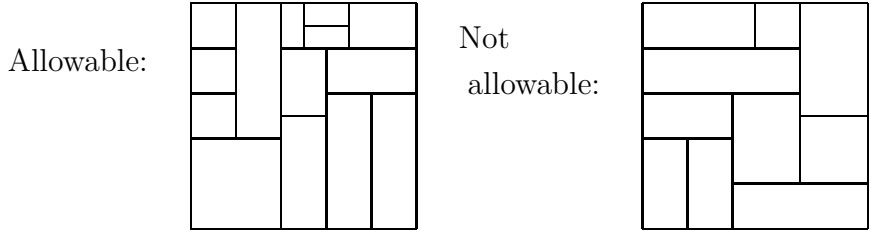
\setlength{\unitlength}{.6mm}
\begin{picture}(175,50)
\put(0,35){\makebox(0,0)[b]{Allowable:}}
\put(25,0){\line(0,1){50}}
\put(25,0){\line(1,0){50}}
\put(25,50){\line(1,0){50}}
\put(75,0){\line(0,1){50}}
\put(25,20){\line(1,0){20}}
\put(25,30){\line(1,0){10}}
\put(25,40){\line(1,0){10}}
\put(25,50){\line(1,0){10}}
\put(35,20){\line(0,1){30}}
\put(45,0){\line(0,1){50}}
\put(45,25){\line(1,0){10}}
\put(45,40){\line(1,0){30}}
\put(50,40){\line(0,1){10}}
\put(50,45){\line(1,0){10}}
\put(55,0){\line(0,1){40}}
\put(55,30){\line(1,0){20}}
\put(60,40){\line(0,1){10}}
\put(65,0){\line(0,1){30}}
\put(90,40){\makebox(0,0)[b]{Not}}
\put(101,30){\makebox(0,0)[b]{allowable:}}
\put(125,0){\line(0,1){50}}
\put(125,0){\line(1,0){50}}
\put(125,50){\line(1,0){50}}
\put(175,0){\line(0,1){50}}
\put(125,20){\line(1,0){20}}
\put(125,30){\line(1,0){35}}
\put(125,40){\line(1,0){35}}
\put(135,0){\line(0,1){20}}
\put(145,0){\line(0,1){30}}
\put(145,10){\line(1,0){30}}
\put(150,40){\line(0,1){10}}
\put(160,10){\line(0,1){40}}
\put(160,25){\line(1,0){15}}
\end{picture}
$$

Note that the notion of {\it allowable} compatible arrangement differs
from the notion of {\it composable} compatible arrangement in
\cite{dawsonpareassociativity} in that Dawson and Par\'e call a
compatible arrangement in a double category $\bbD$ {\it composable} if it
is composable to a single square through the use of both
factorizations and compositions in $\bbD$. So their notion
depends on the ambient double category, not only on the shape of the
arrangement. Any allowable compatible arrangement in our sense is
composable in the sense of Dawson and Par\'e.

\begin{prop}
A compatible arrangement in a double category is allowable if and
only if it can be composed to a single square by a sequence of
horizontal and vertical compositions.
\end{prop}

\begin{pf}
We argue by induction on the number of squares in the arrangement.
The statement is trivially true for arrangements consisting of a single square.
Now let $C$ be a compatible arrangement consisting of two or more squares,
with an assignment into a double category $\bbD$ which is composable by a sequence
of horizontal and vertical compositions of squares.
Consider the last composition used. Without loss of generality, assume that
this is a horizontal composition of two squares $\gamma_1$ and $\gamma_2$
along a vertical morphism $v$, as in
$$
\xymatrix{
\ar[d]\ar[r]\ar@{}[dr]|{\gamma_1} & \ar[d]|{\tb{v}} \ar[r]\ar@{}[dr]|{\gamma_2} & \ar[d]
\\
\ar[r] & \ar[r]&
}
$$
Both $\gamma_1$ and $\gamma_2$ have been obtained by sequences of horizontal and vertical
compositions of squares in $C$, so $v$ is a vertical composition of vertical morphisms
$v_1,\ldots,v_n$ in $C$.
The underlying edges of these vertical morphisms form a cut in the underlying subdivision of
the rectangle for $C$.
The squares on the left side of this cut form
a compatible arrangement, since they form a rectangular subset of a compatible arrangement.
Call this arrangement $C_1$.
It can be composed to $\gamma_1$ by a subsequence of the horizontal vertical compositions
used for $C$. In the same way, the squares on the right side of this cut form a compatible arrangement
$C_2$ which can be composed to $\gamma_2$ by a sequence of horizontal and vertical compositions.
Since both $C_1$ and $C_2$ contain strictly less squares than $C$, the induction hypothesis gives
that they are both allowable compatible arrangements.

Conversely, suppose that a compatible arrangement $C$ of two or more
squares in a double category $\bbD$ is allowable. Then it contains a
horizontal (resp.~vertical) cut into two allowable compatible
arrangements $C_1$ and $C_2$. By induction these arrangements can be
composed to single squares in $\bbD$ by sequences of horizontal and
vertical compositions. Now consider the sequence of horizontal and
vertical compositions used for $C_1$ followed by the one for $C_2$
and then one final vertical (resp. horizontal) composition along the
cut. This shows that $C$ is composable to a single square in $\bbD$
by a sequence of horizontal and vertical compositions of squares.
\end{pf}

For inductive arguments on the number of squares in an allowable
compatible arrangement, we need to know that cutting an allowable
arrangement along {\it any} full length cut produces two smaller
compatible arrangements.

\begin{prop} \label{anycut}
If $CA$ is a compatible arrangement which is allowable, then any
full length cut divides the arrangement into two allowable
compatible arrangements.
\end{prop}
\begin{pf}
We prove this by induction on the number of squares in the
arrangement. It is obviously true for compatible arrangements
consisting of a single square. For an arrangement consisting of
$n\ge 2$ squares, let ${\mathcal C}_1$ be an arbitrary full length
cut as in this proposition and let ${\mathcal C}_2$ be the full
length cut used to establish that $CA$ is allowable. Assume without
loss of generality that ${\mathcal C}_2$ is horizontal. Let $CA_1$
and $CA_2$ be the compatible arrangements obtained by cutting $CA$
along ${\mathcal C}_2$, as in Figure~\ref{horcut}.
\begin{figure}
\setlength{\unitlength}{1mm}
\begin{picture}(35,20)
\put(0,0){\line(1,0){30}} \put(0,0){\line(0,1){20}}
\put(0,10){\line(1,0){30}} \put(0,20){\line(1,0){30}}
\put(30,0){\line(0,1){20}} \put(15,3){\makebox(0,0)[b]{$CA_2$}}
\put(15,13){\makebox(0,0)[b]{$CA_1$}}
\put(33,8){\makebox(0,0)[b]{${\mathcal C}_2$}}
\end{picture}
\label{horcut}
\end{figure}
Note that both of the arrangements $CA_1$ and $CA_2$ are allowable
and contain strictly less than $n$ squares.

If ${\mathcal C}_1$ is vertical, the cut ${\mathcal C}_1$ itself
gets divided by ${\mathcal C}_2$ into two vertical cuts ${\mathcal
C}_{1,1}$ and ${\mathcal C}_{1,2}$, which are full length vertical
cuts for $CA_1$ and $CA_2$ respectively, as in Figure
\ref{dividedcut}.
\begin{figure}
\setlength{\unitlength}{1mm}
\begin{picture}(35,25)
\put(0,0){\line(1,0){30}} \put(0,0){\line(0,1){20}}
\put(0,10){\line(1,0){30}} \put(0,20){\line(1,0){30}}
\put(30,0){\line(0,1){20}} \put(15,0){\line(0,1){20}}
\put(15.5,21){\makebox(0,0)[b]{${\mathcal C}_1$}}
\put(19,13){\makebox(0,0)[b]{${\mathcal C}_{1,1}$}}
\put(19,3){\makebox(0,0)[b]{${\mathcal C}_{1,2}$}}
\put(33,8){\makebox(0,0)[b]{${\mathcal C}_2$}}
\end{picture}
\caption{} \label{dividedcut}
\end{figure}
The cut ${\mathcal C}_{1,1}$ divides $CA_1$ into compatible
arrangements $CA_{1,1}$ and $CA_{1,2}$, and the cut ${\mathcal
C}_{1,2}$ divides  $CA_2$ into compatible arrangements $CA_{2,1}$
and $CA_{2,2}$, as in Figure \ref{subarrangements}.
\begin{figure}
\setlength{\unitlength}{1mm}
\begin{picture}(35,25)
\put(0,0){\line(1,0){30}} \put(0,0){\line(0,1){20}}
\put(0,10){\line(1,0){30}} \put(0,20){\line(1,0){30}}
\put(30,0){\line(0,1){20}} \put(15,0){\line(0,1){20}}
\put(15,22){\makebox(0,0)[b]{${\mathcal C}_1$}}
\put(7.5,3){\makebox(0,0)[b]{$CA_{2,1}$}}
\put(22.5,3){\makebox(0,0)[b]{$CA_{2,2}$}}
\put(7.5,13){\makebox(0,0)[b]{$CA_{1,1}$}}
\put(22.5,13){\makebox(0,0)[b]{$CA_{1,2}$}}
\put(33.5,8){\makebox(0,0)[b]{${\mathcal C}_2$}}
\end{picture}
\caption{} \label{subarrangements}
\end{figure}
By the induction hypothesis, $CA_{1,1}$, $CA_{1,2}$, $CA_{2,1}$, and
$CA_{2,2}$ are all allowable. It is clear that the compatible
arrangement to the left of ${\mathcal C}_1$ gets divided into
$CA_{1,1}$ and $CA_{2,1}$ by the left side of the cut ${\mathcal
C}_2$, so the compatible arrangement to the left of ${\mathcal C}_1$
is allowable. In the same way the compatible arrangement to the
right of ${\mathcal C}_1$ gets cut into $CA_{1,2}$ and $CA_{2,2}$ by
the right side of the cut ${\mathcal C}_2$, so this compatible
arrangement is also allowable, as we wanted to prove.

If ${\mathcal C}_1$ is horizontal, assume without loss of generality
that $CA_1$  contains ${\mathcal C}_1$. By the induction hypothesis,
${\mathcal C}_1$ divides the allowable compatible arrangement $CA_1$
into two allowable compatible arrangements, say $CA_{1,a}$ and
$CA_{1,b}$, as in Figure \ref{newsubarrangements}.
\begin{figure}
\setlength{\unitlength}{1mm}
\begin{picture}(35,30)
\put(0,0){\line(1,0){30}} \put(0,0){\line(0,1){30}}
\put(0,12){\line(1,0){30}} \put(0,21){\line(1,0){30}}
\put(0,30){\line(1,0){30}} \put(30,0){\line(0,1){30}}
\put(15,3){\makebox(0,0)[b]{$CA_2$}}
\put(15,15){\makebox(0,0)[b]{$CA_{1,b}$}}
\put(15,24){\makebox(0,0)[b]{$CA_{1,a}$}}
\put(33,19){\makebox(0,0)[b]{${\mathcal C}_1$}}
\put(33,10){\makebox(0,0)[b]{${\mathcal C}_2$}}
\end{picture}
\caption{} \label{newsubarrangements}
\end{figure}
We derive that ${\mathcal C}_1$ divides the total arrangement $CA$
into two compatible arrangements, $CA_{1,a}$ and $CA_{1,c}$, the
latter of which is divided by ${\mathcal C}_2$ into $CA_{1,b}$ and
$CA_2$. Since both $CA_{1,b}$ and $CA_2$ are allowable, we conclude
that both $CA_{1,c}$ and $CA_{1,a}$ are allowable. This completes
the proof.
\end{pf}

\begin{prop} \label{freeconstructions}
The forgetful functors $T$ and $U$ admit left adjoints $S$ and $R$.
$$\xymatrix@C=4pc{\text{\bf DblGr1-Id} \ar@{}[r]|{\perp} \ar@/^1pc/[r]^{S}
& \ar@/^1pc/[l]^{T} {\bf DblDerSch}  \ar@/^1pc/[r]^{R}
\ar@{}[r]|{\perp} & {\bf DblCat} \ar@/^1pc/[l]^{U} }$$ The left
adjoint $S$ gives the {\it free  double derivation scheme on a
double graph with 1-identities}, and the left adjoint $R$ gives the
{\it free double category on a double derivation scheme}. The
functor $R$ preserves the horizontal and vertical 1-categories.
\end{prop}
\begin{pf}
For a double graph with 1-identities $\bbA$, let $S \bbA$ have vertical and
horizontal 1-categories the free 1-categories on the respective
reflexive graphs. The set of squares remains the same. It is
straightforward to verify that this defines a left adjoint to $T$.

For a double derivation scheme $\bbS$, let $R\bbS$ have vertical and
horizontal 1-categories the vertical and horizontal 1-categories of
$\bbS$ respectively. The nonidentity squares of $R\bbS$ are allowable compatible
arrangements of squares of $\bbS$.  Such compatible arrangements are
composed vertically and horizontally by concatenation. Clearly,
composites of allowable compatible arrangements are allowable. We additionally add horizontal
and vertical identity squares.

If $\xymatrix@1{J\co \bbS \ar[r] & U\bbD}$ is a morphism of double
derivation schemes, then it induces a double functor
$\xymatrix@1{J'\co R\bbS \ar[r] & \bbD}$ which is $J$ on the
horizontal and vertical 1-categories. For an allowable compatible
arrangement $D$, $J'D$ is the composite in $\bbD$ of $J$ applied to
the constituents of $D$. Morphisms $\xymatrix@1{R\bbS \ar[r] &
\bbD}$ restrict to morphisms $\xymatrix@1{\bbS \ar[r] & U\bbD}$, and
it is not hard to check that these two operations are inverse. We
conclude that $R\dashv U$.
\end{pf}

Now that we have free notions, we also define quotients. Note that
the notion of congruence for ordinary categories is an equivalence
relation on the cells of highest dimension, satisfying certain
compatibility properties. We imitate this in our notion of
congruence for a double category.

\begin{defn}
A {\it congruence} on a category {\bf C} is an equivalence relation
on $\bfC(a,b)$ for each $a,b \in \bfC$, such that if $f \sim f'$ and
$g \sim g'$, then $gf \sim g'f'$ whenever the composites exist.
\end{defn}

\begin{defn}
A {\em congruence} on a double derivation scheme $\bbS$ consists of
a congruence on the horizontal 1-category and a congruence on the
vertical 1-category.
\end{defn}

\begin{defn} \label{doublecategorycongruence}
A {\em congruence} on a double category $\bbD$ consists of an
equivalence relation on $\begin{array}{c}\bbD\begin{pmatrix} & f &  \\
j & & k \\ & g & \end{pmatrix}\end{array}$ for each boundary
$\begin{array}{c}\xymatrix{A \ar[r]^f \ar[d]_j & B \ar[d]^k \\ C
\ar[r]_g & D}\end{array}$ such that if $\alpha \sim \alpha', \beta
\sim \beta',$ and $\gamma \sim \gamma'$ then $$\begin{bmatrix}
\alpha & \beta
\end{bmatrix}\sim\begin{bmatrix} \alpha' & \beta' \end{bmatrix}$$
$$\begin{bmatrix} \alpha \\ \gamma
\end{bmatrix}\sim\begin{bmatrix} \alpha' \\ \gamma' \end{bmatrix}$$
whenever the composites exist. Note that the congruence does not concern the
horizontal and vertical morphisms.
\end{defn}

\begin{examp} \label{congruencegeneratedbyK}
Suppose $\xymatrix@1{K\co \bbD \ar[r] & \bbE}$ is a double functor. Then we may define a congruence $\sim_K$ on $\bbD$ by
$$\alpha \sim_K \alpha' :\Longleftrightarrow K(\alpha)=K(\alpha').$$
\end{examp}

\begin{prop}
Let $\bbD$ be a double category equipped with a congruence. If two
allowable compatible arrangements $D_1$ and $D_2$ with the same
underlying subdivision of the rectangle have congruent constituent
squares, then the composites of $D_1$ and $D_2$ in $\bbD$ are
congruent.
\end{prop}

\begin{pf}
By Theorem 1.2 of \cite{dawsonpareassociativity}, any two composites
of a composable compatible arrangement are equal. The compatible
arrangements $D_1$ and $D_2$ are composable since they are
allowable. If we compose each of $D_1$ and $D_2$ using the same
sequence of pairwise compositions, then the pairwise composites in
each step are congruent. An inductive argument shows that total
composites are then also congruent.
\end{pf}

\begin{defn} \label{quotientcategorybycongruence}
Let $\bfC$ be a category and $\sim$ a congruence on $\bfC$. The {\em
quotient category $\bfC/\negthinspace\sim$} has the same objects as
{\bf C} and has homsets
$(\bfC/\negthinspace\sim)(a,b)=\bfC(a,b)/\negthinspace\sim$. The
composition in $\bfC/\negthinspace\sim$ is induced by the
composition in $\bfC$.
\end{defn}

\begin{defn}
Let $\bbS$ be a double derivation scheme and $\sim$ a congruence on
$\bbS$. The {\em quotient double derivation scheme
$\bbS/\negthinspace\sim$} has the same objects and squares as
$\bbS$. The horizontal and vertical 1-categories of
$\bbS/\negthinspace\sim$ are the quotient categories of $(\bfH
\bbS)_0$ and $(\bfV \bbS)_0$.
\end{defn}

\begin{defn} \label{quotientdoublecategorybycongruence}
Let $\bbD$ be a double category and $\sim$ a congruence on $\bbD$.
The {\em quotient double category $\bbD/\negthinspace\sim$} has the
same objects and the same horizontal and vertical 1-categories as
$\bbD$. The set of squares of $\bbD/\negthinspace\sim$ with the
indicated boundary are
$$\begin{array}{c}(\bbD/\negthinspace\sim)\begin{pmatrix} & f &  \\
j & & k \\ & g & \end{pmatrix}\end{array}=\begin{array}{c} \bbD\begin{pmatrix} & f &  \\
j & & k \\ & g & \end{pmatrix}\end{array}\negthinspace/ \sim.$$ The
horizontal and vertical compositions of squares in
$\bbD/\negthinspace\sim$ are induced by the horizontal and vertical
compositions of squares in $\bbD$.
\end{defn}

These are of course not the most general notions of quotient, but
more general quotients can be built from these as follows. All
quotients can be characterized by the usual universal properties.

\begin{defn} \label{quotientcategorybyequivalence}
Let $\bfC$ be a category and $\bfR \subseteq \bfC \times \bfC$ a
subcategory satisfying the usual axioms of an equivalence relation
both on the set of objects and on the set of morphisms. Then the
{\em quotient category $\bfC/\bfR$} is defined as follows. First we
obtain a graph with object set $\Obj \bfC / \Obj\bfR$ and morphism
set $\Mor  \bfC / \Mor \bfR$. We make this into a reflexive graph by
identifying $1_A$ and $1_B$ whenever $A$ and $B$ are identified. Let
$\bfF$ be the free category on this reflexive graph. The {\em
quotient category} $\bfC/\bfR$ is defined as $\bfF / \sim$ where
$\sim$ is the smallest congruence on the free category $\bfF$ such
that induced map of reflexive graphs $\xymatrix@1{\bfC \ar[r] &
\bfF/\sim}$ is a functor.
\end{defn}

Such quotients of categories have been considered in
\cite{boergerdiplomarbeit}. However a counterexample in
\cite{boergerdiplomarbeit}, \cite{merschquotients1}, and
\cite{merschquotients2} shows that the quotient functor may identify
morphisms which are not equivalent. Early work on quotients is found
in \cite{isbellquotients}. More recently, quotients of categories by
generalized congruences have been considered in \cite{bednarczyk}.

For general quotients of double categories, we need quotients of
double derivation schemes as an intermediate notion.

\begin{defn} \label{quotientschemebyequivalence}
Let $\bbS$ be a double derivation scheme and $\bbR \subseteq \bbS
\times \bbS$ a sub-double derivation scheme satisfying the usual
axioms of an equivalence relation on the sets of objects, vertical
morphisms, horizontal morphisms, and squares. Then the {\em quotient
double derivation scheme $\bbS/\bbR$} is defined as follows. The
horizontal and vertical 1-categories are the quotients of the
horizontal and vertical 1-categories of $\bbS$ as in Definition
\ref{quotientcategorybyequivalence}. The squares are
$\Sq(\bbS/\bbR)=(\Sq \bbS )/(\Sq\bbR).$
\end{defn}

\begin{defn}
Let $\bbD$ be a double category and $\bbR \subseteq \bbD \times
\bbD$ a sub-double category satisfying the usual axioms of an
equivalence relation on the sets of objects, vertical morphisms,
horizontal morphisms, and squares. Then the {\em quotient double
category $\bbD/\bbR$} is defined as follows. First we take the
quotient of the underlying double derivation scheme of $\bbD$ by the
underlying double derivation scheme of $\bbR$ as in Definition
\ref{quotientschemebyequivalence}. Let $\bbF$ be the free double
category $\bbF$ on this double derivation scheme. The {\em quotient
double category $\bbD/\bbR$} is defined as $\bbF / \sim$ where
$\sim$ is the smallest congruence on the free double category $\bbF$
such that the induced morphism of double derivation schemes
$\xymatrix@1{\bbD \ar[r] & \bbF /\sim }$ is a double functor.
\end{defn}

For the above definition, see Definitions
\ref{doublecategorycongruence} and
\ref{quotientdoublecategorybycongruence} for quotients of double
categories by congruences. Note that only squares get identified in
the last step, since the horizontal and vertical 1-categories of the
free double category on a double derivation scheme are the same as
the horizontal and vertical 1-categories of the double derivation
scheme.

We will make use of free double categories and their quotients in
our discussion of categorification in Section
\ref{section:categorification} as well as in an explicit description
of certain pushouts of double categories in Theorem
\ref{pushoutinDblCat} and Theorem \ref{pushoutiso}. These are
essential ingredients in the construction of model structures on
{\bf DblCat}. For now it is sufficient to give a colimit formula in
{\bf DblCat}.

\section{Limits and Colimits of Double Categories} \label{section:limitscolimits}

Model structures in general require the existence of limits and
colimits. Moreover, in order to transfer model structures along
certain adjunctions we will need an explicit formula for certain
pushouts of double categories, as in Theorem \ref{pushoutinDblCat}
and Theorem \ref{pushoutiso}. So in this section we discuss limits
and colimits of double categories. We also prove that the
2-categories $\mathbf{DblCat_v}$ and $\mathbf{DblCat_h}$ are
2-cocomplete.

Colimits for categories were described in detail by Gabriel and
Zisman. Their construction was extended in \cite{worytkiewicz2Cat}
to a construction of colimits in {\bf 2-Cat}. We extend this further
to a construction in $\mathbf{DblCat}$ which goes roughly as
follows. To take the colimit of a functor $F$ from an indexing
category $I$ into $\mathbf{DblCat}$, first we take the colimit
$\bbS$ of the underlying double derivation schemes, then we take the
free double category $\bbF$ on $\bbS$, and finally we form the
quotient $\bbF /\sim$ of the free double category $\bbF$ by the
smallest congruence $\sim$ on $\bbF$ such that the induced maps of
double derivation schemes $\xymatrix@1{Fi \ar[r] & \bbF/\sim}$ are
double functors. The quotient double category $\bbF/ \sim$ is the
colimit of the functor $F$. The intermediate notion of double
derivation scheme allows us to deal with the quotients of morphisms
and quotients of squares separately. We present the details in the
following theorems.

\begin{thm} \label{complete} \label{cocomplete}
The category {\bf DblCat} is complete and cocomplete.
\end{thm}
\begin{pf}
The limits of the sets of objects, horizontal morphisms, vertical
morphisms, and squares assemble to form a double category and this
double category is the limit. After all, {\bf DblCat} is a category of algebras.

The category {\bf DblCat} is the category of models in ${\bf Cat}$
of a sketch with finite diagrams, and ${\bf Cat}$ is locally
finitely presentable, so an application of Proposition 1.53 of
\cite{adamekrosicky1994} shows that ${\bf DblCat}$ is locally
finitely presentable. Locally finitely presentable categories are
cocomplete, so {\bf DblCat} is cocomplete.
\end{pf}

Note that the underlying horizontal
and vertical 2-categories of the limit are the limits of the
underlying horizontal and vertical 2-categories, since $\bfH$ and
$\bfV$ admit left adjoints by Proposition \ref{HHVV}. The forgetful
functor $\xymatrix@1{\textbf{2-Cat} \ar[r] & \mathbf{Cat}}$ also
admits a left adjoint, so similar comments hold for the horizontal
and vertical 1-categories.

\begin{thm} \label{2cocomplete}
The 2-categories $\mathbf{DblCat_v}$ and $\mathbf{DblCat_h}$ are
2-co\-complete.
\end{thm}
\begin{pf}
We prove that $\mathbf{DblCat_v}$ is 2-cocomplete; the proof for
$\mathbf{DblCat_h}$ is completely analogous.

The {\it cotensor product} $\{\bfC, \bbE\}$ of a category $\bfC$
with a double category $\bbE$ is $\bbE^{\bbV\bfC}$, since
$$\aligned \mathbf{DblCat_v}\big(\bbD,\bbE^{\bbV\bfC}\big) & \cong
\mathbf{DblCat_v}(\bbD\times\bbV\bfC,\bbE) \\
& \cong \mathbf{DblCat_v}(\bbV\bfC \times \bbD,\bbE) \\
& \cong \mathbf{Cat}(\bfC, \mathbf{DblCat_v}(\bbD,\bbE))
\endaligned$$
is an isomorphism of categories 2-natural in $\bbD$ by Proposition
\ref{DblCatvCartesianclosednessstatement} and Corollary
\ref{likeinternalhom}. Thus $\mathbf{DblCat_v}$ is {\it cotensored},
and by the dual of a statement on page 50 of \cite{kelly3}, the
existence of conical colimits in the 2-category $\mathbf{DblCat_v}$
is equivalent to the existence of ordinary conical colimits in its
underlying 1-category $\mathbf{DblCat}$. But ordinary conical
colimits exist in $\mathbf{DblCat}$ by Theorem \ref{cocomplete}, so
that the 2-category $\mathbf{DblCat_v}$ admits conical colimits.

The {\it tensor product} $\bfC \ast \bbE$ of a category $\bfC$ with
a double category $\bbD$ is $\bbV\bfC \times \bbD$, since
$$
\mathbf{DblCat_v}(\bbV\bfC \times \bbD,\bbE) \cong
\mathbf{Cat}(\bfC, \mathbf{DblCat_v}(\bbD,\bbE))
$$
is an isomorphism of categories 2-natural in $\bbE$ by Corollary
\ref{likeinternalhom}. Thus $\mathbf{DblCat_v}$ is {\it tensored}.

Since $\mathbf{DblCat_v}$ admits conical colimits and tensor
products, we conclude from the dual of Theorem 3.73 in \cite{kelly3}
that $\mathbf{DblCat_v}$ is 2-co\-complete.
\end{pf}

We work towards an explicit description of colimits in
$\mathbf{DblCat}$ which mimics Gabriel and Zisman's calculation of
colimits in {\bf Cat} below.
\begin{thm}[Colimit Formula in {\bf Cat} of \cite{gabrielzismanfractions}] \label{catcolimit}
The colimit of a functor $\xymatrix@1{F\co I \ar[r] &
{\bf Cat}}$ is calculated as follows. Let $\bfF$ be the free category on
the colimit of the underlying reflexive graphs. The colimit of $F$ is
the quotient $\bfF /\sim$ of the free category
$\bfF$ by the smallest congruence $\sim$ on $\bfF$ such that the
induced morphisms of reflexive graphs $\xymatrix@1{Fi \ar[r] &
\mathbf{F}/\sim}$ are functors.
\end{thm}

\begin{lem} \label{horvertcolimit}
The horizontal and vertical 1-categories of a colimit of double
derivation schemes are the colimits of the underlying horizontal and
vertical 1-categories. Similarly, the horizontal and vertical
1-categories of a colimit of double categories are the colimits of
the underlying horizontal and vertical 1-categories.
\end{lem}
\begin{pf}
The right adjoint to the forgetful functor
$$\xymatrix{\mathbf{DblDerSch} \ar[r] & \mathbf{Cat}}$$
$$\xymatrix{\bbS \ar@{|->}[r] & (\bfH \bbS)_0}$$ assigns to a
category $\bfE$ the double derivation scheme $\bbE$ with horizontal
1-category $\bfE$, a unique vertical morphism between any two
objects, and a unique square for each boundary. Similarly, the
forgetful functor $\xymatrix@1{\bbS\, \ar@{|->}[r] & (\bfV \bbS)_0}$
admits a right adjoint. Since left adjoints preserve colimits, the
statement for double derivation schemes follows.

The same argument works for {\bf DblCat} in place of {\bf
DblDerSch}.
\end{pf}

\begin{thm}[Colimit Formula in {\bf DblDerSch}] \label{schemecolimit}
The colimit $\bbS$ of a functor $\xymatrix@1{F\co I \ar[r] & {\bf
DblDerSch}}$ is calculated in the following way. Let $\bbF$ be the
free double derivation scheme on the colimit of the underlying
double graphs with 1-identities. The colimit $\bbS$ of $F$ is the
quotient $\bbF / \sim$ of the free double derivation scheme $\bbF$
by the smallest congruence $\sim$ on $\bbF$ such that the induced
morphisms of double graphs with 1-identities $\xymatrix@1{Fi \ar[r]
& \bbF/ \sim}$ are morphisms of double derivation schemes.
\end{thm}
\begin{pf}
Suppose $\bbS'$ is a double derivation scheme and
$\xymatrix@1{\beta_i\co Fi \ar[r] & \bbS'}$ are natural morphisms of
double derivations schemes. We define a unique factorization
$$\xymatrix{Fi \ar[r]^{\beta_i} \ar[d]_{\alpha_i} & \bbS' \\ \bbS \ar@{.>}[ur] &  }$$
on horizontal and vertical 1-categories by the universal property of
Lemma \ref{horvertcolimit}, and on squares by the universal property
of the colimit of the sets $\Sq Fi$. The set of squares in the free
double derivation scheme on a double graph with identities is the
same as the set of squares in the double graph with identities by
Proposition \ref{freeconstructions}.
\end{pf}

\begin{thm}[Colimit Formula in {\bf DblCat}] \label{colimitformulaDblCat}
The colimit $\bbC$ of a functor
$$\xymatrix@1{F\co I \ar[r] & \mathbf{DblCat}}$$ is calculated as
follows. Let $\bbS$ be the colimit in $\mathbf{DblDerSch}$ of the
underlying double derivation schemes, and $\bbF$ the free double
category on $\bbS$. The colimit $\bbC$ of $F$ is the quotient $\bbF/\sim$ of
$\bbF$ by the smallest congruence $\sim$ such that the induced
natural morphisms of double derivation schemes
\begin{equation}  \label{naturalmorphismsddschemes}
\xymatrix{Fi \ar[r]^-{\alpha_i} & \bbS \ar[r]^p & \ar[r] \bbF & \bbF/\sim}
\end{equation}
 are
double functors. Note that the horizontal and vertical 1-categories
of $\bbS, \bbF$, and $\bbC$ are the same, in particular the
horizontal and vertical 1-categories of $\bbC$ are the colimits of
the horizontal and vertical 1-categories of the $Fi$.
\end{thm}
\begin{pf} Let $\xymatrix@1{q\co \bbS \ar[r] & \bbC }$ denote the morphism of double derivation schemes defined as
the composite of the inclusion $p$ with the quotient double functor
from $\bbF$ to $\bbC$.
Then $q \circ \alpha_i$ is a double functor for all $i \in I$.
Suppose $\bbC'$ is a double category and $\xymatrix@1{\beta_i\co Fi
\ar[r] & \bbC'}$ are natural double functors. Then by Lemma
\ref{schemecolimit} there exists a unique morphism $J$ of double
derivation schemes that makes the upper left triangle commute,
$$\xymatrix@R=3pc@C=3pc{Fi \ar[r]^{\beta_i} \ar[d]_{\alpha_i} &
\bbC'
\\ \bbS \ar@{.>}[ur]|{{}_{\phantom{X}}\exists ! \, J^{\phantom{X}}} \ar[r]_q
& \bbC \rlap{\,.}\ar@{.>}[u]_{\exists !\, L} }$$ The morphism $J$
induces a double functor $\xymatrix@1{K\co \bbF \ar[r] & \bbC'}$
since $\bbF$ is free on $\bbS$. Since $K \circ p \circ \alpha_i = \beta_i$ is a double
functor for all $i$, the induced morphisms of double derivations schemes
\begin{equation*}
\xymatrix@1{Fi \ar[r] & \bbF/ \sim_K}
\end{equation*}
analogous to (\ref{naturalmorphismsddschemes}) are double functors. Here
$\sim_K$ is defined as in Example \ref{congruencegeneratedbyK}. Since $K$ preserves $\sim_K$ and $\sim_K$ contains
$\sim$, the double functor $K$ also preserves $\sim$ and induces a unique functor $L$ which makes the lower right
triangle commute. Therefore the square commutes, and further $L$ is
the unique double functor such that the square commutes by the
uniqueness of the two fillers.
\end{pf}

Recall that filtered colimits in {\bf Cat} are particularly simple
to calculate: the filtered colimit of the underlying reflexive
graphs is already a category and this category is the filtered
colimit in {\bf Cat}. Similarly, one does not need to use free
constructions and quotients to calculate filtered colimits in {\bf
DblCat}.
\begin{thm} \label{filteredcolimits}
A filtered colimit of double categories is calculated by simply
taking the filtered colimits of the underlying reflexive double
graphs.
\end{thm}
\begin{pf}
The filtered colimit of the underlying reflexive double graphs admits all the
associative and unital compositions necessary for a double category
by the corresponding result in {\bf Cat}. The interchange law holds
because it is possible to find representatives of all four squares
in a single stage, where the interchange law is known to hold.
\end{pf}

\section{Nerves of Double Categories} \label{section:Nerves}
Grothendieck's full and faithful nerve $\xymatrix@1{N\co
\mathbf{Cat} \ar[r] & \mathbf{SSet}}$ has been of tremendous use in
higher category theory. One can expect that its $n$-fold version
will similarly be of use. In fact, the authors of
\cite{brownhigginsgroupoidscrossedcomplexes},
\cite{brownhigginsgroupoidscubicalTcomplexes},
\cite{brownhigginscubes}, and \cite{brownhigginstensor} have studied
edge symmetric $n$-fold categories from the point of view of cubical
sets. A double category is a 2-truncated cubical set. We introduce
in this section {\it simplicial} and {\it bisimplicial} nerves of
double categories. The simplicial nerve will be of use in Section
\ref{section:simplicialcategories} where we transfer model
structures on $\mathbf{Cat}^{\Delta^{op}}$ to {\bf DblCat} via a
horizontal categorification-horizontal nerve adjunction. The
bisimplicial nerve will be used in a future article to transfer a
model structure from bisimplicial sets to {\bf DblCat}. The first
nerve we consider is the horizontal nerve, which is really an
internal notion.

\begin{defn}
Let $\bbD=(\bbD_0,\bbD_1)$ be a double category. Then the {\it
horizontal nerve of } $\bbD$ is the simplicial object $N_h\bbD$ in
{\bf Cat} defined as
$$(N_h\bbD)_0=\bbD_0$$
$$(N_h\bbD)_1=\bbD_1$$
$$(N_h\bbD)_n=\underbrace{\bbD_1 \; {}_t\!\times_s \bbD_1 \; {}_t\!\times_s \cdots {}_t\!\times_s \bbD_1}_{
\text{$n$ copies of $\bbD_1$}}.$$
\end{defn}

$$\Obj(N_h\bbD)_n:\hspace{0in} \xymatrix{\ar[r] & \ar[r] & \ar[r] & \ar[r]& \ar[r] & \ar[r] & \ar[r] & } $$
$$\Mor(N_h\bbD)_n:\hspace{0in} \begin{array}{c} \xymatrix{\ar[r] \ar[d] & \ar[r] \ar[d] & \ar[r] \ar[d]
& \ar[r]\ar[d] & \ar[d] \ar[r] &\ar[d] \ar[r] & \ar[r] \ar[d] & \ar[d] \\
\ar[r] & \ar[r] & \ar[r] & \ar[r]& \ar[r] & \ar[r] & \ar[r] &}\end{array}$$
In other words,
$$\Obj (N_h\mathbb{D})_n=\text{\it Mor}_{{\bf Cat}}([n], (\Obj \mathbb{D},\Hor \mathbb{D}))$$
$$\Mor (N_h\mathbb{D})_n=\text{\it Mor}_{{\bf Cat}}([n], (\Ver \mathbb{D},\Sq \mathbb{D})).$$ Composition in
$(N_h\bbD)_n$ is vertical.

There is of course the analogous notion of {\it vertical nerve of}
$\bbD$ denoted $N_v\bbD$.

\begin{examp}
If $\bfC$ is a category, then the simplicial set $N_h(\bbH \bfC)$ is the usual nerve of $\bfC$
viewed as simplicial object in $\mathbf{Cat}$ with discrete categories at each level. In other words, if
we consider the category $\bfC$ as a double category with discrete object category and discrete morphism category,
then the horizontal nerve of $\bfC$ is the classical nerve of $\bfC$. This is the reason we prefer to use
the horizontal nerve rather than the vertical nerve in this paper.
\end{examp}

Like the nerve of a category, the horizontal nerve of a double category
has a representable definition. Recall that $\mathbf{DblCat_v}$ denotes the 2-category of small double
categories, double functors, and vertical natural transformations.

\begin{prop} \label{horizontalrepresentabledefn}
For every double category $\bbD$, the simplicial category
$$[n] \mapsto \mathbf{DblCat_v}(\bbH[n],\bbD)$$ is isomorphic to the
horizontal nerve $N_h\bbD$. Equivalently, the object simplicial set
of the horizontal nerve is
$$[n] \mapsto \mathbf{DblCat}([0] \boxtimes [n],\bbD)$$ and the morphism
simplicial set of the horizontal nerve is
$$[n] \mapsto \mathbf{DblCat}([1] \boxtimes [n],\bbD).$$
\end{prop}
\begin{pf}
The double categories $\bbH[n]$ and $[0] \boxtimes [n]$ are
isomorphic, and vertical natural transformations
$\xymatrix@1{\bbH[n] \ar[r] & \bbD}$ are the same as double functors
$\xymatrix@1{[1] \boxtimes [n] =(\bbH[1])^t \times \bbH[n] \ar[r] &
\bbD}$ as pointed out in \cite{grandisdouble1}.
\end{pf}

Proposition \ref{horizontalrepresentabledefn} makes the
functoriality of $N_h$ immediate. Even more, if we make
$\mathbf{Cat}^{\Delta^{op}}$ into a 2-category with 2-cells the
modifications, then $N_h$ becomes a 2-functor.

\begin{cor}
The horizontal nerve is a 2-functor $$\xymatrix@1{N_h\co
\mathbf{DblCat_v} \ar[r] & \mathbf{Cat}^{\Delta^{op}}}.$$
\end{cor}

In Section \ref{section:categorification} we construct the left
adjoint to the horizontal nerve explicitly, but for now we observe
that a left adjoint exists. Recall the Enriched Lemma from Kan,
which follows from Theorem 4.51 of \cite{kelly3}.

\begin{thm}[Enriched Lemma from Kan] \label{Kansingular}
Let $\mathcal{V}$ be a symmetric monoidal closed category with small
homsets. Suppose $\bfA$ is a small $\mathcal{V}$-category, $\bfB$ is
a cocomplete $\mathcal{V}$-category, and $\xymatrix@1{J\co \bfA
\ar[r] & \bfB}$ a $\mathcal{V}$-functor. Then the enriched left Kan
extension of $J$ along the Yoneda embedding exists and is the
enriched left adjoint of the singular functor
$$\xymatrix@1{J_\ast\co \bfB \ar[r] & \mathcal{V}^{\bfA^{op}}}$$
$$\xymatrix{B \ar@{|->}[r] & \text{\it Mor}_{\bfB}(J(-),B).}$$
\end{thm}

\begin{thm} \label{horizontalcategorificationexists}
The horizontal nerve $\xymatrix@1{N_h\co \mathbf{DblCat} \ar[r] &
\mathbf{Cat}^{\Delta^{op}}}$ admits a left 2-adjoint $c_h$ called
{\it horizontal categorification}.
\end{thm}
\begin{pf}
Let $\mathcal{V}$ be $\mathbf{Cat}$, and let $\xymatrix@1{J\co \bfA
\ar[r] & \bfB}$ be the $\mathbf{Cat}$-functor
$$\xymatrix{\bbH\co\Delta \ar[r] & \mathbf{DblCat_v}.}$$
By Proposition \ref{horizontalrepresentabledefn} the horizontal
nerve is $J_\ast$.

Since the 2-category $\mathbf{DblCat_v}$ is 2-cocomplete by Theorem
\ref{2cocomplete}, and $\Delta$ is small, we may now apply Theorem
\ref{Kansingular} to obtain the left 2-adjoint $c_h$.
\end{pf}

\begin{thm}\label{horizontalnervepreservesfilteredcolimits}
The horizontal nerve $N_h$ preserves filtered colimits.
\end{thm}
\begin{pf}
It follows from Theorem \ref{filteredcolimits} that the category of
horizontal morphisms and squares of a filtered colimit of double
categories is the filtered colimit of the categories of horizontal
morphisms and squares. Since filtered colimits commute with finite
limits, in particular iterated pullbacks, $N_h$ preserves filtered
colimits.
\end{pf}

The horizontal nerve is also well behaved with respect to external products.

\begin{prop} \label{horizontalnerveexternalproduct}
Let $\xymatrix@1{\sigma\co \mathbf{Cat} \ar[r] &
\mathbf{Cat}^{\Delta^{op}}}$ denote the constant functor. Let
$\xymatrix@1{\nu\co \mathbf{Set}^{\Delta^{op}} \ar[r] &
\mathbf{Cat}^{\Delta^{op}}}$ be the inclusion induced by the functor
$\xymatrix@1{\mathbf{Set} \ar[r] & \mathbf{Cat}}$ which takes a set
to the corresponding discrete category. If {\bf A} and {\bf B} are
categories, then $N_h(\mathbf{A} \boxtimes \mathbf{B}) = \sigma
\mathbf{A} \times \nu N\mathbf{B}$. In other words, $N_h(\mathbf{A}
\boxtimes \mathbf{B})_k=\mathbf{A} \times N\bfB_k$ where we view the
set $N\bfB_k$ as a discrete category.
\end{prop}

Like the traditional nerve, the horizontal nerve is fully faithful.

\begin{prop} \label{horizontalfullyfaithful}
The horizontal nerve $\xymatrix@1{N_h\co \mathbf{DblCat_v} \ar[r] &
\mathbf{Cat}^{\Delta^{op}}}$ is fully faithful in the 2-categorical
sense, \ie the functors
$$\xymatrix@C=4pc{\mathbf{DblCat_v}(\bbD,\bbE) \ar[r]^-{(N_h)_{\bbD,\bbE}} & \mathbf{Cat}^{\Delta^{op}}(N_h\bbD,N_h\bbE)}$$
are isomorphisms of categories.
\end{prop}
\begin{pf}
The data of a double functor $F$ and a vertical natural
transformation $\sigma$ are encoded entirely in $N_hF$ and
$N_h\sigma$, so $(N_h)_{\bbD,\bbE}$ is injective on objects and
injective on morphisms.

If $\xymatrix@1{F'\co N_h\bbD \ar[r] & N_h\bbE}$ is a morphism in
$\mathbf{Cat}^{\Delta^{op}}$, then the functors $F'_0$ and $F'_1$
give the data for a double functor $F$, and compatibility with face
and degeneracy maps guarantees compatibility of $F$ with horizontal
composition and units. From $(N_hF)_0=F_0'$ and $(N_hF)_1=F_1'$ it
follows that $N_hF=F'$ using the compatibility with the injective
maps $\xymatrix{e_{i,i+1}\co \{0,1\}  \ar[r] & \{0,\dots,n\}}$
defined by
$$e_{i,i+1}(0)=i$$
$$e_{i,i+1}(1)=i+1$$
for $0 \leq i \leq n-1$. Similarly, if $\sigma'$ is a 2-cell in
$\mathbf{Cat}^{\Delta^{op}}$, we can construct a vertical natural
transformation $\sigma$ from $\sigma'_0$ and $\sigma'_1$ such that
$N_h\sigma=\sigma'$. Then $(N_h)_{\bbD,\bbE}$ is surjective on
objects and surjective on morphisms.
\end{pf}

Also like the traditional nerve, the horizontal nerve is
2-coskeletal, which means that the component $\xymatrix@1{N_h\bbD
\ar[r] & csk_2tr_2 N_h\bbD}$ of the unit for the adjunction
\begin{equation*}
\xymatrix@C=4pc{\mathbf{Cat}^{\Delta^{op}} \ar@{}[r]|{\perp}
\ar@/^1pc/[r]^{tr_2} & \ar@/^1pc/[l]^{csk_2}
\mathbf{Cat}^{\Delta^{op}_2}}
\end{equation*}
is an isomorphism of simplicial objects in $\mathbf{Cat}$. To prove
this, we need a proposition from enriched category theory, the first
part of which is Theorem 5.13 of \cite{kelly3}.

\begin{prop}[Proposition 1.1 of \cite{lackpaoli}] \label{rightKanextensions}
Let $\mathcal{V}$ be a symmetric monoidal closed category which is
complete and cocomplete. If $$\xymatrix@1{\mathcal{A} \ar[r]^I &
\mathcal{B} \ar[r]^J & \mathcal{C} }$$ are  $\mathcal{V}$-functors,
and $J$ is fully faithful, then $J$ is dense provided that $JI$ is
so, and then the identity $JI=JI$ exhibits $J$ as the left Kan
extension of $JI$ along $I$. Furthermore, the singular functor
$\xymatrix@1{\mathcal{C}(J, 1)\co \mathcal{C} \ar[r] &
[\mathcal{B}^{op}, \mathcal{V}]}$ can then be obtained by first
applying the singular functor $$\xymatrix@1{\mathcal{C}(JI,1)\co
\mathcal{C} \ar[r] & [\mathcal{A}^{op}, \mathcal{V}]}$$ and then
right Kan extending along $\xymatrix@1{I\co \mathcal{A}^{op} \ar[r]
& \mathcal{B}^{op}}$.
\end{prop}

\begin{prop} \label{horizontalnerve2coskeletal}
The horizontal nerve of a small double category is 2-coskeletal.
\end{prop}
\begin{pf}
For $I$ and $J$ in Proposition \ref{rightKanextensions}, we take
$$\xymatrix{\Delta_2 \ar[r]^I & \Delta \ar[r]^-J & \mathbf{DblCat_v}}$$
where $I$ is the inclusion of the full subcategory $\Delta_2$ of
$\Delta$ on the objects $[0]$,$[1]$, and $[2]$, and $J=\bbH$ as in
the proof of Theorem \ref{horizontalcategorificationexists}.
Clearly, $J$ is fully faithful in the 2-categorical sense. We denote
by $N_h^2$ the 2-truncation of the horizontal nerve, which is the
singular functor $\mathbf{DblCat_v}(JI,1)$. By the same argument as
in Proposition \ref{horizontalfullyfaithful}, $N_h^2$ is fully
faithful in the 2-categorical sense, which is equivalent to the
density of $JI$ according to Theorem 5.1 (ii) of \cite{kelly3}. By
Proposition \ref{rightKanextensions}, the horizontal nerve
$N_h=\mathbf{DblCat_v}(J,1)$ is 2-naturally isomorphic to the
composite
$$\xymatrix{\mathbf{DblCat_v} \ar[r]^-{N_h^2} & \mathbf{Cat}^{\Delta_2^{op}}
\ar[r]^{Ran_I} & \mathbf{Cat}^{\Delta^{op}} }.$$ Evaluating this
2-natural isomorphism at a double category $\bbD$, we obtain an
isomorphism $$\xymatrix@1{N_h\bbD\ar[r] & csk_2tr_2 N_h\bbD}$$
between $N_h\bbD$ and a 2-coskeletal simplicial object. From the
naturality of the unit, it follows that $N_h\bbD$ is also
2-coskeletal.
\end{pf}

The second nerve we introduce in this section is the bisimplicial
nerve, which we geometrically realize to get a classifying space.
Let $\mathbf{SSet^2}$ denote the category of bisimplicial sets, \ie
functors from $\Delta^{op} \times \Delta^{op}$ into $\mathbf{Set}$.
Since
$$\mathbf{Cat}(\Delta^{op}\times\Delta^{op},\mathbf{Set})\cong\mathbf{Cat}(\Delta^{op},\mathbf{Set}^{\Delta^{op}} )$$
we see that that the category of bisimplicial sets is isomorphic to
the category $\mathbf{SSet}^{\Delta^{op}}$ of simplicial objects in
{\bf SSet}.

\begin{defn}
The {\it bisimplicial nerve} or {\it double nerve} of a double
category $\mathbb{D}$ is the bisimplicial set $N_d\mathbb{D}$ with
$(m,n)$-bisimplices given by composable $m \times n$ arrays
$$
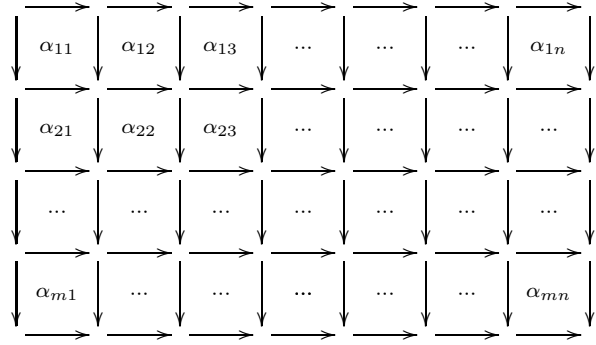
\xymatrix{\ar[r] \ar[d] \ar@{}[dr]|{\alpha_{11}} & \ar[r] \ar[d] \ar@{}[dr]|{\alpha_{12}} & \ar[r] \ar[d]
\ar@{}[dr]|{\alpha_{13}} & \ar[r]\ar[d]  \ar@{}[dr]|{\cdots} & \ar[d] \ar[r] \ar@{}[dr]|{\cdots} & \ar[d] \ar[r] \ar@{}[dr]|{\cdots} & \ar[r] \ar[d] \ar@{}[dr]|{\alpha_{1n}} &
\ar[d]  \\ \ar[r] \ar[d] \ar@{}[dr]|{\alpha_{21}} & \ar[r] \ar[d] \ar@{}[dr]|{\alpha_{22}} & \ar[r] \ar[d] \ar@{}[dr]|{\alpha_{23}} &
\ar[r]\ar[d] \ar@{}[dr]|{\cdots} & \ar[d] \ar[r] \ar@{}[dr]|{\cdots} & \ar[d] \ar[r] \ar@{}[dr]|{\cdots} & \ar[r] \ar[d] \ar@{}[dr]|{\cdots} &
\ar[d] \\ \ar[r] \ar[d] \ar@{}[dr]|{\cdots} & \ar[r] \ar[d] \ar@{}[dr]|{\cdots} & \ar[r] \ar[d] \ar@{}[dr]|{\cdots} &
\ar[r]\ar[d] \ar@{}[dr]|{\cdots} & \ar[d] \ar[r] \ar@{}[dr]|{\cdots} & \ar[d] \ar[r] \ar@{}[dr]|{\cdots} & \ar[r] \ar[d] \ar@{}[dr]|{\cdots} &
\ar[d] \\ \ar[r] \ar[d] \ar@{}[dr]|{\alpha_{m1}} & \ar[r] \ar[d] \ar@{}[dr]|{\cdots} & \ar[r] \ar[d] \ar@{}[dr]|{\cdots} & \ar@{}[dr]|{\cdots}
\ar[r] \ar[d] \ar@{}[dr]|{\cdots} & \ar[d] \ar[r] \ar@{}[dr]|{\cdots} & \ar[d] \ar[r] \ar@{}[dr]|{\cdots} & \ar[r] \ar[d] \ar@{}[dr]|{\alpha_{mn}} & \ar[d]
\\ \ar[r] & \ar[r] & \ar[r] & \ar[r]& \ar[r] & \ar[r] & \ar[r] &}$$
of squares $\alpha_{ij}$ in $\mathbb{D}$. The 0th bisimplicial face
maps correspond to removing the first column respectively row, while
the last bisimplicial face maps correspond to removing the last
column respectively row. The inner bisimplicial face maps correspond
to composing two adjacent columns respectively two adjacent rows.
Bisimplicial degeneracy maps correspond to inserting a trivial
column respectively trivial row of squares.
\end{defn}

In particular $(N_d\mathbb{D})_{0,0}$ is the set of
objects of $\mathbb{D}$, $(N_d\mathbb{D})_{0,n}$ consists of paths
of $n$ horizontal morphisms, and $(N_d\mathbb{D})_{n,0}$ consists of
paths of $n$ vertical morphisms. The simplicial set $(N_d\mathbb{D})_{0,\ast}$ is the nerve of the horizontal 1-category
of $\bbD$ and the simplicial set $(N_d\mathbb{D})_{\ast,0}$ is the nerve of the vertical 1-category of $\bbD$.

\begin{defn} \label{classifyingspacefunctor}
The {\it classifying space functor $B$} is the composite
$$\xymatrix{\mathbf{DblCat} \ar[r]^-{N_d} & \mathbf{SSet^2} \ar[r]^{\diag}
& \mathbf{SSet} \ar[r]^{|\cdot|} & \mathbf{Top},}$$ where $\diag$ is
the functor induced by the diagonal and $|\cdot|$ is the geometric
realization.
\end{defn}

The traditional nerve, the horizontal nerve, and the double nerve are related as follows.

\begin{prop}
The functor
$$\xymatrix{\mathbf{DblCat} \ar[r]^{N_h} & {\bf Cat}^{\Delta^{op}} \ar[r]^{N_*} & {\bf SSet}^{\Delta^{op}} }$$
$$\xymatrix{\mathbb{D} \ar@{|->}[r] & ([n] \mapsto N((N_h\mathbb{D})_n))}$$
is naturally isomorphic to $N_d$.
\end{prop}

\begin{cor}
Let $\bfC$ be a 2-category. Consider the bisimplicial set obtained
by taking the nerves of the hom-categories of $\bfC$, viewing the
resulting {\bf SSet}-category as a simplicial object in {\bf Cat}
with constant object set, and then composing this functor
$\xymatrix@1{\Delta^{op} \ar[r] & \mathbf{Cat}}$ with the
traditional nerve functor. This bisimplicial set is naturally
isomorphic to $N_d(\bfC)$ if we view $\bfC$ as a double category
with trivial vertical morphisms.
\end{cor}

\section{Horizontal Categorification}\label{section:categorification}

We now construct a left adjoint $c_h$ to the horizontal nerve $N_h$
by analogy with the usual nerve $N$. Horizontal categorification
$c_h$ is appropriately compatible with external products, as we show
in Example \ref{horcatcatssetfromformula} and Proposition
\ref{horcatcatsset}.

We recall from \cite{gabrielzismanfractions} the left adjoint
$\xymatrix@1{c\co\mathbf{SSet} \ar[r] & \mathbf{Cat}}$ to the nerve
functor $N$. For a simplicial set $X$, the category $cX$ is the {\it
fundamental category of $X$}, or {\it categorification of $X$}. It
is the free category on the reflexive graph $(X_0,X_1)$ modulo the
smallest congruence such that for every $\tau \in X_2$ with edges
$$\xymatrix{& \ar@{}[d]|{\overset{\phantom{X}}{\tau}} \ar[dr]^g & \\ \ar[ur]^f \ar[rr]_h & &  }$$
we have $g \circ f \sim h$. The following proof is our guideline for
the left adjoint $c_h$ to $N_h$.

\begin{prop} \label{classicalcategorification}
Categorification $c$ is left adjoint to the nerve functor $N$.
\end{prop}
\begin{pf}
We need to construct a natural bijection
$$\mathbf{Cat}(cX,\bfA)\cong\mathbf{SSet}(X,N\bfA).$$
Suppose we have a map $\xymatrix@1{G\co X \ar[r] & N\bfA}$ of
simplicial sets. The 1-truncation is a morphism of reflexive graphs,
so there is a unique functor $J$ making the upper left triangle
commute.
$$\xymatrix@R=3pc@C=3pc{(X_0,X_1) \ar[r]^-{(G_0,G_1)} \ar[d] & \bfA \\ \text{FreeCat}(X_0,X_1)
\ar@{.>}[ur]|{{}_{\phantom{X}}\exists ! \, J^{\phantom{X}}} \ar[r] & cX
\ar@{.>}[u]_{\exists ! \, G'} }$$
Since $J$ comes from a morphism of simplicial sets, the functor $J$ takes congruent morphisms
to equal ones. Therefore there exists a unique functor $G'$ making the lower right triangle commute.

For the converse, given a functor $\xymatrix@1{G'\co cX \ar[r] &
\bfA}$, we compose it with the morphism of reflexive graphs
$$\xymatrix{(X_0,X_1) \ar[r] & \text{FreeCat}(X_0,X_1) \ar[r] & cX}$$
to obtain a morphism of reflexive graphs $\xymatrix@1{(G_0,G_1)\co
(X_0,X_1) \ar[r] & \bfA}$. We define $\xymatrix@1{G_2\co X_2 \ar[r]
& (N\bfA)_2}$ on $\sigma \in X_2$ as
$$\xymatrix@C=4pc{\ar[r]^{G_1(d_2\sigma)} & \ar[r]^{G_1(d_0\sigma)}&}.$$
By definition $G_2$ is compatible with $d_2$ and $d_0$, and the
quotient in the definition of $cX$ makes $G_2$ compatible with
$d_1$. This, together with the simplicial identities relating $X_1$
and $X_2$
$$d_0s_0=id_{X_1} \hspace{.5in} d_2s_0=s_0d_1 \hspace{.5in} d_0s_1=s_0d_0 \hspace{.5in} d_2s_1=id_{X_1},$$
implies that $G_2$ is also compatible with the degeneracies $s_0$
and $s_1$. Since $N\bfA$ is 2-coskeletal this morphism
$(G_0,G_1,G_2)$ of 2-truncated simplicial sets induces a morphism
$\xymatrix@1{G\co X \ar[r] & N \bfA}$ of simplicial sets:
$$Mor(tr_2X,tr_2N\bfA)\cong Mor(X,csk_2tr_2N\bfA)\cong Mor(X,N\bfA).$$

The two procedures $G \mapsto G'$ and $G' \mapsto G$ are inverse to
one another.
\end{pf}

\begin{rmk} \label{morphismdeterminedby1truncation}
Any morphism of simplicial sets $\xymatrix@1{G\co X \ar[r] & N \bfA}$
is completely determined by its 1-truncation $(G_0,G_1)$ as follows. We let
 $\xymatrix{e_{i,i+1}\co \{0,1\}  \ar[r] &
\{0,\dots,n\}}$ be the injective map defined by
$$e_{i,i+1}(0)=i$$
$$e_{i,i+1}(1)=i+1$$
for $0 \leq i \leq n-1$, and we let
$\xymatrix{e_i \co \{0\}  \ar[r] &
\{0,\dots,n\}}$ be the injective map defined by
$$e_i(0)=i$$
for $0 \leq i \leq n$. If $\sigma$ is an $n$-simplex, then
$G(\sigma)$ is the string of $n$ morphisms in $\bfA$
$$\xymatrix@C=4pc{ \ar[r]^-{G(e_{0,1}^*(\sigma))} &  \ar[r]^-{G(e_{1,2}^*(\sigma))}
&  \cdots & \cdots \ar[r]^-{G(e_{n-1,n}^*(\sigma))} & }$$
where the source and target of $G(e_{i,i+1}^*(\sigma))$ are $G(e_i^*(\sigma))$ and $G(e_{i+1}^*(\sigma))$.
\end{rmk}

We turn next to the left adjoint of the horizontal nerve. We will
exhibit two proofs that the horizontal categorification of the
product of a category with a simplicial set is an external product
of the category with the fundamental category of the simplicial set.
This is done in Example \ref{horcatcatssetfromformula} using the
definition of horizontal categorification, while it is done in
Proposition \ref{horcatcatsset} using weighted colimits.

\begin{defn} \label{horizontalcategorification}
Let $X \in \mathbf{Cat}^{\Delta^{op}}$. We define a double category $c_hX$
called the {\it horizontal categorification} or {\it fundamental
double category of $X$} as follows. First we define a double derivation
scheme $\bbS$ with vertical 1-category $X_0$ and with horizontal
1-category the fundamental category of the simplicial set $\Obj X$.
The squares of $\bbS$ are the morphisms of $X_1$. We equip the free
double category $\bbF$ on the double derivation scheme $\bbS$ with
the smallest congruence $\sim$ such that
\begin{enumerate}
\item
If $\alpha, \beta \in \Mor X_1$ are composable in $X_1$,
then the vertical composite $\vcomp{\alpha}{\beta}$ in $\bbF$ is congruent to the composite of
$\beta$ and $\alpha$ in $X_1$,
\item
For all $\tau \in \Mor X_2$ with boundary
$$\xymatrix{& \ar@{}[d]|{\overset{\phantom{X}}{\tau}} \ar[dr]^\beta & \\ \ar[ur]^\alpha \ar[rr]_\gamma & &  }$$
we have  $$[\alpha \beta ] \sim \gamma,$$
\item
For any vertical morphism $j$, that is, for any $j\in\Mor X_0$, the
horizontal identity $i_j^h$ is congruent to the degeneracy of $j$ in
$\Mor X_1$.
\item
For any $f \in \Obj X_1$, the vertical identity square $i^v_f$ on
the image of $f$ in the horizontal 1-category of $\bbS$ is congruent
to the identity on $f$ in the category $X_1$.
\end{enumerate}
We define $c_hX$ as the quotient of $\bbF$ by the congruence $\sim$. The horizontal and
vertical 1-categories of $c_hX$ are the horizontal and vertical 1-categories of $\bbS$.
\end{defn}

\begin{rmk}
In the definition of horizontal categorification it is not necessary
to mod out by additional relations to make the identity squares
functorial. If $g \circ f \sim h$ in the horizontal 1-category
because of $\tau \in \Obj X_2$, then the identity morphism on $\tau$
in the category $X_2$ implies we have $i^v_{h} \sim [i^v_f \,
i^v_g]$ (the face maps are functors and we have (ii) and (iv)). For
vertically composable morphisms $j$ and $k$, we have
$$i^h_{\vcomp{j}{k}}\sim\vcomp{i^h_j}{i^h_k}$$ because degeneracy is a
functor and by (i) and (iii).
\end{rmk}

\begin{examp} \label{horcatsset}
If $X$ is a simplicial set, then $c_h \nu X=\bbH cX$. By definition, the horizontal
1-category is $cX$, and the vertical 1-category is the discrete category $X_0$. Since $X_1$ is also
discrete, there are no nontrivial squares.
\end{examp}

\begin{examp} \label{horcatcatssetfromformula}
Recall from Proposition \ref{horizontalnerveexternalproduct} that
$\xymatrix@1{\sigma\co \mathbf{Cat} \ar[r] &
\mathbf{Cat}^{\Delta^{op}}}$ denotes the constant functor and
$\xymatrix@1{\nu\co \mathbf{Set}^{\Delta^{op}} \ar[r] &
\mathbf{Cat}^{\Delta^{op}}}$ denotes the inclusion. If {\bf A} is a
category and $Y$ is a simplicial set, then the horizontal
categorification of the simplicial category $\sigma \mathbf{A}
\times \nu Y$ is $\mathbf{A} \boxtimes cY$. In fact, the horizontal
1-category of the double derivation scheme $\bbS$ is
$$c(\Obj (\sigma\mathbf{A} \times \nu Y))=c(\Obj \bfA \times Y)=(\bfH(\bfA \boxtimes cY))_0.$$
The vertical 1-category of $\bbS$ is
$$(\sigma \mathbf{A} \times \nu Y)_0=\bfA \times Y_0 = (\bfV(\bfA \boxtimes cY))_0.$$
The squares of $\bbS$ are
$$\Mor (\sigma \bfA \times \nu Y)_1 = (\Mor\bfA) \times Y_1.$$ The congruence
on $\bbF$ corresponds precisely to the relations in $\mathbf{A}
\boxtimes cY$ for pairwise compositions of squares and identity
squares. We present an alternative conceptual proof of this example
in Proposition \ref{horcatcatsset}.
\end{examp}

\begin{prop}
Horizontal categorification $c_h$ is left adjoint to the horizontal nerve $N_h$.
\end{prop}
\begin{pf}
We use the notation of Definition \ref{horizontalcategorification}
and construct a natural bijection $$\mathbf{DblCat}(c_hX,\bbD) \cong
[\Delta^{op},\mathbf{Cat}](X,N_h\bbD).$$ Suppose $\xymatrix@1{G\co X
\ar[r] & N_h\bbD}$ is a morphism of simplicial objects in {\bf Cat}.
This induces a morphism of double derivation schemes
$\xymatrix@1{\bbS \ar[r] & \bbD}$ and a unique double functor $J$
making the upper left triangle commute,
$$\xymatrix@R=3pc@C=3pc{\bbS \ar[r] \ar[d] & \bbD
\\ \bbF \ar[r] \ar@{.>}[ur]|{{}_{\phantom{X}}\exists ! \, J^{\phantom{X}}}
& c_hX \rlap{\,.}\ar@{.>}[u]_{\exists ! \, G'} }$$
Since $G$ is a morphism of simplicial objects in $\mathbf{Cat}$, $J$ takes congruent squares
to equal squares, and there exists a unique double functor $G'$ making the lower right triangle commute.

On the other hand, given a double functor $\xymatrix@1{G'\co c_hX
\ar[r] & \bbD}$ we compose it with
$$\xymatrix{(X_0,X_1) \ar[r] & \bbF \ar[r] & c_hX}$$
to obtain a morphism $(G_0,G_1)$ of 1-truncated simplicial objects
in $\mathbf{Cat}$. By the same argument as in the proof of
Proposition \ref{classicalcategorification}, we obtain maps of
simplicial sets
$$\xymatrix{G^{\text{\it Obj}}\co\Obj X \ar[r] & \Obj N_h\bbD}$$
$$\xymatrix{G^{\text{\it Mor}}\co\Mor X \ar[r] & \Mor N_h\bbD}.$$
The maps $(G^{\text{\it Obj}}_0, G^{\text{\it Mor}}_0)$ and
$(G^{\text{\it Obj}}_1, G^{\text{\it Mor}}_1)$ are already known
to be functors, and $X$ is a simplicial object in $\mathbf{Cat}$.

We claim that $(G^{\text{\it Obj}}_2, G^{\text{\it Mor}}_2)$ is also
a functor. If $\sigma,\sigma'\in \Mor X_2$ are composable we denote
their composite as $\left[ \begin{array}{c} \sigma \\ \hline \sigma'
\end{array} \right]$. We similarly denote composites in $(N_h\bbD)_2$. To indicate a horizontal path of squares in
$\bbD$, we separate the squares by a dot. Then $G_2$ preserve
compositions
$$\aligned
G_2\left[ \begin{array}{c} \sigma \\ \hline \sigma'
\end{array} \right]&=G_1d_2\left[ \begin{array}{c} \sigma \\ \hline \sigma'
\end{array} \right]\cdot G_1d_0\left[ \begin{array}{c} \sigma \\ \hline \sigma'
\end{array} \right] \\
&=\left[ \begin{array}{c} G_1d_2\sigma \\ \hline G_1d_2\sigma'
\end{array} \right]\cdot\left[ \begin{array}{c} G_1d_0\sigma \\ \hline G_1d_0 \sigma'
\end{array} \right] \\
&=\left[ \begin{array}{c} G_1d_2\sigma \cdot G_1d_0\sigma \\ \hline
G_1d_2\sigma' \cdot G_1d_0\sigma'
\end{array} \right] \\
&=\left[ \begin{array}{c} G_2\sigma \\ \hline G_2\sigma'
\end{array} \right],
\endaligned$$
and $G_2$ preserves units similarly.

Since $N_h\bbD$ is 2-coskeletal by Proposition
\ref{horizontalnerve2coskeletal}, the 2-truncated morphism
$(G_0,G_1,G_2)$ induces a morphism $\xymatrix@1{X \ar[r] &
N_h\bbD}$, and by Remark \ref{morphismdeterminedby1truncation}, this
morphism must be $(G^{\text{\it Obj}},G^{\text{\it Mor}})$ from
above.

The two procedures $G \mapsto G'$ and $G' \mapsto G$ are inverse to
one another.
\end{pf}

\begin{prop} \label{catembedding}
Consider $\mathbf{Cat}$ embedded into $\mathbf{Cat}^{\Delta^{op}}$
as the constant simplicial objects, and consider $\mathbf{Cat}$
embedded vertically into $\mathbf{DblCat}$. Then the adjunction $c_h
\dashv N_h$ restricts to the identity adjunction on these full
subcategories.
\end{prop}

We now move towards a conceptual proof of Example \ref{horcatcatssetfromformula} in Proposition
\ref{horcatcatsset}.

\begin{rmk} \label{recallingcopower}
Recall that if $S$ is a set and $A$ is an object of a category, then
the {\it copower} $S \cdot A$ is the coproduct of $A$ with itself
$S$ times. In some categories, the copower has a simple description.
For example, if $\bfC$ is a category, then the copower in
$\mathbf{Cat}$ is
\begin{equation*}
S \cdot \bfC=\coprod_S \bfC=S \times \bfC.
\end{equation*}
If $X$ is a simplicial set and $Y \in [\Delta^{op},\mathbf{Cat}]$,
then $X\cdot Y$ is the simplicial object in $\mathbf{Cat}$
\begin{equation*}
\xymatrix{[n] \ar@{|->}[r] & X_n \cdot Y_n =\coprod_{X_n} Y_n = X_n
\times Y_n},
\end{equation*}
 which is the same as $\nu X \times Y$.
\end{rmk}

\begin{lem} \label{productisweightedlimit}
If $X$ and $Y$ are simplicial objects in $\mathbf{Cat}$, then $X
\times Y$ is the weighted colimit $X*G$ of the
$\mathbf{Cat}$-functor $$\xymatrix{G\co \Delta \ar[r] &
[\Delta^{op},\mathbf{Cat}]}$$
$$\xymatrix{[n] \ar@{|->}[r] & Y \times \nu \Delta[n]}$$ with weighting
$\xymatrix@1{X\co \Delta^{op} \ar[r] & \mathbf{Cat}}$.
\end{lem}
\begin{pf}
Since $(\mathbf{Cat},\times)$ is symmetric monoidal closed, it
follows from a general fact that $[\Delta^{op}, \mathbf{Cat}]$ has a
tensor product
$$(Y \otimes Z)_n:=Y_n \times Z_n=(Y\times Z)_n$$
and an internal hom
$$\aligned
\left[Y,Z\right]_n :&=[\Delta^{op},\mathbf{Cat}](\Delta[n]\cdot Y,Z) \\
&\cong [\Delta^{op},\mathbf{Cat}](Y  \times \nu \Delta[n],Z)
\endaligned$$
for all $Y,Z \in [\Delta^{op},\mathbf{Cat}]$.

 For any $Z \in
[\Delta^{op},\mathbf{Cat}]$,
$$[\Delta^{op},\mathbf{Cat}](G([n]),Z)$$ is the $n$-th
category of the internal hom $[Y,Z]$. Thus we have a natural
isomorphism
$$[\Delta^{op},\mathbf{Cat}](X,[\Delta^{op},\mathbf{Cat}](G(-),Z)) \cong [\Delta^{op},\mathbf{Cat}](X \times Y,Z)$$
and $X \times Y$ satisfies the universal property of the weighted
colimit $X*G$.
\end{pf}

We finish the conceptual proof of Example
\ref{horcatcatssetfromformula}.

\begin{prop} \label{horcatcatsset}
If {\bf A} is a category and $Y$ is a simplicial set, then the
horizontal categorification of the simplicial category $\sigma
\mathbf{A} \times \nu Y$ is $\mathbf{A} \boxtimes cY$ where $cY$ is
the traditional categorification of $Y$.
\end{prop}
\begin{pf}
By Lemma \ref{productisweightedlimit}, $\sigma \bfA \times \nu Y$ is
the weighted colimit $\sigma \bfA * G$ of
$$\xymatrix{G\co\Delta \ar[r] &
[\Delta^{op},\mathbf{Cat}]}$$
$$\xymatrix{[n] \ar@{|->}[r] & \nu Y \times \nu \Delta[n]}$$
with weighting $\xymatrix@1{\sigma \bfA\co \Delta^{op} \ar[r] &
\mathbf{Cat}}$.

Let $\xymatrix@1{J\co \Delta \ar[r] & \mathbf{DblCat}}$ be the
horizontal embedding. Then by Theorem 4.51 of \cite{kelly3}, for
each $Z \in [\Delta^{op},\mathbf{Cat}]$, $c_h(Z)\cong  Z*J.$ Hence
\begin{equation} \label{weightedcategorification1}
\aligned
c_h(\sigma \mathbf{A} \times \nu Y) &=c_h(\sigma \bfA * G) \\
& \cong (\sigma \bfA * G)*J \\
& \cong \sigma \bfA * (G*J)
\endaligned
\end{equation}
by the general Fubini Theorem, which is equation (3.23) in
\cite{kelly3}. The functor $\xymatrix@1{G*J\co \Delta \ar[r] &
\mathbf{DblCat}}$ in the last line takes $[n]$ to
$$G([n])*J \cong c_h(G([n]))=c_h(\nu Y \times \nu \Delta[n]).$$
From Example \ref{horcatsset} and the fact that $c$ preserves finite
products, we have
$$c_h(\nu Y \times \nu \Delta[n])=\bbH c(Y \times \Delta[n])\cong
\bbH cY \times \bbH[n].$$ We conclude that
(\ref{weightedcategorification1}) has the form
\begin{equation} \label{weightedcategorification2}
c_h(\sigma \bfA \times \nu Y) \cong \sigma \bfA * (\bbH cY \times
\bbH[-]).
\end{equation}

We claim that the right hand side of
(\ref{weightedcategorification2}) is isomorphic to $\bbV \bfA \times
\bbH cY$. In fact, Proposition \ref{likeinternalhom} and the
adjunction $sk_0 \dashv tr_0$ give, for all $\bbE \in
\mathbf{DblCat_v}$
$$\aligned
\mathbf{DblCat_v}(\bbV \bfA \times \bbH cY, \bbE)) & \cong
\mathbf{Cat}(\bfA,\mathbf{DblCat_v}(\bbH cY,\bbE)) \\
& \cong \mathbf{Cat}(\bfA,tr_0 \,\mathbf{DblCat_v}(\bbH
cY \times \bbH[-],\bbE)) \\
& \cong [\Delta^{op},\mathbf{Cat}](\sigma
\bfA,\mathbf{DblCat_v}(\bbH cY \times \bbH[-],\bbE)).
\endaligned$$
The claim follows now from the definition of weighted colimit. Hence, (\ref{weightedcategorification2})
implies that
$$c_h(\sigma \bfA \times \nu Y) \cong \bbV \bfA \times \bbH cY=\bfA
\boxtimes cY.$$
\end{pf}

The vertical categorification of a simplicial object $X$ in $\mathbf{Cat}$ is the transpose
of $c_hX$.

\section{Model Structures Arising from $\mathbf{Cat}^{\Delta^{op}}$} \label{section:simplicialcategories}

Now that we have the adjunction $c_h \dashv N_h$ in place we can use
it to transfer model structures from $\mathbf{Cat}^{\Delta^{op}}$ to
$\mathbf{DblCat}$ using Kan's Lemma on Transfer (Theorem \ref{Kan}).
This theorem says that one can lift a model structure across an
adjunction under certain smallness conditions, which guarantee
functorial factorizations. This is our first method for constructing
model structures on {\bf DblCat}. In Section
\ref{section:topologies} we will adopt the point of view of double
categories as internal categories and apply the results of
\cite{everaertinternal}. In Section \ref{2monadstructure} we will
consider {\bf DblCat} as a category of algebras for a 2-monad and
use \cite{lack2monads}.

The category $\mathbf{Cat}^{\Delta^{op}}$ has four model structures
of interest to us. These arise as diagram structures and Reedy
structures associated to two cofibrantly generated model structures
on {\bf Cat}: the categorical structure and the Thomason structure.
In Sections \ref{structuresonCat}-\ref{SubsectionOnKan} we review
some material for the reader's convenience: model structures on {\bf
Cat}, their associated diagram structures, smallness arguments, and
Kan's Lemma on Transfer. After these preliminaries, we turn to our
new results. In Sections \ref{subsection:Thomasontransfer} and
\ref{subsection:categoricaltransfer} we transfer the diagram
structures to {\bf DblCat} across the horizontal
categorification-horizontal nerve adjunction, and show that the
transferred structures on $\mathbf{DblCat}$ extend the Thomason
structure and categorical structure on the vertically embedded
subcategory $\mathbf{Cat}$. In the proofs of our transfer results we
crucially need to know the behavior of certain pushouts, and these
are treated in Theorems \ref{pushoutinDblCat} and \ref{pushoutiso}
of the Appendix. We show in Section \ref{subsection:NoReedytransfer}
that the Reedy categorical structure cannot transfer.

Recall the notion of cofibrantly generated model category.
\begin{defn}
A model category $\bfC$ is {\it cofibrantly generated} if there exist sets of morphism $I$ and $J$ in $\bfC$
such that
\begin{enumerate}
\item
The domains of $I$ are small with respect to $I$-cell as defined in Definitions \ref{defnsmallness} and \ref{transfinitecompositionIcell},
\item
The domains of $J$ are small with respect to $J$-cell,
\item
The class of fibrations is precisely the class of morphisms with the right lifting
property with respect $J$,
\item
The class of acyclic fibrations is precisely the class of morphisms with the right lifting
property with respect to $I$.
\end{enumerate}
In this case, $I$ is the set of {\it generating cofibrations} and $J$ is the set of {\it
generating acyclic cofibrations}.
\end{defn}

\subsection{Model Structures on $\mathbf{Cat}$} \label{structuresonCat}
In the {\it Thomason structure} on {\bf Cat} in \cite{thomasonCat} a
functor $F$ is a weak equivalence (respectively fibration) if and
only if $\Ex^2NF$ is a weak equivalence (respectively fibration) of
simplicial sets. The functor $\Ex$ is superfluous for weak
equivalences, as Thomason proved that $F$ is a weak equivalence if
and only if $NF$ is. The functor $\xymatrix@1{\Ex\co \mathbf{SSet}
\ar[r] & \mathbf{SSet}}$ is the left adjoint to barycentric
subdivision $\xymatrix@1{\Sd\co \mathbf{SSet} \ar[r] &
\mathbf{SSet}}$, which we recall below. The Thomason structure is
cofibrantly generated. The generating cofibrations are the
inclusions of categorical boundaries
$$\xymatrix@1{c\Sd^2\partial \Delta[m] \ar[r] & c\Sd^2\Delta[m]}$$ while the
generating acyclic cofibrations are
the inclusions of categorical horns $$\xymatrix@1{c\Sd^2\Lambda^k[m] \ar[r] & c\Sd^2\Delta[m]}.$$

We now recall the definition of {\it barycentric subdivision}
$\Sd$.\label{subdivisionobjects} The simplicial sets $\Sd\Delta[m]$
and $\Sd\Lambda^k[m]$ are respectively the nerves of the posets of
nondegenerate simplices of $\Delta[m]$ and $\Lambda^k[m]$. The
ordering is the face relation. Thus a $q$-simplex of $\Sd\Delta[m]$
is a tuple $(v_0, \dots, v_q)$ of nondegenerate simplices (faces) of
$\Delta[m]$ such that $v_i$ is a face of $v_{i+1}$ for all $0\leq
i\leq q-1$. Such a tuple is a $q$-simplex of $\Sd\Lambda^k[m]$ if
and only if all $v_0, \dots, v_q$ are in $\Lambda^k[m]$. A
$p$-simplex $u$ is a face of a $q$-simplex $v$ in $\Sd\Delta[m]$ if
and only if
$$\{u_0,\dots, u_p\} \subseteq \{v_0,\dots, v_q\}.$$ A $p$-simplex
$u$ of $\Sd\Delta[m]$ is nondegenerate if and only if all $u_i$ are
distinct.

The {\it barycentric subdivision} of a simplicial set $Y$ is by
definition
$$\underset{\Delta[n] \rightarrow Y}{\colim}
\Sd\Delta[n]$$ where the colimit is indexed over the category of
simplices of $Y$. It follows from Page 311 of \cite{thomasonCat}
that $c\Sd^2\Delta[m]$ and $c\Sd^2\Lambda^k[m]$ are respectively the
posets of nondegenerate simplices of $\Sd\Delta[m]$ and
$\Sd\Lambda^k[m]$ and the generating acyclic cofibration
$\xymatrix@1{c\Sd^2\Lambda^k[m] \ar[r] & c\Sd^2\Delta[m]}$ is the
inclusion of these posets.

The other model structure on {\bf Cat} is the {\it categorical
structure} of \cite{joyaltierney}. In the categorical structure a
functor is a weak equivalence if and only if it is an equivalence of
categories. A functor $\xymatrix@1{F\co \bfA \ar[r] & \bfB}$ is a
fibration if for each isomorphism $g\co b \cong Fa$ in $\bfB$ there
is an isomorphism $f\co a' \cong a$ in $\bfA$ such that $Fa'=b$ and
$Ff=g$. These fibrations of categories are also called {\it
isofibrations}. A cofibration is a functor that is injective on
objects. The categorical structure on {\bf Cat} is also cofibrantly
generated. There are three generating cofibrations:
$$\xymatrix{\emptyset \, \ar@{^{(}->}[r] & \{1\}}$$
$$\xymatrix{\{0,1\} \, \ar@{^{(}->}[r] & \{0 \rightarrow 1\}}$$
$$\xymatrix{\{0 \rightrightarrows 1\} \ar[r] & \{0 \rightarrow 1\}}$$
and one generating acyclic cofibration:
$$\xymatrix{\{1\} \, \ar@{^{(}->}[r] & \{0 \cong 1\}=\bfI}.$$

\subsection{Diagram Model Structures on
$\mathbf{Cat}^{\Delta^{op}}$}\label{diagramstructures} Given a model
category $\bfM$ and a small category $\bfC$, one might hope that the
category $\bfM^{\bfC}$ of functors $\xymatrix@1{\bfC \ar[r] & \bfM}$
is also a model category with levelwise weak equivalences and
levelwise fibrations. By this we mean that a natural transformation
is a weak equivalence (respectively fibration) if and only if each
of its components is. Unfortunately, this definition does not always
give rise to a model structure on $\bfM^{\bfC}$. However, if $\bfM$
is a {\it cofibrantly generated} model category, Theorem
\ref{diagramgenerators} guarantees that this definition does indeed
give rise to a model structure on $\bfM^{\bfC}$, which is even
cofibrantly generated.

\begin{thm}[11.6.1 in \cite{hirschhorn}] \label{diagramgenerators}
Let $\bfC$ be a small category and $\bfM$ a cofibrantly generated model category with $I$ the set of
generating cofibrations and $J$ the set of generating acyclic cofibrations. Then $\bfM^{\bfC}$ is a
cofibrantly generated model category with levelwise weak equivalences and levelwise fibrations. The
generating cofibrations are natural transformations of the form
\begin{equation*}
\coprod_{\bfC(C,-)}A \overset{\underset{\bfC(C,-)}{\coprod}
f}{\xymatrix@1@C=4pc{\ar[r] &}} \coprod_{\bfC(C,-)} B
\end{equation*}
for $\xymatrix@1{f\co A \ar[r] & B}$ in $I$. The generating acyclic
cofibrations are defined similarly with $f$ in $J$. A morphism in
$\bfM^{\bfC}$ is a cofibration if it is a retract of a transfinite
composition of pushouts of generating cofibrations. The components
of a cofibration are also cofibrations.
\end{thm}

Thus, the category $\mathbf{Cat}^{\Delta^{op}}$ inherits two model
structures from Section \ref{structuresonCat}. In the {\it diagram
Thomason structure} on $\mathbf{Cat}^{\Delta^{op}}$, a natural
transformation $\alpha$ is a weak equivalence (respectively
fibration) if and only if $N\alpha_i$ is a weak equivalence
(respectively fibration) of simplicial sets for each $i \geq 0$. In
the {\it diagram categorical structure} on
$\mathbf{Cat}^{\Delta^{op}}$, a natural transformation $\alpha$ is a
weak equivalence (respectively fibration) if and only if $\alpha_i$
is an equivalence of categories (respectively isofibration) for all
$i \geq 0$.

If $\bfC$ is a {\it Reedy category}, then a model structure on
$\bfM$ also induces a {\it Reedy model structure} on $\bfM^{\bfC}$
(see for example \cite{hirschhorn} or \cite{hovey}). The category
$\Delta$ is a Reedy category, so the Thomason and categorical
structures on {\bf Cat} also give rise to two more model structures
on $\mathbf{Cat}^{\Delta^{op}}$. However, we do not study these in
more detail because of the following Theorem and also because of
Theorem \ref{NoReedy}.

\begin{thm}[Theorem 15.6.4 of \cite{hirschhorn}]
If $\bfC$ is a Reedy category and $\bfM$ is a cofibrantly generated
model category, then the identity functor of $\bfM^\bfC$ is a left
Quillen equivalence from the cofibrantly generated diagram model
structure to the Reedy model structure, and a right Quillen
equivalence in the opposite direction.
\end{thm}

\subsection{Smallness}

We will need some knowledge about smallness to use Kan's Lemma on
Transfer. We recall some of the relevant notions from \cite{hovey}.
Appropriate smallness conditions also allow us to conclude that a
transfinite composition of weak equivalences is a weak equivalence.

\begin{defn}
Let $\kappa$ be a cardinal. An ordinal $\lambda$ is {\it $\kappa$-filtered} if it is a limit ordinal
and, if $A \subseteq \lambda$ and $|A| \leq \kappa$, then $\text{sup }A < \lambda$.
\end{defn}

\begin{defn} \label{defnsmallness}
Let $\bfC$ be a category with all small colimits and $\kappa$ a
cardinal. An object $A$ of $\bfC$ is called {\it $\kappa$-small} if
for all $\kappa$-filtered ordinals $\lambda$ and all
colimit-preserving functors $\xymatrix@1{X\co \lambda \ar[r] &
\bfC}$ the map of sets
\begin{equation} \label{smallness}
\underset{\beta<\lambda}{\text{colim}}\;\bfC(A,X_\beta) \xymatrix{ \ar[r] & }
\bfC(A,\underset{\beta<\lambda}{\text{colim}}\;X_\beta)
\end{equation}
is a bijection. An object $A$ is said to be {\it small} if it is $\kappa$-small for some cardinal $\kappa$.
An object $A$ is said to be {\it finite} if it is $\kappa$-small for a finite cardinal
$\kappa$, \ie  for {\it any} limit ordinal $\lambda$ and colimit-preserving functor $X$, the map
(\ref{smallness}) is a bijection. We say the concepts hold {\it relative to a class of morphisms $\bfD$ in $\bfC$}
if they hold true for all $X$ with $\xymatrix@1{X_\beta \ar[r] & X_{\beta+1}}$ in $\bfD$ for all $\beta+1 <\lambda$.
\end{defn}

For example, categories are small as follows, and we conclude similarly that double categories are small.

\begin{prop} \label{categoriesaresmall}
Any category $\bfA$ is $\kappa$-small where $$\kappa=|\Obj
\bfA|+|\Mor \bfA|+|\Mor \bfA \,\mbox{}_s\!\!\times_t \Mor \bfA|.$$ In
particular, if $\Mor \bfA$ is a finite set, then $\bfA$ is finite as
an object of {\bf Cat}.
\end{prop}
\begin{pf}
Let $\xymatrix@1{\bfX\co \lambda \ar[r] & \mathbf{Cat}}$ be a
colimit-preserving functor from a $\kappa$-filtered ordinal
$\lambda$. Recall that ordinals are filtered categories and filtered
colimits of categories are formed by simply taking the filtered
colimits of the object set and the morphism set.

Suppose $\xymatrix@1{F\co \bfA \ar[r] & \text{colim}\;\bfX}$ is a
functor. For each $A \in \Obj \bfA$ and $f \in \Mor \bfA$ there are
ordinals $\alpha_1(A)$ and $\alpha_2(f)$ such that $F(A)$ and $F(f)$ are in
the image of $\bfX_{\alpha_1(A)}$ and $\bfX_{\alpha_2(A)}$. Let
$\beta$ be the supremum of all the $\alpha_1(A)$ and $\alpha_2(f)$.
Then $\beta < \lambda$ and we obtain maps of sets
$$\xymatrix@1{G^{\Obj}\co\Obj \bfA \ar[r] & \Obj \bfX_{\beta} }$$
$$\xymatrix@1{G^{\Mor}\co\Mor \bfA \ar[r] & \Mor \bfX_{\beta} }$$
which factor the functor $F$. There exists for each $f \in \Mor
\bfA$ an index $\gamma(f)$ such that $s(G(f))=G(s(f))$ and
$t(G(f))=G(t(f))$ in $\bfX_{\gamma(f)}$. For each $A \in \Obj \bfA$
there is an index $\delta(A)$ such that $G(1_A)=1_{G(A)}$ in
$\bfX_{\delta(A)}$. For each $(\ell,k) \in \Mor \bfA \,\mbox{}_s\!\!\times_t
\Mor \bfA$ there exists an index $\epsilon(\ell,k)$ such that
$G(\ell \circ k)=G(\ell) \circ G(k)$ in $\bfX_{\epsilon(\ell,k)}$.
Let $\zeta$ be the supremum of all these indices
$\gamma,\delta,\epsilon$. Then $\zeta < \lambda$ and $G$ induces a
functor $\xymatrix@1{\bfA \ar[r] & \bfX_{\zeta}}$ which factors $F$.
Hence (\ref{smallness}) is onto.

Suppose $\xymatrix@1{M\co\bfA \ar[r] & \bfX_{\alpha}}$ and
$\xymatrix@1{N\co\bfA \ar[r] & \bfX_{\beta}}$ are functors that
become equal in the colimit. Then for each $A\in \Obj \bfA$ and each
$f \in \Mor \bfA$ there are indices $\gamma(A)$ and $\delta(f)$ such
that $M(A)=N(A)$ and $M(f)=N(f)$ in $\bfX_{\gamma(A)}$ and
$\bfX_{\delta(f)}$ respectively. Let $\zeta<\lambda$ be the supremum
of all these indices $\gamma(A)$ and $\delta(f)$. Then $M$ and $N$
become equal at the stage $\zeta$ and the map (\ref{smallness}) is
injective.
\end{pf}

\begin{prop} \label{doublecategoriesaresmall}
Let $\bbD$ be a double category and $s^h,s^v,t^h,t^v$ the horizontal and vertical source and target maps.
Then $\bbD$ is $\kappa$-small where $$\aligned \kappa=&|\Obj \bbD|+|\Hor \bbD|+|\Hor \bbD \, \,\mbox{}_{s^h}\!\!\times_{t^h}
\Hor \bbD|
\\&+|\Ver \bbD|+|\Ver \bbD \,\,\mbox{}_{s^v}\!\!\times_{t^v} \Ver \bbD|\\&+|\Sq \bbD|+|\Sq \bbD \, \,\mbox{}_{s^v}\!\!\times_{t^v}
\Sq \bbD|\\&+|\Sq \bbD \, \,\mbox{}_{s^h}\!\!\times_{t^h} \Sq \bbD|. \endaligned$$
In particular, if $\Sq\bbD$ is a finite set, then $\bbD$ is finite as an object of {\bf DblCat}.
\end{prop}
\begin{pf}
We first obtain a map of the underlying quadruple of sets, and then we go out far enough to make it into a double functor
by considering the various compositions and identities as in Proposition \ref{categoriesaresmall}.
\end{pf}
Note that this proposition easily generalizes to $n$-fold categories.

One useful application of finiteness is to transfinite compositions of weak equivalences.

\begin{defn} \label{transfinitecompositionIcell}
If $\bfC$ is a category with all small colimits, $\lambda$ is an
ordinal, $\bfD$ is a class of morphisms in $\bfC$, and
$\xymatrix@1{X\co\lambda \ar[r] & \bfC}$ is a colimit preserving
functor such that $\xymatrix@1{X_\beta \ar[r] & X_{\beta+1}}$ is in
$\bfD$ for all $\beta+1 <\lambda$, then the morphism
$$\xymatrix@1{X_0 \ar[r] & \text{colim}\;X}$$ is called a {\it transfinite composition of morphisms in $\bfD$}. If
$I$ is a class of morphisms in $\bfC$, then a transfinite composition of pushouts of elements of $I$ is called a
{\it relative $I$-cell complex}. The class of relative $I$-cell complexes is denoted {\it $I$-cell}.
\end{defn}

\begin{prop}[7.4.2 of \cite{hovey}]
Suppose $\bfC$ is a cofibrantly generated model category in which the domains and codomains of the generating
cofibrations and generating acyclic cofibrations are finite. Then every transfinite composition of weak
equivalences is a weak equivalence.
\end{prop}

\begin{examp} \label{simplicialtransfinitecompositions}
In both the Thomason structure and the categorical structure on {\bf Cat}, every transfinite composition of weak
equivalences is a weak equivalence, as the domains and codomains of the generating cofibrations and
generating acyclic cofibrations only have finitely many morphisms. Since weak equivalences and
colimits in ${\bf Cat}^{\Delta^{op}}$ are levelwise, every transfinite composition of weak
equivalences in the diagram structures is also a weak equivalence.
\end{examp}

\subsection{Kan's Lemma on Transfer} \label{SubsectionOnKan}
Our first main tool for constructing model structures on {\bf
DblCat} is Kan's Lemma on Transfer. The form we will use is
Corollary \ref{KanCorollary}.

\begin{thm}[Kan's Lemma on Transfer, 11.3.2 in \cite{hirschhorn}] \label{Kan}
Let $\mathbf{C}$ be a cofibrantly generated model category with
generating cofibrations $I$ and generating acyclic cofibrations $J$.
Suppose {\bf D} is complete and cocomplete, and that
\begin{equation} \label{adjunction}
\xymatrix@C=4pc{{\bf C} \ar@{}[r]|{\perp} \ar@/^1pc/[r]^{F} &
\ar@/^1pc/[l]^{G} {\bf D}}
\end{equation}
is an adjunction. Assume the following.
\begin{enumerate}
\item
For every $i\in I$, $\dom Fi$ is small with respect to $FI$-cell.
For every $j\in J$, $\dom Fj$ is small with respect to $FJ$-cell.
\item
The functor $G$ maps every relative $FJ$-complex to a weak
equivalence in {\bf C}.
\end{enumerate}
Then there exists a cofibrantly generated model structure on {\bf D}
with generating cofibrations $FI$ and generating acyclic
cofibrations $FJ$. Further, $f$ is a weak equivalence in {\bf D} if
and only $G(f)$ is a weak equivalence in {\bf C}, and $f$ is a
fibration in {\bf D} if and only $G(f)$ is a fibration in {\bf C}.
\end{thm}

Along the lines of \cite{worytkiewicz2Cat}, we have the following
corollary.

\begin{cor} \label{KanCorollary}
Let $\mathbf{C}$ be a cofibrantly generated model category with
generating cofibrations $I$ and generating acyclic cofibrations $J$.
Suppose {\bf D} is complete and cocomplete, and that $F \dashv G$ is
an adjunction as in (\ref{adjunction}). Assume the following.
\begin{enumerate}
\item
For every $i\in I$ and $j \in J$, the objects $\dom Fi$ and $\dom
Fj$ are small with respect to the entire category {\bf D}.
\item
For any ordinal $\lambda$ and any colimit preserving functor
$\xymatrix@1{X\co\lambda \ar[r] & {\bf C}}$ such that
$\xymatrix@1{X_{\beta} \ar[r] & X_{\beta+1}}$ is a weak equivalence,
the transfinite composition
$$\xymatrix{X_0 \ar[r] & \mbox{\text{colim}}\;X}$$
is a weak equivalence.
\item
$G$ preserves filtered colimits.
\item
If $j'$ is a pushout of $F(j)$ in {\bf D} for $j \in J$, then
$G(j')$ is a weak equivalence in {\bf C}.
\end{enumerate}
Then there exists a cofibrantly generated model structure on {\bf D}
with generating cofibrations $FI$ and generating acyclic
cofibrations $FJ$. Further, $f$ is a weak equivalence in {\bf D} if
and only $G(f)$ is a weak equivalence in {\bf C}, and $f$ is a
fibration in {\bf D} if and only $G(f)$ is a fibration in {\bf C}.
\end{cor}
\begin{pf}
Clearly, (i) of Theorem \ref{Kan} follows from the hypotheses. To
see (ii), we recall that a relative $FJ$-complex is a transfinite
composition of pushouts of morphisms $Fj$ where $j \in J$. If a
relative $FJ$-complex $f$ is a transfinite composition of
$\xymatrix@1{Y\co\lambda \ar[r] & \bfD}$, then $Gf$ is the
transfinite composition of $X=G \circ Y$. Since $Gf$ is a
transfinite composition of weak equivalences, $Gf$ is also a weak
equivalence. Hence $G$ takes relative $FJ$-complexes to weak
equivalences in $\bfC$.
\end{pf}

\subsection{Transfer of the Diagram Thomason Structure on
$\mathbf{Cat}^{\Delta^{op}}$}\label{subsection:Thomasontransfer}

With these preliminaries and our free constructions on double
categories, we can transfer the diagram Thomason structure to {\bf
DblCat}. Recall the diagram Thomason structure on
$\mathbf{Cat}^{\Delta^{op}}$ from Section \ref{diagramstructures}.

\begin{thm} \label{diagramtransfer}
There is a cofibrantly generated model structure on {\bf DblCat}
such that a double functor $K$ is a weak equivalence (respectively
fibration) if and only if $N_hK$ is levelwise a weak equivalence
(respectively fibration) in the Thomason structure on {\bf Cat}.
\end{thm}
\begin{pf}
We apply \ref{KanCorollary} to the adjunction $F=c_h\dashv N_h=G$.
First we point out that
\begin{equation*}
\aligned c_h\Big(\coprod_{\Delta^{op}([n],-)}c\Sd^2\Lambda^k[m]\Big) &=c_h(c\Sd^2\Lambda^k[m]\times \Delta[n])\\
&=(c\Sd^2\Lambda^k[m]) \boxtimes c\Delta[n] \\
&=(c\Sd^2\Lambda^k[m]) \boxtimes [n]
\endaligned
\end{equation*}
by Example \ref{horcatcatssetfromformula} or Proposition
\ref{horcatcatsset} (for simplicity we suppress $\sigma$ and $\nu$).
Similarly,
\begin{equation*}
c_h\Big(\coprod_{\Delta^{op}([n],-)}c\Sd^2\Delta[m]\Big)=(c\Sd^2\Delta[m])
\boxtimes [n]
\end{equation*}
and the horizontal categorification of the generating acyclic
cofibrations $j$ in Theorem \ref{diagramgenerators} are the
inclusions $f \boxtimes 1_{[n]}$ for the inclusions
$$\xymatrix@1{f\co c\Sd^2\Lambda^k[m] \ar[r] & c\Sd^2\Delta[m]}$$
and $n\geq 0$.
\begin{enumerate}
\item
The double categories $(c\Sd^2\partial\Delta[m]) \boxtimes [n]$ and
$(c\Sd^2\Lambda^k[m]) \boxtimes [n]$ have a finite number of
squares, hence they are finite by Proposition
\ref{doublecategoriesaresmall}.
\item
A transfinite composition of weak equivalences in $\mathbf{Cat}^{\Delta^{op}}$ is a weak equivalence
by Example \ref{simplicialtransfinitecompositions}.
\item
The horizontal nerve $N_h$ preserves filtered colimits by Theorem \ref{horizontalnervepreservesfilteredcolimits}.
\item
Consider the pushout in {\bf DblCat},
$$\xymatrix{(c\Sd^2\Lambda^k[m]) \boxtimes  [n] \ar[r]
\ar[d]_{c_h(j)=f\boxtimes 1_{[n]}} & \bbD \ar[d]^{j'} \\
(c\Sd^2\Delta[m]) \boxtimes [n] \ar[r] & \mathbb{P}\rlap{\,.}}$$
Then by Proposition \ref{horizontalnerveexternalproduct} $N_hc_h(j)$
is the acyclic cofibration $j$, and by Theorem \ref{pushoutiso} the
diagram
$$
\xymatrix@R=3pc{ (c\Sd^2\Lambda^k[m]) \times \Delta[n] \ar[r]
\ar[d]_{N_hc_h(j)=j} & N_h\bbD \ar[d]^{N_h(j')}  \\
(c\Sd^2\Delta[m]) \times \Delta[n] \ar[r] & N_h\bbP}
$$
is a pushout in $\mathbf{Cat}^{\Delta^{op}}$. Hence $N_h(j')$ is an
acyclic cofibration, and in particular a weak equivalence in
$\mbox{\bf Cat}^{\Delta^{op}}$.
\end{enumerate}
\end{pf}

We may compare the transferred diagram Thomason model structure on
$\mathbf{DblCat}$ to the Thomason model structure on $\mathbf{Cat}$
as follows.

\begin{prop} \label{ThomasonRightQuillenFunctor}
The functor $\xymatrix@1{(\bfV\text{-})_0 \co \mathbf{DblCat} \ar[r]
& \mathbf{Cat}}$ maps the weak equivalences and fibrations of the
transferred diagram Thomason model structure on $\mathbf{DblCat}$ to
weak equivalences and fibrations in the Thomason model structure on
$\mathbf{Cat}$. In particular, $(\bfV\text{-})_0$ is a right Quillen
functor.
\end{prop}
\begin{pf}
As functors, $(\bfV\text{-})_0$ is the same as $(N_h\text{-})_0$.
\end{pf}

\begin{cor} \label{ThomasonQuillenAdjunction}
The adjunction $$\xymatrix@C=4pc{\mathbf{Cat} \ar@{}[r]|{\perp}
\ar@/^1pc/[r]^-{\bbV} & \ar@/^1pc/[l]^-{(\bfV-)_0}
\mathbf{DblCat}}$$ is a Quillen adjunction.
\end{cor}
\begin{pf}
An adjunction is a Quillen adjunction if and only if the right
adjoint is a right Quillen functor, so the Corollary follows
immediately from Proposition \ref{ThomasonRightQuillenFunctor}.
\end{pf}

\begin{prop} \label{ThomasonExtension}
The functor $\xymatrix@1{\bbV\co \mathbf{Cat} \ar[r] &
\mathbf{DblCat}}$ preserves and reflects weak equivalences,
fibrations, and cofibrations. In other words, a functor $F$ is a
weak equivalence (respectively fibration, respectively cofibration)
in the Thomason model structure on $\mathbf{Cat}$ if and only if
$\bbV F$ is a weak equivalence (respectively fibration, respectively
cofibration) in the transferred diagram Thomason model structure on
$\mathbf{DblCat}$. As a consequence, the transferred diagram
Thomason model structure on $\mathbf{DblCat}$ extends the Thomason
model structure on $\mathbf{Cat}$ as a vertically embedded
subcategory.
\end{prop}
\begin{pf}
For a functor $F$, the morphism $N_h\bbV F$ of simplicial objects in
$\mathbf{Cat}$ is $F$ in every degree. Thus $\bbV$ preserves and
reflects weak equivalences and fibrations.

By Corollary \ref{ThomasonQuillenAdjunction}, $\bbV$ preserves
cofibrations. It also reflects cofibrations as follows. If $F$ is a
functor such that $\bbV F$ is a cofibration, and $G$ is an acyclic
fibration in $\mathbf{Cat}$, then any diagram in $\mathbf{DblCat}$
$$\xymatrix@1{\bbV\bfC \ar[d]_{\bbV F} \ar[r]^{K} & \bbV\mathbf{C'} \ar[d]^{\bbV G} \\
\bbV\bfD \ar[r]_{L} & \bbV\mathbf{D'}}$$ admits a lift, as $\bbV G$
is an acyclic fibration by the above. Since $\bbV$ is fully
faithful, this lift gives us a lift in $\mathbf{Cat}$. Hence $F$ is
a cofibration.
\end{pf}

\subsection{Transfer of the Diagram Categorical Structure on $\mbox{\bf Cat}^{\Delta^{op}}$} \label{subsection:categoricaltransfer}

Our preparations allow us to also quickly transfer the diagram
categorical structure. Recall the diagram categorical structure on
$\mathbf{Cat}^{\Delta^{op}}$ from Section \ref{diagramstructures}.
In Section \ref{subsection:simpliciallysurjective} we will show that
the vertical analogue of this model structure on {\bf DblCat}
coincides with the model structure induced by the simplicially
surjective topology $\tau$ on {\bf Cat} using the methods of
\cite{everaertinternal}. An important reason for interest in the
equality of these two structures lies in the fact that the second
construction yields an explicit form for the cofibrant replacement,
which is not at all transparent using only the transferred
structure.

\begin{thm} \label{categoricaltransfer}
There is a cofibrantly generated model structure on {\bf DblCat}
such that a double functor $K$ is a weak equivalence (respectively
fibration) if and only if $N_hK$ is levelwise a weak equivalence
(respectively fibration) in the categorical structure on {\bf Cat}.
\end{thm}
\begin{pf}
We apply \ref{KanCorollary} to the adjunction $F=c_h\dashv N_h=G$.
All generating acyclic cofibrations $j$ in Theorem
\ref{diagramgenerators} for the categorical diagram structure on
$\mathbf{Cat}^{\Delta^{op}}$ are natural transformations of the form
\begin{equation*}
\coprod_{\Delta^{op}([n],-)}\{1\}
\overset{\underset{\Delta^{op}([n],-)}{\coprod}
f}{\xymatrix@1@C=4pc{\ar[r] &}} \coprod_{\Delta^{op}([n],-)} \bfI
\end{equation*}
where $f$ is the inclusion $\xymatrix@1{\{1\} \ar[r] & \bfI}$ and
$[n]$ is an object of $\Delta^{op}$. These generating acyclic
cofibrations have horizontal categorification
$$\xymatrix@1{\{1\} \boxtimes [n] \ar[r] & \bfI \boxtimes [n]}$$
by Example \ref{horcatcatssetfromformula} or Proposition \ref{horcatcatsset} (for simplicity we suppress $\sigma$ and $\nu$).
\begin{enumerate}
\item
The double categories $\emptyset\boxtimes[n]$,
$\{0,1\}\boxtimes[n]$, $\{0\rightrightarrows1\}\boxtimes[n]$, and
$\{1\}\boxtimes[n]$ have a finite number of squares, hence they are
finite by Proposition \ref{doublecategoriesaresmall}.
\item
A transfinite composition of weak equivalences in $\mbox{\bf Cat}^{\Delta^{op}}$ is a weak equivalence
by Example \ref{simplicialtransfinitecompositions}.
\item
The horizontal nerve $N_h$ preserves filtered colimits by Theorem \ref{horizontalnervepreservesfilteredcolimits}.
\item
Consider the pushout in {\bf DblCat},
$$\xymatrix{\{1\} \boxtimes [n] \ar[r]
\ar[d]_{c_h(j)} & \bbD \ar[d]^{j'} \\
\bfI \boxtimes [n] \ar[r] & \mathbb{P}\rlap{\,.}}$$ Then by
Proposition \ref{horizontalnerveexternalproduct} $N_hc_h(j)$ is the
acyclic cofibration $j$, and by Theorem \ref{pushoutiso} the diagram
$$\xymatrix@R=3pc{\{1\} \times  \Delta[n] \ar[r]
\ar[d]_{N_hc_h(j)=j} & N_h\bbD  \ar[d]^{N_h(j')} \\
\bfI \times \Delta[n] \ar[r] & N_h\bbP}$$ is a pushout in
$\mathbf{Cat}^{\Delta^{op}}$. Hence $N_h(j')$ is an acyclic
cofibration, and in particular a weak equivalence in
$\mathbf{Cat}^{\Delta^{op}}$.
\end{enumerate}
\end{pf}

We may compare the transferred diagram categorical model structure
on $\mathbf{DblCat}$ to the categorical model structure on
$\mathbf{Cat}$ as follows.

\begin{prop} \label{CategoricalRightQuillenFunctor}
The functor $\xymatrix@1{(\bfV\text{-})_0 \co \mathbf{DblCat} \ar[r]
& \mathbf{Cat}}$ maps the weak equivalences and fibrations of the
transferred diagram categorical model structure on $\mathbf{DblCat}$
to weak equivalences and fibrations in the categorical model
structure on $\mathbf{Cat}$. In particular, $(\bfV\text{-})_0$ is a
right Quillen functor.
\end{prop}
\begin{pf}
As functors, $(\bfV\text{-})_0$ is the same as $(N_h\text{-})_0$.
\end{pf}

\begin{cor} \label{CategoricalQuillenAdjunction}
The adjunction $$\xymatrix@C=4pc{\mathbf{Cat} \ar@{}[r]|{\perp}
\ar@/^1pc/[r]^-{\bbV} & \ar@/^1pc/[l]^-{(\bfV-)_0}
\mathbf{DblCat}}$$ is a Quillen adjunction.
\end{cor}
\begin{pf}
This follows immediately from Proposition
\ref{CategoricalRightQuillenFunctor}.
\end{pf}

\begin{prop} \label{CategoricalExtension}
The functor $\xymatrix@1{\bbV\co \mathbf{Cat} \ar[r] &
\mathbf{DblCat}}$ preserves and reflects weak equivalences,
fibrations, and cofibrations. In other words, a functor $F$ is a
weak equivalence (respectively fibration, respectively cofibration)
in the categorical model structure on $\mathbf{Cat}$ if and only if
$\bbV F$ is a weak equivalence (respectively fibration, respectively
cofibration) in the transferred diagram categorical model structure
on $\mathbf{DblCat}$. As a consequence, the transferred diagram
categorical model structure on $\mathbf{DblCat}$ extends the
categorical model structure on $\mathbf{Cat}$ as a vertically
embedded subcategory.
\end{prop}
\begin{pf}
The proof is completely analogous to the proof of Proposition
\ref{ThomasonExtension}.
\end{pf}

\subsection{No Transfer of the Reedy Categorical Structure on $\mbox{\bf
Cat}^{\Delta^{op}}$} \label{subsection:NoReedytransfer} In this
subsection we consider the category $\mbox{\bf Cat}^{\Delta^{op}}$
of simplicial objects in {\bf Cat} equipped with the Reedy model
structure associated with the categorical model structure on {\bf
Cat}. The weak equivalences in this Reedy model structure are the
levelwise equivalences of categories and the fibrations are the
Reedy fibrations. (For further details, see \cite{hirschhorn}.) In
this section we will show that it is impossible to transfer this
model structure to {\bf DblCat} via the adjunction $c_h \dashv N_h\
$, where $N_h$ is the horizontal nerve and $c_h$ is the horizontal
categorification. We will need the following theorem.
\begin{thm}[Theorem 1 in \cite{joyalstreetpullback}]
\label{joyalstreet} For a given functor $\xymatrix@1{G\co\bfB \ar[r]
& \bfC}$, the canonical comparison functor from the pullback of $F$
along $G$ to the pseudo pullback of $F$ along $G$ is an equivalence
of categories for all functors $\xymatrix@1{F\co\bfA \ar[r] & \bfC}$
if and only if $G$ is an isofibration.
\end{thm}
Now we turn to the objective of this subsection.

\begin{thm} \label{NoReedy}
There does not exist a model structure on {\bf DblCat} such that a
double functor $K$ is a weak equivalence (respectively fibration) if
and only if $N_hK$ is a weak equivalence (respectively fibration) in
the Reedy model structure on $\mathbf{Cat}^{\Delta^{op}}$ associated
to the categorical structure on {\bf Cat}.
\end{thm}
\begin{pf}
Suppose that such a transferred model structure on {\bf DblCat} does
exist. Then $(c_h,N_h)$ is a Quillen pair. Let $\bbD$ be a double
category and consider a Reedy fibrant replacement $\xymatrix@1{r\co
N_h\bbD \ar[r] & V_\bullet}$ in $\mathbf{Cat}^{\Delta^{op}}$, that
is, $V_\bullet$ is a Reedy fibrant object and $r$ is an acyclic
cofibration in the Reedy structure. Our strategy is to prove that
the existence of such a transferred model structure implies a false
statement, namely that (\ref{canonicalfunctor}) is an equivalence of
categories for every double category $\bbD$. We then exhibit a
double category for which (\ref{canonicalfunctor}) is not an
equivalence.

Since $V_\bullet$ is Reedy fibrant, the map $(d_0,d_1)$ is an
isofibration in $\mathbf{Cat}$, which implies that the composite
functors $d_i$ $$\xymatrix@C=3pc{V_1\ar[r]_-{(d_0,d_1)}
\ar@/^1.5pc/[rr]^{d_i} & V_0\times V_0 \ar[r]_-{\pi_i} & V_0}$$ are
themselves isofibrations (the two projections $\xymatrix@1{\pi_i\co
V_0\times V_0 \ar[r] & V_0}$ are clearly isofibrations).  By Theorem
\ref{joyalstreet} this implies that the canonical functor
\begin{equation} \label{isofibrationforreplacement}
\xymatrix@1{V_1\times_{V_0}V_1 \ar[r] &
\overset{\vspace{0ex}}{\underset{\vspace{1ex}}{
V_1\stackrel{\mbox{\scriptsize ps}}{\times}_{V_0}V_1}}}
\end{equation}
from the pullback to the pseudo pullback is an equivalence of
categories.

Next we similarly show that (\ref{isofibrationforcategorification})
is an equivalence of categories. By Definition
\ref{horizontalcategorification} $(c_hV)_0=V_0$, the set of objects
of $V_1$ is contained in the set of horizontal morphisms of $c_hV$,
and the set of morphisms of $V_1$ is contained in the set of squares
of $c_hV$. We claim that the functor
\begin{equation}\label{fibration}
\xymatrix@C=3pc{ (c_hV)_1\ar[r]^-{(d_0,d_1)} &
(c_hV)_0\times(c_hV)_0 }
\end{equation}
is also an isofibration. First note that $\xymatrix@1{V_1\ar[r] &
V_0\times V_0}$ is an isofibration precisely when any diagram in
$c_hV$ of the special form
$$
\xymatrix{ \ar[d]_{v}^{\|\wr} &\ar[d]_{\wr\|}^{w}
\\
\ar[r]_g&
}
$$
with $v,w$ isomorphisms in $V_0$ and $g$ an object of $V_1$, can be
filled with an isomorphism $\alpha$ in $V_1$ (a ``generating
vertically invertible square'' in $c_hV$),
$$
\xymatrix@C=3pc{ \ar[d]_{v}^{\|\wr} \ar[r]^{g'} \ar@{}[dr]|{\alpha}
& \ar[d]_{\wr\|}^{w}
\\
\ar[r]_g&\rlap{\quad.} }
$$
Now consider a general such diagram in $c_h(V)$. This has the form
\begin{equation}\label{tbf}
\xymatrix{ \ar[d]_{v}^{\|\wr} &                 & &
&\ar[d]_{\wr\|}^{w}
\\
\ar[r]_{g_1} & \ar[r]_{g_2} & \ar@{..}[r] & \ar[r]_{g_n}&\rlap{\quad,}
}
\end{equation}
where the bottom edge is an equivalence class of a path of
composable horizontal morphisms by Definition
\ref{horizontalcategorification}. We next insert vertical identity
morphisms and fill in the individual squares to obtain the following
allowable compatible arrangement
$$
\xymatrix@C=3pc{ \ar[d]_{v}^{\|\wr} \ar[r]^{g_1'}
\ar@{}[dr]|{\alpha_1} & \ar@{=}[d]\ar[r]^{g_2'}
\ar@{}[dr]|{\alpha_2} &
\ar@{=}[d]\ar@{..}[r]&\ar@{=}[d]\ar[r]^{g_n'} \ar@{}[dr]|{\alpha_n}
& \ar[d]_{\wr\|}^w
\\
\ar[r]_{g_1} & \ar[r]_{g_2} & \ar@{..}[r] & \ar[r]_{g_n}&\rlap{\quad.}}
$$
(this is possible because $\xymatrix@1{V_1\ar[r] & V_0\times V_0}$
is an isofibration). The equivalence class in $c_h(V)$ of this
allowable compatible arrangement gives the required filling for
(\ref{tbf}). So (\ref{fibration}) is indeed an isofibration.
Reasoning as for (\ref{isofibrationforreplacement}), this implies
that the following functor is an equivalence of categories.
\begin{equation} \label{isofibrationforcategorification}
\xymatrix{ (c_hV)_1\times_{(c_hV)_0}(c_hV)_1\ar[r] &
\overset{\vspace{0ex}}{ \underset{\vspace{.1ex}}{ (c_hV)_1
\stackrel{\mbox{\scriptsize ps}}{\times}_{(c_hV)_0} (c_hV)_1}}}
\end{equation}

We claim that the unit $\eta_V$ is a weak equivalence. First we
observe that every component of the counit
$\xymatrix@1{\varepsilon\co c_hN_h \ar@{=>}[r] &
\mbox{Id}_{\mbox{\scriptsize\bf DblCat}}}$ is an isomorphism of
double categories, since the nerve functor $N_h$ is fully faithful.
Also, one of the triangle identities states that
$N_h\varepsilon_\bbD\cdot\eta_{N_h\bbD}=\mbox{id}$, so
$\eta_{N_h\bbD}$ is an isomorphism. The naturality of the unit
$\eta$ therefore gives us a commutative diagram
$$
\xymatrix@C=2.5em{ N_h\bbD = N_hc_hN_h\bbD \ar[d]_{N_hc_hr} &
\ar[l]_-{\eta_{N_h\bbD}} N_h\bbD \ar[d]^r
\\
N_hc_hV & \ar[l]^-{\eta_V} V
}
$$
in which the morphism $\eta_{N_h\bbD}$ is a levelwise equivalence
(since it is an isomorphism), and the morphism $r$ is a levelwise
equivalence (by hypothesis). The morphism $N_hc_hr$ is one as well,
since $c_hr$ is a weak equivalence ($c_hr$ is an acyclic cofibration
since $c_h$ was assumed to be a left Quillen functor and $r$ is an
acyclic cofibration). By the 2-out-of-3 property, it follows that
$\eta_V$ is a levelwise equivalence of categories, as claimed.

Consider the following commutative diagram in {\bf Cat}.
$$
\xymatrix{ (N_h\bbD)_2\ar[d]_{r_2}^\wr
\ar[r]^-{\text{Segal}}_-{\cong} \ar@{}[dr]|{(D)}
 & (N_h\bbD)_1\times_{(N_h\bbD)_0}(N_h\bbD)_1 \ar[d]^{(r_1,r_1)} \ar@{}[dr]|{(C)} \ar[r]
 & (N_h\bbD)_1\stackrel{\mbox{\scriptsize
 ps}}{\times}_{(N_h\bbD)_0}(N_h\bbD)_1
 \ar[d]^{(r_1,r_1)}
\\
V_2\ar[d]_{(\eta_V)_2}^\wr \ar@{}[dr]|{(B)} \ar[r]
 & V_1\times_{V_0}V_1 \ar[r]^\sim \ar[d]^{((\eta_V)_1,(\eta_V)_1)} \ar@{}[dr]|{(A)}
 & V_1\stackrel{\mbox{\scriptsize ps}}{\times}_{V_0}V_1 \ar[d]^{((\eta_V)_1,(\eta_V)_1)}
\\
(N_hc_hV)_2\ar[r]_-{\text{Segal}}^-{\cong} &
(N_hc_hV)_1\times_{(N_hc_hV)_0}(N_hc_hV)_1\ar[r]^\sim
 & (N_hc_hV)_1\stackrel{\mbox{\scriptsize ps}}{\times}_{(N_hc_hV)_0}(N_hc_hV)_1
}
$$
Note that $(B)$ and $(D)$ commute by the definition of Segal maps,
while the commutativity of $(A)$ and $(C)$ follows from the
universal property of the pseudo pullbacks. The vertical functors
$r_2$ and $(\eta_V)_2$ are equivalences of categories, since $r$ and
$\eta_V$ are weak equivalences from above. The bottom edge of $(C)$
is an equivalence, since it is (\ref{isofibrationforreplacement}).
The bottom edge of $(A)$ is an equivalence as it is
(\ref{isofibrationforcategorification}) (recall that
$(N_h\bbE)_0=\bbE_0$ and $(N_h\bbE)_1=\bbE_1$ for any double
category $\bbE$).

We claim that the top edge of $(C)$ is an equivalence of categories.
Since $r_0$ and $r_1$
are equivalences, the vertical functor
$$\xymatrix{(r_1,r_1)\co(N_h\bbD)_1\stackrel{\mbox{\scriptsize ps}}{\times}_{(N_h\bbD)_0}(N_h\bbD)_1
\ar@<-.6ex>[r] & V_1\stackrel{\mbox{\scriptsize
ps}}{\times}_{V_0}V_1 }$$ is an equivalence. Moreover, since
$\eta_V$ is a levelwise equivalence, the 2-out-of-3 property and the
commutativity of $(A)$ imply that the functor
$$
\xymatrix{V_1\times_{V_0}V_1 \ar[r] & (N_hc_hV)_1
\times_{V_0}(N_hc_hV)_1}
$$ is an equivalence.
Also, the commutativity of $(B)$ and the 2-out-of-3 property imply
that $\xymatrix@1{V_2 \ar[r] & V_1\times_{V_0}V_1}$ is an
equivalence. The commutativity of $(D)$ then implies that
$$
\xymatrix{(N_h\bbD)_1\times_{(N_h\bbD)_0}(N_h\bbD)_1 \ar[r] &
V_1\times_{V_0}V_1}
$$
is an equivalence. Finally, the commutativity of $(C)$ implies that
the canonical map
\begin{equation} \label{isofibrationfornerve}
\xymatrix{(N_h\bbD)_1\times_{(N_h\bbD)_0}(N_h\bbD)_1 \ar[r] &
\underset{\hspace{0ex}}{(N_h\bbD)_1\stackrel{\mbox{\scriptsize
ps}}{\times}_{(N_h\bbD)_0}(N_h\bbD)_1}}
\end{equation}
is an equivalence of categories, as claimed.

The map (\ref{isofibrationfornerve}) is nothing but
\begin{equation} \label{canonicalfunctor}
\xymatrix{\bbD_1\times_{\bbD_0}\bbD_1 \ar[r] &
\underset{\hspace{0ex}}{\bbD_1\stackrel{\mbox{\scriptsize
ps}}{\times}_{\bbD_0}\bbD_1}.}
\end{equation}
The objects of
$\underset{\hspace{0ex}}{\bbD_1\stackrel{\mbox{\scriptsize
ps}}{\times}_{\bbD_0}\bbD_1}$ are diagrams of the form
\begin{equation} \label{objectsofpseudopullback}
\xymatrix{& \ar[r]^g \ar[d]^{\cong} & \\ \ar[r]_f & &}
\end{equation} and morphisms of $\underset{\hspace{0ex}}{\bbD_1\stackrel{\mbox{\scriptsize
ps}}{\times}_{\bbD_0}\bbD_1}$ are diagrams of the form
\begin{equation} \label{morphismsofpseudopullback}
\xymatrix{\ar[r]^f \ar[d] \ar@{}[dr]|\alpha & \ar[d] & \ar[l]_\cong \ar[r]^g \ar@{}[dr]|\beta \ar[d] & \ar[d] \\
\ar[r]_{f'} & & \ar[l]^\cong \ar[r]_{g'} &}
\end{equation}
where the middle square is a {\it commutative square of vertical
morphisms}. The objects and morphisms of
$\bbD_1\times_{\bbD_0}\bbD_1$ are those of
(\ref{objectsofpseudopullback}) and
(\ref{morphismsofpseudopullback}) where the isomorphisms are
identities. The canonical functor (\ref{canonicalfunctor}) is given
by this inclusion.

We now exhibit a double category $\bbD$ where the canonical functor
(\ref{canonicalfunctor}) is not an equivalence, which then implies
that our original assumption on the existence of a transferred Reedy
structure on {\bf DblCat} is false. Let $\bbD$ be the double
category with four distinct objects $A,B,C,D$ and only the following
nontrivial morphisms.
\begin{equation}
\xymatrix{& C \ar[r]^g \ar[d]^{\cong} & D \\ A \ar[r]_f & B &}
\end{equation}
There are no nontrivial squares. Suppose that the canonical functor
(\ref{canonicalfunctor}) is essentially surjective. Then there exist
objects $X,Y,Z$ and morphisms as in the following diagram
\begin{equation}
\xymatrix{A \ar[r]^f \ar[d] \ar@{}[dr]|\alpha & B \ar[d] & C
\ar[l]_\cong \ar[r]^g \ar@{}[dr]|\beta \ar[d] & D \ar[d]
\\ X \ar[r]_{f'} & Y & Y \ar@{=}[l] \ar[r]_{g'} & Z}
\end{equation}
with $\alpha$ and $\beta$ vertically invertible squares. However,
since all squares in $\bbD$ are trivial, we conclude that $B=Y=C$, a
contradiction. Hence, the canonical functor (\ref{canonicalfunctor})
is not essentially surjective and is not an equivalence.

We conclude that it is impossible to transfer the categorical Reedy
model structure on $\mbox{\bf Cat}^{\Delta^{op}}$ to {\bf DblCat}.
\end{pf}

\section{Model Structures Arising from Grothendieck Topologies} \label{section:topologies}

Until now we have considered model structures transferred from
$\mathbf{Cat}^{\Delta^{op}}$. But one can also view double
categories as internal categories, and for these a homotopy theory
has already been developed. Model structures on internal categories
in a category $\bfC$ satisfying certain hypotheses have been studied
by Everaert, Kieboom, and Van der Linden in \cite{everaertinternal}.
As they point out, there are various notions of internal equivalence
of internal categories. The notions full and faithful representably
make sense for internal functors as in Definition
\ref{fullyfaithful}, but notions of essential surjectivity depend on
a class of morphisms $\mathcal{E}$ in $\bfC$. If this class of
morphisms is the class $\mathcal{E}_\mathcal{T}$ of
$\mathcal{T}$-epimorphisms for a Grothendieck topology $\mathcal{T}$
on $\bfC$, then the internal equivalences are the weak equivalences
for a model structure on $\Cat(\bfC)$ in good cases. The classes
$\we(\mathcal{T})$, $\fib(\mathcal{T})$, and $\cof(\mathcal{T})$ are
defined in \cite{everaertinternal} so that the following theorem
holds. We will recall the classes $\we(\mathcal{T})$,
$\fib(\mathcal{T})$, and $\cof(\mathcal{T})$ below.

\begin{thm}[5.5 of \cite{everaertinternal}] \label{topologymodel}
Let $\bfC$ be a finitely complete category such that $\Cat(\bfC)$ is
finitely complete and finitely cocomplete and $\mathcal{T}$ is a
Grothendieck topology on $\bfC$. If the class $\we(\mathcal{T})$ of
$\mathcal{T}$-equivalences has the 2-out-of-3 property and $\bfC$
has enough $\mathcal{T}$-projectives, then
$$(\Cat(\bfC),\we(\mathcal{T}),
\fib(\mathcal{T}),\cof(\mathcal{T}))$$ is a model
category.\footnote{As Tim Van der Linden pointed out to us, the
factorizations in this model structure are functorial if there
exists a functor $\xymatrix@1{P\co\bfC \ar[r] & \bfC}$ and a natural
transformation $\xymatrix@1{\eta\co P \ar@{=>}[r] & Id_\bfC}$ such
that $P(C)$ is $\mathcal{T}$-projective and $\eta_C$ is a
$\mathcal{T}$-epimorphism for all objects $C$ of $\bfC$. This is the
case with the $\tau$-topology, the $\tau'$-topology, and the trivial
topology we consider in this paper.}
\end{thm}
We apply this theorem to $\bfC=\mathbf{Cat}$ for various
Grothendieck topologies in this section.  In Section
\ref{subsection:simpliciallysurjective} we show that the model
structure associated to the simplicially surjective topology $\tau$
is the same as the transferred diagram categorical structure. This
second construction using \cite{everaertinternal} is advantageous,
as it gives us more information about the model structures, such as
simple descriptions of cofibrations and cofibrant replacements. We
will show in Section \ref{2monadstructure} that the model structure
associated to the categorically surjective basis $\tau'$ in Section
\ref{subsection:categoricallysurjective} turns out to be the same as
the model structure on {\bf DblCat} viewed as a category of algebras
over a 2-monad. We also show that the trivial topology induces the
trivial model structure associated to the 2-category
$\mathbf{DblCat_h}$.

\subsection{Homotopy Theory of Internal Categories as in
\cite{everaertinternal}}
First we recall the notions and results of \cite{everaertinternal}
for the special case of internal categories in $\bfC=\mathbf{Cat}$.
\begin{defn}
Let $\xymatrix@1{\iiso\co\Cat(\mbox{\bf Cat}) \ar[r] &
\Grpd(\mbox{\bf Cat})}$ be the right adjoint to the inclusion
$\xymatrix@1{\Grpd(\mbox{\bf Cat}) \ar[r] & \Cat(\mbox{\bf Cat})}$.
For $\bbB\in\Cat(\mbox{\bf Cat})$, this means that
$\iiso(\mathbb{B})_1$ has objects the invertible horizontal
morphisms of $\mathbb{B}$ and morphisms the horizontally invertible
squares. The category $\iiso(\mathbb{B})_1$ is a subcategory of
$\bbB_1$. Composition is the vertical composition of squares.
\end{defn}

\begin{defn} \label{mappingpathobjectdefinition}
If $\xymatrix@1{F\co\mathbb{A} \ar[r] & \mathbb{B}}$ is a double
functor, then the {\it mapping path object} is the category
$(\bbP_F)_0$ defined as the pullback below,
$$
\xymatrix@C=3pc@R=3pc{ (\bbP_F)_0 \ar[r]^{\overline{F}_0}
\ar[d]_{\overline{t}} & \iiso(\bbB)_1 \ar[d]^{t}
\\
\bbA_0 \ar[r]_{F_0} & \bbB_0\rlap{\,.}
}
$$
The objects of $(\bbP_F)_0$ are $(a,\xymatrix@1{f\co b
\ar[r]^{\cong} & F_0a})$ for $a$ an object of $\bbA$ and $f$ a
horizontal isomorphism of $\bbB$. The morphisms are pairs
\begin{equation} \label{morphismsofp}
\left(\raisebox{1.7em}{\xymatrix{a\ar[d]_k
\\a'}},\raisebox{1.8em}{\xymatrix{b\ar[r]^-{\cong}
\ar[d]_j\ar@{}[dr]|\alpha & F_0a\ar[d]^{F_0k}\\b'\ar[r]_-\cong &
F_0a'}} \right)
\end{equation}
where $k$ is a vertical morphism in $\bbA$ and $\alpha$ is a
horizontally invertible square in $\bbB$. Composition in
$(\bbP_F)_0$ comes from the vertical composition in $\bbA$ and
$\bbB$. The functor $\xymatrix@1{t\co\iiso(\bbB)_1 \ar[r] & \bbB_0}$
is the target for horizontal composition, exactly as in Definition
\ref{defnofdoublecategory}.
\end{defn}

\begin{defn}
Let $\mathcal{T}$ be a topology on ${\bf Cat}$. We denote by
$Y_{\mathcal{T}}$ the composition of the Yoneda embedding $Y$ with
the sheafification functor.
$$\xymatrix{\mathbf{Cat} \ar[r]_-Y \ar@/^1.5pc/[rr]^{Y_{\mathcal{T}}}
& \mathbf{Set}^{\mathbf{Cat}^{op}} \ar[r] &
Sh(\mathbf{Cat},\mathcal{T})}$$ A functor $\xymatrix@1{p\co\bfE
\ar[r] & \bfB}$ is a {\it $\mathcal{T}$-epimorphism} if
$Y_{\mathcal{T}}(p)$ is an epimorphism. In this case, we often
simply say that $p$ is {\it $\mathcal{T}$-epi}. We denote the class
of $\mathcal{T}$-epimorphisms by $\mathcal{E}_{\mathcal{T}}$.
\end{defn}

To show that a functor is $\mathcal{T}$-epi, we will use the
following characterization of $\mathcal{T}$-epimorphisms.

\begin{prop}[Corollary III.7.5 and III.7.6 in \cite{maclanemoerdijk1992},
Proposition 2.12 in \cite{everaertinternal}]
\label{epicharacterization} Let $\mathcal{T}$ be a topology on a
small category. A morphism $\xymatrix@1{p\co E \ar[r] & B}$ is
$\mathcal{T}$-epi if and only if for every morphism
$\xymatrix@1{g\co X \ar[r] & B}$ there exists a covering sieve
$\{\xymatrix@1{f_i\co U_i \ar[r] & X}\}_i$ and a family of morphisms
$\{\xymatrix@1{u_i\co U_i \ar[r] & E}\}_i$ such that for every $i
\in I$ the diagram
\begin{equation} \label{epicharacterizationdiagram}
\xymatrix{U_i \ar[d]_{u_i} \ar[r]^{f_i} & X \ar[d]^g \\
E \ar[r]_p & B }
\end{equation}
commutes.
\end{prop}

\begin{rmk} \label{epicharacterizationremark}
Suppose $K$ is a basis for the topology $\mathcal{T}$ in Proposition
\ref{epicharacterization} and such $g$ and $p$ are given. Then there
exists a covering {\it sieve} $\{\xymatrix@1{f_i\co U_i \ar[r] &
X}\}_i$ in $\mathcal{T}$ and a family of morphisms
$\{\xymatrix@1{u_i\co U_i \ar[r] & E}\}_i$ making
(\ref{epicharacterizationdiagram}) commute if and only if there
exists a covering {\it family} $\{\xymatrix@1{g_j\co V_j \ar[r] &
X}\}_j$ in $K$ and a family of morphisms $\{\xymatrix@1{v_j\co V_j
\ar[r] & E}\}_j$ making (\ref{epicharacterizationdiagram}) commute.
Thus, in Proposition \ref{epicharacterization} one could
equivalently replace the phrase ``covering sieve'' by the phrase
``covering family in a given basis''.
\end{rmk}

\begin{pf}
A sieve $S$ is a covering sieve in the topology $\mathcal{T}$
generated by the basis $K$ if and only if it contains a covering
family $R$ from the basis $K$. Suppose such a covering sieve
$\{f_i\}_i$ with morphisms $\{u_i\}_i$ is given. Then this covering
sieve contains a covering family in $L$ for which
(\ref{epicharacterizationdiagram}) commutes. Conversely, given such
a covering family $\{g_j\}_j$ with morphisms $\{v_j\}_j$, we may
take the sieve $$\{g_j \circ w| w \text{ a morphism such that }g_j
\circ w\text{ exists}\}_j$$ generated by the family $\{v_j\}_j$.
Then the family $\{v_j \circ w\}$ makes
(\ref{epicharacterizationdiagram}) commute.
\end{pf}

\begin{examp} \label{retractionsareTepi}
Suppose $\mathcal{T}$ is a topology on a small category. If a
morphism $p$ admits a right inverse $q$, then $p$ is a
$\mathcal{T}$-epimorphism. To see this using Proposition
\ref{epicharacterization}, take any covering family $\{f_i\}$ of $X$
and the morphisms $u_i=qgf_i$.
\end{examp}

\begin{defn}  \label{fullyfaithful}
Let $\xymatrix@1{s\co \iiso(\bbB)_1 \ar[r] & \bbB_0}$ be the source
map for horizontal composition, as in Definition
\ref{defnofdoublecategory}. A double functor $\xymatrix@1{F\co\bbA
\ar[r] & \bbB}$ is {\it essentially $\mathcal{T}$-surjective} if the
functor
$$\xymatrix@C=3pc{(\bbP_F)_0 \ar[r]^-{\overline{F}_0} & \iiso(\bbB)_1 \ar[r]^-{s} & \bbB_0}$$
given by
$$s\circ\overline{F}_0(a,\xymatrix@1{b\ar[r]^-\cong & F_0a})=b,$$
and
$$
s\circ\overline{F}_0\left(
\raisebox{1.7em}{\xymatrix{a\ar[d]_k\\a'}},\raisebox{1.8em}{\xymatrix{b\ar[r]^-{\cong}
\ar[d]_j\ar@{}[dr]|\alpha & F_0a\ar[d]^{F_0k}\\b'\ar[r]_-\cong &
F_0a'}} \right)=\raisebox{1.8em}{\xymatrix{b\ar[d]_j\\b'}}
$$
is a $\mathcal{T}$-epimorphism. If $F$ is additionally {\it fully
faithful} in the sense of \cite{bungepare79}, \ie if
$$\xymatrix@C=3pc{\bbA_1 \ar[r]^{F_1} \ar[d]_{(s,t)} & \bbB_1
\ar[d]^{(s,t)}
\\ \bbA_0 \times \bbA_0 \ar[r]_{F_0 \times F_0} & \bbB_0 \times \bbB_0 }$$
is a pullback square in ${\bf Cat}$, then $F$ is called a {\it
$\mathcal{T}$-equivalence}. We denote the class of
$\mathcal{T}$-equivalences by $\we(\mathcal{T})$. Note that a double
functor $F$ is fully faithful if and only if its restrictions to the
two functors
$$\xymatrix@1{(\Obj \bbA , \Hor \bbA) \ar[r] & (\Obj \bbB , \Hor \bbB) }$$
$$\xymatrix@1{(\Ver \bbA , \Sq \bbA) \ar[r] & (\Ver \bbB , \Sq \bbB)}$$
are both fully faithful.
\end{defn}

\begin{rmk} \label{traditionalesssurj}
If $\bfA$ and $\bfB$ are 1-categories, then a functor
$\xymatrix@1{F\co\bfA \ar[r] & \bfB}$ is essentially surjective in
the usual sense if and only if $s \circ \overline{F}_0$ is
surjective. The functor $F$ is fully faithful in the sense of
\ref{fullyfaithful} if and only if it is fully faithful in the usual
sense. The notions of essential surjectivity and fully faithfulness
can be found in any standard reference on category theory, such as
Pages 19 and 115 of \cite{borceux1} or Pages 14, 15, and 93 of
\cite{maclaneworking}.
\end{rmk}

\begin{defn} \label{projectivedefinition}
Let $\mathcal{E}$ be a class of functors. We say that a category
$\bfP$ is {\it projective with respect to the class $\mathcal{E}$}
if for every functor $\xymatrix@1{G\co \bfQ \ar@{->>}[r] & \bfR}$ in
$\mathcal{E}$ and every functor $\xymatrix@1{H\co\bfP \ar[r] &
\bfR}$ there exists a functor $\xymatrix@1{F\co\bfP \ar[r] & \bfQ}$
such that $GF=H$,
$$\xymatrix@C=3pc@R=3pc{& \bfP \ar[d]^H \ar@{-->}[dl]_{\exists F} \\ \bfQ \ar@{>>}[r]_G
& \bfR\rlap{\,.}}$$ The double arrow head signifies that $G$ is in
$\mathcal{E}$. If $\mathcal{E}$ is the class
$\mathcal{E}_\mathcal{T}$ of $\mathcal{T}$-epi functors, then a
projective category $\bfP$ is called {\it
$\mathcal{T}$-projective}.\footnote{This was also called {\it
$\mathcal{E}_{\mathcal{T}}$-projective} in \cite{everaertinternal}.}
\end{defn}

\begin{defn}
We say that there are {\it enough $\mathcal{T}$-projectives} in
${\bf Cat}$ if for every category $\bfC$ there exists a
$\mathcal{T}$-projective category $\bfP$ and a $\mathcal{T}$-epi
functor $\xymatrix@1{\bfP \ar@{->>}[r] & \bfC}$.
\end{defn}

\begin{defn} \label{definitionoffibration}
A double functor $\xymatrix@1{F\co\bbE \ar[r] & \bbB}$ is a {\it
$\mathcal{T}$-fibration} if the induced morphism $(r_F)_0$ in the
diagram below is a $\mathcal{T}$-epimorphism,
\begin{equation}\label{fibrationdefinition}
\xymatrix{\iiso(\bbE)_1 \ar@{-->}[dr]|{(r_F)_0}
\ar@/^1pc/[drr]^{\mbox{\scriptsize\sf iso}\,(F)_1}
\ar@/_1pc/[ddr]_{t} & & \\ & (\bbP_F)_0
\ar[d]_{\overline{t}} \ar[r]^{\overline{F}_0} & \iiso(\bbB)_1 \ar[d]^{t} \\
& \bbE_0 \ar[r]_{F_0} & \bbB_0\rlap{\,.}}
\end{equation}
\end{defn}

\begin{rmk} \label{isofibration}
If $\bfE$ and $\bfB$ are 1-categories, then $(r_F)_0$ is surjective
if and only if $F$ is an isofibration. Recall from Section
\ref{structuresonCat} that a functor $\xymatrix@1{F\co\bfE \ar[r] &
\bfB}$ is said to be an {\it isofibration} if for any object $e$ of
$\bfE$ and any isomorphism $b \cong Fe$ in $\bfB$, there exists a
lift to an isomorphism $b' \cong e$ in $\bfE$.
\end{rmk}

\begin{examp}[Example 5.2 of \cite{everaertinternal}]
\label{everyobjectfibrant} If $\mathbf{1}$ denotes the terminal
double category, then the unique double functor $\xymatrix@1{\bbD
\ar[r] & \mathbf{1}}$ is a $\mathcal{T}$-fibration for every
topology $\mathcal{T}$ on $\mathbf{Cat}$. Hence, in the model
structure of Theorem \ref{topologymodel}, every object is fibrant.
\end{examp}
\begin{pf}
In diagram (\ref{fibrationdefinition}), the functor $\overline{t}$
is the identity, so that $(r_F)_0$ is simply $t$. A right inverse
assigns horizontal identities to objects, and horizontal identity
squares to vertical morphisms. By Example \ref{retractionsareTepi},
the functor $(r_F)_0$ is a $\mathcal{T}$-epimorphism.
\end{pf}

\begin{prop}[Proposition 5.6 of \cite{everaertinternal}] \label{acyclicfibration}
A double functor $\xymatrix@1{F\co\bbE \ar[r] & \bbB}$ is an acyclic
$\mathcal{T}$-fibration if and only if it is fully faithful and
$F_0$ is a $\mathcal{T}$-epi functor.
\end{prop}

\begin{defn} \label{definitionofcofibration}
A double functor is a {\it $\mathcal{T}$-cofibration} if it has the left lifting
property with respect to all acyclic $\mathcal{T}$-fibrations.
\end{defn}

\begin{prop}[Proposition 5.9 of \cite{everaertinternal}] \label{cofibration}
A double functor $\xymatrix@1{J\co\bbA \ar[r] & \bbX}$ is a
$\mathcal{T}$-cofibration if and only if $J_0$ has the left lifting
property with respect to all $\mathcal{T}$-epi functors.
\end{prop}

\begin{cor} \label{cofibrantobject}
A double category $\bbX$ is cofibrant in the $\mathcal{T}$-model
structure if and only if $\bbX_0$ is $\mathcal{T}$-projective.
\end{cor}
\begin{pf}
By Proposition \ref{cofibration}, $\bbX$ is cofibrant if and only if
for any $\mathcal{T}$-epi functor $G$ and any functor $H$, a lift
$\ell$
$$\xymatrix{\emptyset \ar[r] \ar[d] & \bfQ \ar@{>>}[d]^G
\\ \bbX_0 \ar[r]_H \ar@{-->}[ur]^{\ell} & \bfR}$$ exists, or
equivalently, if $\bbX_0$ is $\mathcal{T}$-projective.
\end{pf}

\begin{rmk} \label{cofibrantreplacementprocedure}
These results allow us to construct a {\it cofibrant replacement}
$\bbE$ for a double category $\bbB$ in the $\mathcal{T}$-model
structure. Let $\bbE_0$ be a $\mathcal{T}$-projective category and
$\xymatrix@1{K_0\co\bbE_0 \ar@{>>}[r] & \bbB_0}$ a
$\mathcal{T}$-epimorphism (we will explicitly give $\bbE_0$ and
$K_0$ in the $\tau$- and $\tau'$-structures in Remarks
\ref{taucofibrantreplacement} and \ref{tau'cofibrantreplacement}).
Let $\bbE_1$ be the following pullback in {\bf Cat},
$$\xymatrix@C=3pc{\bbE_1 \ar[r]^{K_1} \ar[d]_{(s,t)} & \bbB_1 \ar[d]^-{(s,t)}
\\ \bbE_0 \times \bbE_0 \ar[r]_{K_0 \times K_0} & \bbB_0 \times \bbB_0\rlap{\,.} }$$
Then the double graph $\bbE$ carries a unique double category
structure such that $K=(K_0,K_1)$ is a double functor by Lemma 5.14
of \cite{everaertinternal}. Since $K$ is fully faithful and $K_0$ is
$\mathcal{T}$-epi, $\xymatrix@1{K\co\bbE \ar[r] & \bbB}$ is an
acyclic $\mathcal{T}$-fibration by Proposition
\ref{acyclicfibration}. By Corollary \ref{cofibrantobject}, $\bbE$
is a cofibrant double category, and hence a cofibrant replacement
for $\bbB$ in the $\mathcal{T}$-model structure.
\end{rmk}

\begin{rmk} \label{topologycontainment}
We conclude from Proposition \ref{epicharacterization} that if
$\mathcal{T}' \subseteq \mathcal{T}$ are Grothendieck topologies,
then every $\mathcal{T}'$-epi\-morph\-ism is a
$\mathcal{T}$-epi\-morph\-ism. Thus we conclude from Definitions
\ref{fullyfaithful}, \ref{definitionoffibration}, and
\ref{definitionofcofibration} that if $\mathcal{T}' \subseteq
\mathcal{T}$ then $\we(\mathcal{T}') \subseteq \we(\mathcal{T})$,
$\fib(\mathcal{T}') \subseteq \fib(\mathcal{T})$, and
$\cof(\mathcal{T}') \supseteq \cof(\mathcal{T})$.\footnote{We thank
Joachim Kock for posing this question.}
\end{rmk}

\subsection{Model Structure from the Simplicially Surjective Basis} \label{subsection:simpliciallysurjective}
For a category $\bfC$, we write $\bfC_k$ for the $k$-th set of the
nerve $N\bfC$. Similarly for a functor $F$ we write $(NF)_k=F_k$. We
say that a functor $F$ is {\it simplicially surjective} if $F_k$ is
surjective for all $k\geq 0$. We prove that the associated topology
on {\bf Cat} induces a model structure on {\bf DblCat} which
coincides with the vertical analogue of the transferred diagram categorical structure of
Section \ref{subsection:categoricaltransfer}. This second
construction gives additional information about (the vertical analogue of) the transferred
diagram categorical structure, including an explicit form for the
cofibrant replacement functor.

\begin{lem}
For a category $\bfC$ define $$K(\bfC):=\{\{\xymatrix@1{F\co\bfD
\ar[r] & \bfC}\}|\text{ $F$ a simplicially surjective functor }\}.$$
Then $K$ is a basis for a Grothendieck topology $\tau$ on ${\bf
Cat}$.
\end{lem}
\begin{pf}
\begin{enumerate}
\item
If $F$ is an isomorphism, then $NF$ is an isomorphism and each $F_k$
is bijective.
\item
If $\{F\}\in K(\bfC)$ and $\xymatrix@1{G\co \mathbf{C'} \ar[r] &
\bfC}$ is any functor, consider the pullback $\xymatrix@1{\pi_2\co
\bfD \times_{\bfC} \mathbf{C'} \ar[r] & \mathbf{C'}}$ in ${\bf Cat}$
of $F$ along $G$. Since the nerve functor preserves limits, $N\pi_2$
is the pullback of $NF$ along $NG$. Then $N\pi_{2}$ is simplicially
surjective, since limits of simplicial sets are formed pointwise.
\item
If $G \circ F$ exists and $F_k$ and $G_k$ are surjective for all
$k\geq 0$, then clearly $G_k \circ F_k$ is surjective for all $k\geq
0$, and $\{G \circ F\}$ is a covering.
\end{enumerate}
\end{pf}

\begin{lem}
A functor $\xymatrix@1{p\co \bfE \ar[r] & \bfB}$ is $\tau$-epi for
the Grothendieck topology $\tau$ if and only if $p$ is simplicially
surjective.
\end{lem}
\begin{pf}
If $p$ is $\tau$-epi, then take $g=1_{\bfB}$ in Proposition
\ref{epicharacterization} with Remark
\ref{epicharacterizationremark} to obtain $pu=f$ for some covering
family $\{f\}$ in $K$. Then $f_k$ is surjective for all $k \geq 0$.
Hence $p$ is simplicially surjective.

If $p$ is simplicially surjective, then $\{p\}$ is a covering family
in $K$, and so is the pullback $\pi_2$ of $p$ along $g$. Applying
Proposition \ref{epicharacterization} with Remark
\ref{epicharacterizationremark} again, we see that $p$ is
$\tau$-epi.
\end{pf}

Recall that the objects of the $k$-th category
$((\bbA_0)_k,(\bbA_1)_k)=(N_v\bbA)_k$ of the vertical nerve are
composable strings of $k$ vertical morphisms, and the morphisms are
vertically composable strings of $k$ squares. The composition is
horizontal composition of vertical strings of squares. Fully
faithful double functors and $\tau$-equivalences have a useful
characterization in terms of the vertical nerve.

\begin{prop} \label{nervecharacterizationoffullyfaithful}
A double functor $\xymatrix@1{F\co \bbA \ar[r] & \bbB}$ is fully
faithful if and only for every $k \geq 0$ the functor
$$\xymatrix{((F_0)_k,(F_1)_k)\co ((\bbA_0)_k,(\bbA_1)_k) \ar[r] &
((\bbB_0)_k,(\bbB_1)_k)}$$ is fully faithful.
\end{prop}
\begin{pf}
Since the nerve functor preserves pullbacks, and pullbacks of
simplicial sets are formed pointwise, it follows from Definition
\ref{fullyfaithful} that $F$ is fully faithful if and only if each
$((F_0)_k,(F_1)_k)$ is fully faithful.
\end{pf}

\begin{prop} \label{tauequivalence}
A double functor $\xymatrix@1{F\co\bbA \ar[r] & \bbB}$ is a
$\tau$-equivalence if and only if for every $k \geq 0$ the functor
$$\xymatrix{((F_0)_k,(F_1)_k)\co((\bbA_0)_k,(\bbA_1)_k) \ar[r] &
((\bbB_0)_k,(\bbB_1)_k)}$$ is an equivalence of categories.
\end{prop}
\begin{pf}
The double functor $F$ is essentially $\tau$-surjective if and only
if $s \circ \overline{F}_0$ is $\tau$-epi. But this occurs if and
only if $(s \circ \overline{F}_0)_k$ is surjective for each $k$,
which is equivalent to the essential surjectivity of
$((F_0)_k,(F_1)_k)$ by Remark \ref{traditionalesssurj}. Fully
faithfullness follows from Proposition
\ref{nervecharacterizationoffullyfaithful}.
\end{pf}

\begin{cor} \label{tau23}
The class $\we(\tau)$ of $\tau$-equivalences has the 2-out-of-3
property.
\end{cor}

\begin{prop} \label{tauprojectives}
${\bf Cat}$ has enough $\tau$-projectives.
\end{prop}
\begin{pf}
We first construct a $\tau$-projective category $\bfP$ from a
category $\bfC$. Let
\begin{equation*}
\bfP:=\coprod_{n \geq 0} \bfC_n \cdot [n]=\left( \coprod_{c \in
\bfC_0} [0]\right)\coprod \left(\coprod_{n \geq 1} \coprod_{(f_1,
\dots, f_n) \in \bfC_n} [n]\right)
\end{equation*}
where $[n]=\{0,1,\dots,n\}$ is the $(n+1)$-element ordinal viewed as
a category and $\bfC_n \cdot [n]$ denotes the copower of the
category $[n]$ with the set $\bfC_n$, as recalled in Remark
\ref{recallingcopower}. Suppose we have functors
$$
\xymatrix{
& \bfP \ar[d]^H
\\
\bfQ \ar@{->>}[r]_G & \bfR  \rlap{\,,} }
$$ and $G$ is $\tau$-epi. We denote $H$ on the $(f_1, \dots,
f_n)$-summand of $\bfP$ by
$$\xymatrix{H_{(f_1, \dots, f_n)}\co [n] \ar[r] & \bfR}$$
and by $H_c$ on the $c$-summand. If $H_{(f_1, \dots, f_n)}(j-1 \leq
j)=r_j$ for $1 \leq j \leq n$, then  there exists $(q_1, \dots, q_n)
\in \bfQ_n$ such that $G_n(q_1, \dots, q_n)=(r_1, \dots, r_n)$ since
$G$ is $\tau$-epi. We define a functor
$$\xymatrix{F_{(f_1, \dots, f_n)}\co [n] \ar[r] & \bfQ}$$
$$F_{(f_1, \dots, f_n)}(j-1 \leq j):=q_j$$
for $1 \leq j \leq n$. Similarly, for $c\in\bfC_0$, there exists $d
\in \bfQ_0$ such that $G_0(d)=H_c(0)$. We define $F_c(0)=d$. Putting
these $F$'s together, we obtain a functor $\xymatrix@1{F\co \bfP
\ar[r] & \bfQ}$ such that $GF = H$, and we conclude that $\bfP$ is
$\tau$-projective.

Next we construct a $\tau$-epi functor $\xymatrix@1{L\co \bfP
\ar@{->>}[r] & \bfC}$. On the $(f_1, \dots, f_n)$-summand of $\bfP$
define $L$ as
$$L_{(f_1, \dots, f_n)}(j-1 \leq j):=f_j$$ for $1\leq j \leq n$.
Similarly on the $c$ summand, $L_c(0):=c$. We claim that for each $k
\geq 0$, $\xymatrix@1{L_k\co \bfP_k \ar[r] & \bfC_k}$ is surjective.
Note that
\begin{equation*}
\bfP=\left( \coprod_{c \in \bfC_0} [0]_k\right)\coprod
\left(\coprod_{n \geq 1} \coprod_{(f_1, \dots, f_n) \in \bfC_n}
[n]_k\right).
\end{equation*}
If $k\geq1$, and $(f_1, \dots, f_k) \in \bfC_k$, then $L_k$ maps $(0
\leq 1, 1 \leq 2, \dots, k-1\leq k)$ in the $(f_1, \dots,
f_k)$-component of $\bfP_k$ to $(f_1, \dots, f_k)$. Similarly if $c
\in \bfC_0$, then $L_0$ maps $0$ in the $c$-component of $\bfP_0$ to
$c$. Hence $L_k$ is surjective for all $k \geq 0$, $L$ is
$\tau$-epi, and {\bf Cat} has enough $\tau$-projectives.
\end{pf}

\begin{thm}
The simplicially surjective topology $\tau$ on {\bf Cat} determines
a model structure
$$(\Cat(\mathbf{Cat}),\we(\tau),\fib(\tau),\cof(\tau)).$$
\end{thm}
\begin{pf}
The category $\Cat(\mathbf{Cat})$ is complete and cocomplete by
Theorem \ref{complete}. The class of $\tau$-equivalences has the
2-out-of-3 property by Corollary \ref{tau23} and {\bf Cat} has
enough $\tau$-projectives by Proposition \ref{tauprojectives}, so we
can apply Theorem \ref{topologymodel}.
\end{pf}

We now give a more explicit description of the fibrations, acyclic
fibrations, cofibrant objects, and fibrant objects.

\begin{prop} \label{taufibrationexplicitly}
Let $\xymatrix@1{F\co\bbE \ar[r] & \bbB}$ be a double functor.
\begin{enumerate}
\item
$F$ is a $\tau$-fibration if and only if for each $k\geq0$ the
functor
$$\xymatrix{((F_0)_k,(F_1)_k)\co ((\bbE_0)_k,(\bbE_1)_k) \ar[r] &
((\bbB_0)_k,(\bbB_1)_k)}$$ is an isofibration.
\item
$F$ is an acyclic $\tau$-fibration if and only if for each $k\geq0$
the functor $((F_0)_k,(F_1)_k)$ is fully faithful and surjective on
objects.
\end{enumerate}
\end{prop}

\begin{pf}
\begin{enumerate}
\item
Applying the nerve to Diagram (\ref{fibrationdefinition}), we see
that $F$ is a $\tau$-fibration if and only if $(r_F)_{0k}$ is
surjective for all $k \geq 0$. By Remark \ref{isofibration}, this is
the case if and only if for each $k \geq 0$ the functor
$((F_0)_k,(F_1)_k)$ is an isofibration. Here
$(\iiso(\bbB)_1)_k=\iiso((\bbB_0)_k,(\bbB_1)_k)$ is the category
with objects composable strings of $k$ vertical morphisms and with
morphisms vertical strings of vertically composable squares that are
each horizontally invertible.
\item
From Proposition \ref{nervecharacterizationoffullyfaithful}, $F$ is
fully faithful if and only if each $((F_0)_k,(F_1)_k)$ is fully
faithful. Since $F_0$ is $\tau$-epi if and only if $(F_0)_k$ is
surjective for each $k\geq 0$, the statement follows from
Proposition \ref{acyclicfibration}.
\end{enumerate}
\end{pf}

\begin{cor} \label{simpliciallysurjectivecoincidence}
The model structure on {\bf DblCat} induced by the simplicially
surjective topology $\tau$ on {\bf Cat} coincides with the model
structure obtained by transferring the diagram categorical structure
across the vertical categorification-vertical nerve adjunction $c_v
\dashv N_v$. The transfer across $c_v \dashv N_v$ is completely
analogous to the transfer across $c_h \dashv N_h$ in Section
\ref{subsection:categoricaltransfer}.
\end{cor}
\begin{pf}
From Propositions \ref{tauequivalence} and
\ref{taufibrationexplicitly} we see that the weak equivalences and
fibrations of the two model structures coincide.
\end{pf}

\begin{rmk} \label{taucofibrantreplacement}
We can now easily construct a cofibrant replacement $\bbE$ for a
double category $\bbB$ in the $\tau$-model structure. Let $\bbE_0$
be the $\tau$-projective category associated to $\bbB_0$ with
projection $K_0:=L$ as in the proof of Proposition
\ref{tauprojectives}. Then $\bbE$ and $K$ as defined in Remark
\ref{cofibrantreplacementprocedure} are a cofibrant replacement for
$\bbB$ in the $\tau$-model structure.
\end{rmk}

\begin{prop}
Let $F$ be a $\tau$-equivalence. Then $BF$, as in Definition
\ref{classifyingspacefunctor}, is a weak homotopy equivalence of
topological spaces.
\end{prop}
\begin{pf}
By Proposition \ref{tauequivalence}, $F$ is a $\tau$-equivalence if
and only if $(N_vF)_k$ is an equivalence of categories for each $k
\geq 0$. Since $((N_vF)_k)_{\ell}=(N_dF)_{k\ell}$, we see that
$(N_dF)_{k*}$ is a weak equivalence of simplicial sets for each
$k\geq0$. Hence $\diag(N_dF)$ is a weak equivalence of simplicial
sets, and $BF=|\diag(N_dF)|$ is a weak homotopy equivalence.
\end{pf}

\begin{rmk}
For each $m \in \mathbb{N}$, the assignment
$$\bfC \mapsto K_m(\bfC):=\{\{\xymatrix@1{F\co\bfD \ar[r] & \bfC}\}|
F_k \text{ surjective for all } 0\leq k \leq m \}$$ is a basis for a
Grothendieck topology $\tau_m$ on {\bf Cat}. We obtain a
$\tau_m$-model structure as above, though $\tau_m$-equivalences will
not necessarily be weak homotopy equivalences of classifying spaces.
\end{rmk}

\subsection{Model Structure from the Categorically Surjective Basis}
\label{subsection:categoricallysurjective} A functor is said to be
{\it categorically surjective} if it is surjective on objects and
full. It is straightforward to check that a basis for a Grothendieck
topology on {\bf Cat} is given by declaring a covering family to be a single
categorically surjective functor. We call this topology $\tau'$. In
this section we study the model structure on {\bf DblCat} induced by
$\tau'$. In Section \ref{2monadstructure} we show that this model
structure is the model structure on {\bf DblCat} viewed as a
category of algebras over a 2-monad.

As before we start with a characterization of the $\tau'$-epi
functors. We will use this to prove a 2-out-of-3 property for the
$\tau'$-equivalences.

\begin{prop} \label{tau'epis}
A functor $\xymatrix@1{p\co \bfE \ar[r] &\bfB}$ is $\tau'$-epi if
and only if there is a subcategory $\xymatrix@1{\bfH \,
\ar@{^{(}->}[r] & \bfE}$ such that $\xymatrix@1{p|_\bfH \co
\bfH\ar[r] & \bfB}$ is surjective on objects and full. Thus, a
$\tau'$-epi functor is not necessarily categorically surjective.
\end{prop}
\begin{pf}
Suppose that $p$ is $\tau'$-epi. Then by Proposition
\ref{epicharacterization} and Remark \ref{epicharacterizationremark}
there is a commutative square
$$
\xymatrix{ \bfU\ar[d]_u\ar[r]^{f} & \bfB\ar[d]^{1_{\bfB}}
\\
\bfE\ar[r]_p & \bfB\rlap{\,,} }
$$
where $f$ is surjective on objects and full. For each pair of
objects $x,y$ in $\bfU$, we have a commutative triangle
$$
\xymatrix{ \bfU(x,y)\ar[rr]^{f(x,y)}\ar[dr]_{u(x,y)} && \bfB(fx,fy)
\\
& \bfE(ux,uy) \ar[ur]_{p(ux,uy)}&. }
$$
Since $f(x,y)$ is surjective, so is $p(ux,uy)$. Let $\bfH=\im(u)$.
Thus, $p|_\bfH$ is surjective on objects and full.

Conversely, let $\xymatrix@1{\bfH\, \ar@{^{(}->}[r]^\ell & \bfE}$ be
a subcategory such that $p|_{\bfH}$ is surjective on objects and
full, and let $\xymatrix@1{g\co \bfX\ar[r] & \bfB}$ be any functor.
Consider the commutative diagram
$$
\xymatrix{
\bfH\times_{\bfB}\bfX\ar[r]^-{p'}\ar@{}[dr]|{\mbox{\scriptsize pb}}
\ar[d]_s & \bfX\ar[dd]^g
\\
\bfH\ar[d]_\ell\ar[dr]^{p|_{\bfH}} &
\\
\bfE \ar[r]_p & \bfB\rlap{\,.} }
$$
Then $p'$ is surjective on objects and full since $p|_\bfH$ is.
Further, $gp'=p\ell s$. By Proposition \ref{epicharacterization} and
Remark \ref{epicharacterizationremark}, it follows that $p$ is
$\tau'$-epi.
\end{pf}

Even though the $\tau'$-epi functors do not coincide with the
categorically surjective functors, they do give rise to the same
projective objects.

\begin{cor} \label{sameprojectives}
A category $\bfP$ is $\tau'$-projective if and only if it is
projective with respect to categorically surjective functors.
\end{cor}
\begin{pf}
We use the same notation as in Definition
\ref{projectivedefinition}. If $\bfP$ is $\tau'$-projective, then
$\bfP$ is projective with respect to categorically surjective
functors because every categorically surjective functor is
$\tau'$-projective by Proposition \ref{tau'epis}. For the converse,
suppose $\bfP$ is projective with respect to categorically
surjective functors, and suppose $G$ is $\tau'$-epi. Then by
Proposition \ref{tau'epis} again, there exists an inclusion
$\xymatrix@1{\ell\co \mathbf{Q'} \ar[r] & \bfQ}$ such that $G\ell$
is a categorically surjective functor. Thus there exists an $F'$
such that $G\ell F'=H$. If we let $F=\ell F'$ then we see that
$\bfP$ is $\tau'$-projective.
\end{pf}

Even better, we can characterize the $\tau'$-projective objects.

\begin{prop} \label{tauprimeprojectivecharacterization}
A category is $\tau'$-projective if and only if it is a free category on a directed graph.
\end{prop}
\begin{pf}
Let $\bfC$ be a free category on a directed graph $\Gamma$; we show
that $\bfC$ is $\tau'$-projective.

Suppose $\xymatrix@1{G \co \bfQ \ar@{->>}[r] &\bfR}$ is a
$\tau'$-epi and $\xymatrix@1{H\co\bfC \ar[r] & \bfR}$ is a functor.
Let $\mathbf{Q'}\subseteq\bfQ$ be a subcategory such that
$\xymatrix@1{G|_{\mathbf{Q'}} \co \mathbf{Q'} \ar[r] & \bfR}$ is
surjective on objects and full. Let $\xymatrix@1{U\co \mbox{\bf
Cat}\ar[r] & \mbox{\bf Graph}}$ denote the forgetful functor. Then
there is a map of graphs which makes the following diagram commute,
$$
\xymatrix@C=4pc{
& \Gamma \ar@{-->}[dl] \ar[d]^{U(H)|_{\Gamma}}\\
U\mathbf{Q'} \ar[r]_{U(G|_{\mathbf{Q'}})} & U\bfR}
$$
and induces a functor $F$ such that
$$
\xymatrix{ &\bfC\ar@{-->}[dl]_F \ar[d]^{H}
\\
\mathbf{Q'} \ar[r]_{G|_{\mathbf{Q'}}} & \bfR}
$$
commutes. Hence the diagram
$$
\xymatrix@R=.5pc{
    &&\bfC\ar[dl]\ar[dd]^H\\
    &\mathbf{Q'}\ar@{^{(}->}[dl] \\
\bfQ \ar@{^{(}->>}[rr]_G && \bfR
}
$$
commutes and $\bfC$ is $\tau'$-projective.

For the converse, suppose $\bfP$ is a $\tau'$-projective category.

Let $\bfQ$ be the free category on the underlying graph of $\bfP$ with the identity
arrows omitted.
Note that every arrow of $\bfQ$ is either a path $\langle f_1,\ldots,f_n\rangle$
of non-identity arrows $f_i$ in $\bfP$ or an empty path $\langle\rangle_A$ (forming the identity arrow $1_A$ on
the object $A$ in $\bfQ$).
We define an identity-on-objects functor $\xymatrix@1{G\colon \bfQ\ar[r] & \bfP}$
by $G(\langle\rangle_A)=1_A$ and $G(\langle f_1,\ldots,f_n\rangle)=f_n\circ\cdots\circ f_1$. The functor
$G$ is clearly categorically surjective.

Consider the diagram
$$
\xymatrix{
&\bfP\ar[d]^{1_\bfP}\ar@{-->}[dl]_F
\\
\bfQ\ar@{>>}[r]_-{G} & \bfP.
}
$$
Since $\bfP$ is $\tau'$-projective, there exists a functor $F$ which makes the diagram commute as indicated.
When $F(f)=\langle f_1,\ldots,f_n\rangle$, we say that $F(f)$ {\it has length $n$}.
Note that since $F$ is a functor, $F(1_A)$ has length $0$. By the commutativity of the diagram,
the length of $F(f)$ is greater than or equal to $1$ for any non-identity arrow $f$.

Now let $f$ be any non-identity arrow in $\bfP$, and suppose $F(f)=\langle f_1,\ldots,f_n\rangle$ as well as
$F(f_i)=\langle f_{i1},\ldots,f_{im_i}\rangle$. Since $GF =1_\bfP$, we
know that $f=f_n\circ\cdots\circ f_1$. Since $F$ is a functor, we also know that
$F(f)=F(f_n)\circ\cdots\circ F(f_1)=\langle f_{11},\ldots,f_{nm_n}\rangle$. Thus
$$\langle f_1,\ldots,f_n\rangle=\langle f_{11},\ldots,f_{nm_n}\rangle.$$
Since all the $f_i$ and $f_{ij}$ are non-identity arrows, it follows that
$m_i=1$ and $f_{i1}=f_i$ for all $i\in\{1,\ldots, n\}$. Summarizing,  if $F(f)=\langle f_1,\ldots,f_n\rangle$,
then $F(f_i)$ has length $1$ for each $i$.

Let $\bfC$ be the free category on the graph with the objects of $\bfP$ as vertices and as edges those arrows $f$ of $\bfP$
for which the length of $F(f)$ is $1$. By the argument above, the functor $F$ factors through a
functor $\xymatrix@1{\tilde{F}\colon \bfP\ar[r] & \bfC}$ as in the diagram below
$$
\xymatrix{
\bfC\ar[d]_{\mbox{\scriptsize incl}} & \ar[l]_-{\tilde F} \bfP \ar[d]^{1_\bfP}\ar@{-->}[dl]_F
\\ \bfQ\ar[r]_-{G} & \bfP.
}
$$
Let $\tilde{G}$ be the restriction of $G$ to $\bfC$. It is obvious that $\tilde{G}\tilde{F}=1_\bfP$.

Now consider the other composition, $\tilde{F}\tilde{G}$. This is obviously the identity on objects.
For a morphism $\langle f_1,\ldots,f_n\rangle$ in $\bfC$, we have
$$\aligned
\tilde{F}\tilde{G}(\langle f_1,\ldots,f_n\rangle) &=\tilde{F} (f_n\circ\cdots\circ f_1) \\
&=\tilde{F}(f_n)\circ\cdots\circ \tilde{F}(f_1) \\
&=\langle f_n\rangle\circ\cdots\circ\langle f_1\rangle \\
&=\langle f_1,\ldots,f_n\rangle,
\endaligned$$ where the second to last equality follows from the fact that the $f_i$ are edges in the graph on which $\bfC$
is free, in other words $\tilde{F}(f_i)$ has length $1$. So $\tilde{F}\tilde{G}(\langle f_1,\ldots,f_n\rangle)= \langle f_1,\ldots,f_n\rangle$. We conclude that $\bfP\cong \bfC$, so $\bfP$ is a free category on a directed graph.

\end{pf}

\begin{prop} \label{tauprimeepiimpliesnervesurjective}
If a functor $\xymatrix@1{p\co \bfA \ar[r] & \bfB}$ is $\tau'$-epi,
then for all $k \geq 0$, $\xymatrix@1{p_k\co \bfA_k \ar[r] &
\bfB_k}$ is surjective.
\end{prop}
\begin{pf}
By Proposition \ref{epicharacterization} and Remark
\ref{epicharacterizationremark} there exists a functor $f$
surjective on objects and full such that $1_\bfB \circ f=p \circ u$
for some $u$. Since $f_k=p_k \circ u_k$ is surjective, so is $p_k$.
\end{pf}

\begin{prop} \label{verticalnervesoftauprimequivalences}
If a double functor $\xymatrix@1{F\co \bbA \ar[r] & \bbB}$ is a
$\tau'$-equivalence, then for every $k \geq 0$ the functor
$$\xymatrix{((F_0)_k,(F_1)_k)\co ((\bbA_0)_k,(\bbA_1)_k) \ar[r] &
((\bbB_0)_k,(\bbB_1)_k)}$$ is an equivalence of categories.
\end{prop}
\begin{pf}
Since $F$ is fully faithful, $((F_0)_k,(F_1)_k)$ is fully faithful
for all $k \geq 0$ by Proposition
\ref{nervecharacterizationoffullyfaithful}. Since $F$ is essentially
$\tau'$-surjective, $s \circ \overline{F}_0$ is $\tau'$-epi and
hence $(s \circ \overline{F}_0)_k=(s)_k \circ (\overline{F}_0)_k$ is
surjective for all $k \geq 0$ by Proposition
\ref{tauprimeepiimpliesnervesurjective}. Remark
\ref{traditionalesssurj} then implies that $((F_0)_k,(F_1)_k)$ is
essentially surjective for all $k \geq 0$.
\end{pf}

\begin{lem}
Suppose $\xymatrix@1{\bbA\ar[r]^F & \bbB \ar[r]^G & \bbC}$ are
double functors and two of $GF,G,$ or $F$ are $\tau'$-equivalences.
Then the third double functor is fully faithful.
\end{lem}
\begin{pf}
By Proposition \ref{verticalnervesoftauprimequivalences} the
vertical nerves of the two $\tau'$-equivalences are levelwise
equivalences of categories. Hence the vertical nerve of the third
double functor is also levelwise an equivalence of categories, and
in particular levelwise fully faithful. By Proposition
\ref{nervecharacterizationoffullyfaithful}, this implies that the
third functor is fully faithful.
\end{pf}

\begin{lem}
Suppose $\xymatrix@1{\bbA\ar[r]^F & \bbB \ar[r]^G & \bbC}$ are
double functors and $GF$ and $F$ are $\tau'$-equivalences. Then $G$
is essentially $\tau'$-surjective.
\end{lem}
\begin{pf}
We need to show that $s\circ\overline{G}_0$ is $\tau'$-epi. Let
$\bfH_F\subseteq(\bbP_F)_0$ and $\bfH_{GF}\subseteq(\bbP_{GF})_0$ be
subcategories such that $s\circ\overline{F}_0|_{\bfH_F}$ and
$s\circ\overline{(GF)}_0|_{\bfH_{GF}}$ are surjective on objects and
full. Define a full subcategory $\bfH_G$ of
$(\bbP_G)_0=\bbB_0\times_{\bbC_0}\iiso(\bbC)_1$ by applying $F_0$ to
the first coordinate of $\bfH_{GF}$ as follows. For any object
$(a,c\stackrel{\cong}{\rightarrow}G_0F_0a)$ in $\bfH_{GF}$, we have
an object $(F_0a,c\stackrel{\cong}{\rightarrow}G_0(F_0a))$ in
$\bfH_G$. For any morphism
$$
\left(
\raisebox{1.7em}{\xymatrix{a\ar[d]_k\\a'}},\raisebox{1.8em}{\xymatrix{c\ar[r]^-{\cong}
\ar[d]_j\ar@{}[dr]|\alpha & G_0F_0a\ar[d]^{G_0F_0k}\\c'\ar[r]_-\cong
& G_0F_0a'}} \right) $$ in $\bfH_{GF}$ we have a morphism
$$
\left(
\raisebox{1.7em}{\xymatrix{F_0a\ar[d]_{F_0k}\\F_0a'}},\raisebox{1.8em}{\xymatrix{c\ar[r]^-{\cong}
\ar[d]_j\ar@{}[dr]|\alpha &
G_0(F_0a)\ar[d]^{G_0(F_0k)}\\c'\ar[r]_-\cong & G_0(F_0a')}} \right)
$$
in $\bfH_{G}$. Then we see as follows that
$s\circ\overline{G}_0|_{\bfH_G}\co\bfH_G\rightarrow\bbC_0$ is
surjective on objects and full. If $c\in\bbC_0$, there exists an
object $(a,c\stackrel{\cong}{\rightarrow}G_0F_0a)\in\bfH_{GF}$, with
$(F_0a,c\stackrel{\cong}{\rightarrow}G_0(F_0a))\in\bfH_G$ and
$s\overline{G}_0((F_0a,c\stackrel{\cong}{\rightarrow}G_0F_0a))=c$.
So $s\circ\overline{G}_0|_{\bfH_G}$ is surjective on objects. If
$c\stackrel{j}{\rightarrow}c'$ is a morphism in $\bbC_0$ and
$(F_0a,c\stackrel{\cong}{\rightarrow} G_0(F_0a))$ and
$(F_0a',c'\stackrel{\cong}{\rightarrow} G_0(F_0a'))$ are objects of
$\bfH_G$, then there exists a morphism
$$\left(
\raisebox{1.7em}{\xymatrix{
a\ar[d]_k\\a'}},\raisebox{1.8em}{\xymatrix{c\ar[r]^-{\cong}
\ar[d]_j\ar@{}[dr]|\alpha & G_0F_0a\ar[d]^{G_0F_0k}\\c'\ar[r]_-\cong
& G_0F_0a'}} \right)
$$
in $\bfH_{GF}$ which gives rise to a morphism
$$\left(
\raisebox{1.7em}{\xymatrix{
F_0 a\ar[d]_{F_0k}\\ F_0 a'}},\raisebox{1.8em}{\xymatrix{c\ar[r]^-{\cong}
\ar[d]_j\ar@{}[dr]|\alpha &
G_0(F_0a)\ar[d]^{G_0(F_0k)}\\c'\ar[r]_-\cong & G_0(F_0a')}} \right)
$$
in $\bfH_G$ that maps to $j$ under $s\circ\overline{G}_0|_{\bfH_G}$.
We conclude that $s\circ\overline{G}_0|_{\bfH_G}$ is surjective on
objects and full and therefore $s\circ\overline{G}_0$ is
$\tau'$-epi. This implies that $G$ is essentially
$\tau'$-surjective.
\end{pf}

\begin{lem}
Suppose $\xymatrix@1{\bbA\ar[r]^F&\bbB\ar[r]^G&\bbC}$ are double
functors and $GF$ and $G$ are $\tau'$-equivalences. Then $F$ is
essentially $\tau'$-surjective.
\end{lem}

\begin{pf}
We need to show that $s\circ\overline{F}_0$ is $\tau'$-epi. Let
$\bfH_{GH}\subseteq(\bbP_{GF})_0$ be a subcategory such that
$s\circ\overline{GF}_0|_{\bfH_{GF}}$ is surjective on objects and
full. Define a full subcategory $\bfH_F$ of
$(\bbP_F)_0=\bbA_0\times_{\bbB_0}\iiso(\bbB)_1$ with object set
$$
\Obj
\bfH_F:=\{(a,b\stackrel{\cong}{\rightarrow}F_0a)|(a,G_0(b\stackrel{\cong}{\rightarrow}F_0a))\in\Obj
\bfH_{GF}\}.
$$
Then we can see as follows that
$s\circ\overline{F}_0|_{\bfH_F}\co\bfH_F \rightarrow \bbB_0$ is
surjective on objects and full. If $b\in \bbB_0$, then
$G_0b\in\bbC_0$, and there is an object
$(a,G_0b\stackrel{\cong}{\rightarrow} G_0F_0a)\in\bfH_{GF}$, because
$s\circ(\overline{GF})_0|_{\bfH_{GF}}$ is surjective on objects.
However, $((G_0)_0,(G_1)_0)$ is fully faithful, \ie $G$ restricted
to the objects of $\bbB$ and the horizontal morphisms of $\bbB$ is a
fully faithful functor of categories. So there is a unique
isomorphism $b\stackrel{\cong}{\rightarrow} F_0a$ whose image under
$G$ is $G_0b\stackrel{\cong}{\rightarrow} G_0F_0a$. Hence,
$(a,b\stackrel{\cong}{\rightarrow}F_0a)\in\bfH_F$ and this object
maps to $b$ under $s\circ\overline{F}_0|_{\bfH_F}$. We conclude that
$s\circ\overline{F}_0|_{\bfH_F}$ is surjective on objects.

Moreover, if $b\stackrel{j}{\rightarrow}b'$ is a morphism in
$\bbB_0$ and $(a,b\stackrel{\cong}{\rightarrow}F_0a)$ and
$(a',b'\stackrel{\cong}{\rightarrow}F_0a')$ are objects of $\bfH_F$,
then $G_0j$ is a morphism of $\bbC_0$, and since
$s\circ(\overline{GF})_0|_{\bfH_{GF}}$ is full, there is a morphism
of the form
$$
\left( \raisebox{1.7em}{\xymatrix{
a\ar[d]_k\\a'}},\raisebox{1.8em}{\xymatrix{G_0b \ar[r]^-{G(\cong)}
\ar[d]_{G_0j}\ar@{}[dr]|\alpha &
G_0F_0a\ar[d]^{G_0F_0k}\\G_0b'\ar[r]_-{G(\cong)} & G_0F_0a'}}
\right)
$$
in $\bfH_{GF}$. However, the functor $((G_0)_1,(G_1)_1)$ is fully
faithful, so there is a unique square $\beta$, such that
$$
G\left(\raisebox{1.8em}{\xymatrix{b\ar[d]_j\ar[r]^-\cong\ar@{}[dr]|\beta & F_0a\ar[d]^{F_0k}\\
b'\ar[r]_-\cong & F_0a' }}\right)=\alpha.
$$
Moreover, $\beta$ is also horizontally invertible. Hence
$$
\left( \raisebox{1.7em}{\xymatrix{
a\ar[d]_k\\a'}},\raisebox{1.8em}{\xymatrix{b\ar[r]^-{\cong}
\ar[d]_{j}\ar@{}[dr]|\beta &F_0a\ar[d]^{F_0k}\\b'\ar[r]_-{\cong} &
F_0a'}} \right)
$$
is a morphism in $\bfH_F$ which maps to $j$ under
$s\circ\overline{F}_0$. We conclude that
$s\circ\overline{F}_0|_{\bfH_F}$ is surjective on objects and full,
so $s\circ\overline{F}_0$ is $\tau'$-epi and $F$ is essentially
$\tau'$-surjective.
\end{pf}

\begin{lem}
Suppose that $\xymatrix{\bbA\ar[r]^F&\bbB\ar[r]^G&\bbC}$ are double
functors, and $F$ and $G$ are $\tau'$-equivalences. Then $G\circ F$
is essentially $\tau'$-surjective.
\end{lem}

\begin{pf}
We need to show that $s\circ(\overline{GF})_0$ is $\tau'$-epi. Let
$\bfH_F\subseteq(\bbP_F)_0$ and $\bfH_G\subseteq(\bbP_G)_0$ be
subcategories such that $s\circ\overline{F}_0|_{\bfH_F}$ and
$s\circ\overline{G}_0|_{\bfH_G}$ are surjective on objects and full.
Let $\bfH_{GF}$ be the full subcategory of
$(\bbP_{GF})_0=\bbA_0\times_{\bbC_0}\iiso(\bbC)_1$, with objects
$$\aligned
\Obj\bfH_{GF}:=&\{(a,c\stackrel{\cong}{\rightarrow}G_0b\stackrel{\cong}{\rightarrow}G_0F_0a)\,|
\\
&(b,c\stackrel{\cong}{\rightarrow}G_0b)\in\bfH_G,(a,b\stackrel{\cong}{\rightarrow}F_0a)\in\bfH_F\}.
\endaligned
$$
Suppose that $c\in\bbC_0$, then there are objects
$(b,c\stackrel{\cong}{\rightarrow}G_0b)\in \bfH_G$ and
$(a,b\stackrel{\cong}{\rightarrow}F_0a)\in\bfH_F$, because $s\circ
\overline{G}_0|_{\bfH_G}$ and $s\circ \overline{F}_0|_{\bfH_F}$ are
surjective on objects. Thus
$(a,c\stackrel{\cong}{\rightarrow}G_0b\stackrel{\cong}{\rightarrow}G_0F_0a)\in\bfH_{GF}$,
and this object maps to $c$ under $s\circ
\overline{GF}_0|_{\bfH_{GF}}$.

Next, suppose that $c\stackrel{j}{\rightarrow}c'$ is a morphism of
$\bbC_0$ and
$(a,c\stackrel{\cong}{\rightarrow}G_0b\stackrel{\cong}{\rightarrow}G_0F_0a)$
and
$(a',c'\stackrel{\cong}{\rightarrow}G_0b'\stackrel{\cong}{\rightarrow}G_0F_0a')$
are objects of $\bfH_{GF}$. Then there exist morphisms
$$
\left(\raisebox{1.7em}{\xymatrix{b\ar[d]_{k_b}\\b'}},
\raisebox{1.8em}{\xymatrix{c\ar[d]_j\ar[r]^-\cong\ar@{}[dr]|\alpha & G_0b\ar[d]^{G_0k_b}\\
c'\ar[r]_-\cong & G_0b'}}\right)
$$
in $\bfH_G$, and
$$
\left(\raisebox{1.7em}{\xymatrix{a\ar[d]_{k_a}\\a'}},
\raisebox{1.8em}{\xymatrix{b\ar[d]_{k_b}\ar[r]^-\cong\ar@{}[dr]|\beta & F_0a\ar[d]^{F_0k_a}\\
b'\ar[r]_-\cong & F_0a'}}\right)
$$
in $\bfH_F$,
and therefore
$$
\left(\raisebox{1.7em}{\xymatrix{a\ar[d]_{k_a}\\a'}},\raisebox{1.8em}{\xymatrix{c\ar[d]_j\ar[r]^-\cong\ar@{}[dr]|\alpha
&
    G_0b\ar[d]|{G_0k_b}\ar@{}[dr]|{G\beta}\ar[r]^-\cong & G_0F_0a\ar[d]^{G_0F_0k_a}\\ c'\ar[r]_-\cong &
    G_0b'\ar[r]_-\cong & G_0F_0a'}}\right)
$$
is a morphism of $\bfH_{GF}$ that maps to
$c\stackrel{j}{\rightarrow}c'$ under
$s\circ(\overline{GF})_0|_{\bfH_{GF}}$. So we have proved that
$s\circ(\overline{GF})_0|_{\bfH_{GF}}$ is surjective on objects and
full. We conclude that $GF$ is $\tau'$-epi and essentially
$\tau'$-surjective.
\end{pf}

The previous four lemmas are summarized in the following theorem.

\begin{thm} \label{tau'23}
The class $\we(\tau')$ of $\tau'$-equivalences has the 2-out-of-3
property.
\end{thm}

\begin{prop} \label{tau'projectives}
{\bf Cat} has enough $\tau'$-projectives.
\end{prop}
\begin{pf}
Suppose $\bfC$ is a category.  Let $\bfP$ be the free category on the underlying
graph of $\bfC$. This is $\tau'$-projective by Proposition \ref{tauprimeprojectivecharacterization}.
The functor $\xymatrix{\bfP\ar[r] & \bfC}$, which
is the identity on objects and defined by composition on paths of
morphisms, is surjective on objects and full, so it is $\tau'$-epi.
Thus {\bf Cat} has enough $\tau'$-projectives.
\end{pf}

\begin{thm}
The categorically surjective topology $\tau'$ determines a model
structure
$$(\Cat(\mathbf{Cat}),\we(\tau'),\fib(\tau'),\cof(\tau')).$$
\end{thm}
\begin{pf}
The category $\Cat(\mathbf{Cat})$ is complete and cocomplete by
Theorem \ref{complete}. The class of
$\tau'$-equivalences has the 2-out-of-3 property by Corollary
\ref{tau'23} and {\bf Cat} has enough $\tau'$-projectives by
Proposition \ref{tau'projectives}, so we can apply Theorem
\ref{topologymodel}.
\end{pf}

\begin{prop} \label{tau'cofibrantobjects}
A double category $\bbX$ is cofibrant in the $\tau'$-model structure
if and only if $\bbX_0$ is a free category on a directed graph.
\end{prop}
\begin{pf}
By Corollary \ref{cofibrantobject} a double category is cofibrant if
and only if $\bbX_0$ is projective with respect to $\tau'$-epi
functors. But by Corollary \ref{tauprimeprojectivecharacterization}, $\bbX_0$ is
projective with respect to $\tau'$-epis if and only if it is
a free category on a directed graph.
\end{pf}

\begin{rmk} \label{tau'cofibrantreplacement}
We can now easily construct a cofibrant replacement $\bbE$ for a
double category $\bbB$ in the $\tau'$-model structure. Let $\bbE_0$
be the free category on the underlying graph of $\bbB_0$, and
$\xymatrix@1{K_0\co \bbE_0 \ar[r] & \bbB_0}$ the functor which is
the identity on objects and composition on paths of morphisms. Then
$\bbE_0$ is $\tau'$-projective, and $K_0$ is a $\tau'$-epimorphism
as in the proof of Proposition \ref{tau'projectives}. Then $\bbE$
and $K$ as defined in Remark \ref{cofibrantreplacementprocedure} are
a cofibrant replacement for $\bbB$ in the $\tau'$-model structure.
\end{rmk}

As an immediate consequence of Propositions \ref{tauequivalence} and
\ref{verticalnervesoftauprimequivalences}, we see that every
$\tau'$-equivalence is a $\tau$-equivalence. This also follows from
Remark \ref{topologycontainment}, since the categorically surjective
$\tau'$-topology is contained in the simplicially surjective
$\tau$-topology. An interesting question is whether or not a
condition slightly stronger than simplicial surjectivity but also
slightly weaker than categorical surjectivity would give rise to a
model structure with weak equivalences between those of the
$\tau'$-structure and the $\tau$-structure. For example, such a
condition on a functor is to be {\it $U$-split}. However, this
condition recovers the $\tau'$-topology instead of something new. In
fact, this condition only gives a different basis for the
$\tau'$-topology which will be of use in Section
\ref{2monadstructure}.
\begin{defn}
Let $\xymatrix@1{U\co \mathbf{Cat} \ar[r] & \mathbf{Graph}}$ be the
forgetful functor from categories to directed graphs. We say that a
functor $p$ is {\it $U$-split} if there exists a morphism $q$ of
directed graphs such that $(Up) \circ q=id$.
\end{defn}

\begin{lem} \label{Usplitcharacterization}
A functor $\xymatrix@1{p\co \bfE \ar[r] &\bfB}$ is $U$-split if and
only if there is a subcategory $\xymatrix@1{\bfH \, \ar@{^{(}->}[r]
& \bfE}$ such that $\xymatrix@1{p|_\bfH \co \bfH\ar[r] & \bfB}$ is
surjective on objects and full.
\end{lem}
\begin{pf}
Suppose $p$ is $U$-split. Then there exists a morphism of directed
graphs $q$ such that $Up \circ q=id$. Let $\bfH$ be the full
subcategory $\bfE$ whose objects are in the image of $q$. Then
$p|_\bfH$ is surjective on objects and full, as one sees using the
directed graph section $q$.

Conversely, suppose there exists a subcategory $\bfH$ of $\bfE$ such
that $p|_\bfH$ is categorically surjective. Then $p|_{\bfH}$ is
$U$-split, and $id=U(p|\bfH) \circ q=Up \circ q$ so that $p$ is also
$U$-split.
\end{pf}

\begin{prop}
The assignment
$$\bfC \mapsto L(\bfC):=\{\{\xymatrix@1{F\co\bfD \ar[r]
& \bfC}\}|\text{ $F$ is $U$-split } \}$$ is a basis for the
$\tau'$-topology on {\bf Cat}.
\end{prop}
\begin{pf}
We omit the proof that this is a basis.

Recall that a sieve is a covering sieve in the topology induced by a
basis if and only if it contains a covering family from the basis.
If $S$ is a $\tau'$-covering sieve, it contains a categorically
surjective functor, and hence a $U$-split functor by Lemma
\ref{Usplitcharacterization}, so that $S$ is also a covering sieve
in the topology induced by $L$.

Conversely, suppose $S$ is a sieve on $\bfB \in \mathbf{Cat}$ that
contains a $U$-split functor $\xymatrix@1{p\co \bfE \ar[r] & \bfB}$
and $p|_\bfH$ is categorically surjective. Then $p \circ i \in S$
for the inclusion $\xymatrix@1{i\co \bfH \ar[r] & \bfE}$, and $S$ is
a covering sieve in the $\tau'$-topology.
\end{pf}

\subsection{Model Structure from the Trivial Topology}\label{trivialmodelstructure}

On the underlying category of any 2-cat\-e\-go\-ry $\mathcal{K}$
with finite limits and finite colimits there is the {\it trivial
model structure} as proved in \cite{lack2monads} using pseudo
limits. A weak equivalence (respectively fibration) in this model
structure is a morphism $\xymatrix@1{f\co A \ar[r] & B}$ such that
$\xymatrix@1{\mathcal{K}(E,f)\co \mathcal{K}(E,A) \ar[r] &
\mathcal{K}(E,B)}$ is a weak equivalence (respectively fibration)
for all $E$ in the categorical model structure on {\bf Cat}. Thus
$f$ is a weak equivalence if and only if there is a morphism
$\xymatrix@1{g\co B \ar[r] & A}$ such that $gf$ and $fg$ are
isomorphic via 2-cells to the respective identities. A morphism $f$
is a fibration, or {\it isofibration}, if and only if for all
morphisms $\xymatrix@1{a\co E \ar[r] & A}$ and $\xymatrix@1{b\co E
\ar[r] & B}$ and any invertible 2-cell $\beta\co b \cong fa$, there
exists a morphism $\xymatrix@1{a'\co E \ar[r] & A}$ and an
invertible 2-cell $\alpha\co a' \cong a$ with $fa'=b$ and
$f\alpha=\beta$. If the 2-category $\mathcal{K}$ is merely a
1-category, then the trivial model structure agrees with the usual
trivial model structure: the weak equivalences are exactly the
isomorphisms and all morphisms are both fibrations and cofibrations.
The trivial model structure on the underlying category of a
2-category $\mathcal{K}$ is compatible with the
$\mathbf{Cat}$-enrichment as proved in \cite{lack2monads}.

Thus {\bf DblCat} admits three trivial model structures, depending
on whether we take as 2-cells the horizontal natural
transformations, the vertical natural transformations, or only
trivial 2-cells. When we say {\it trivial model structure on} {\bf
DblCat} we mean the one arising from the 2-category
$\mathbf{DblCat_h}=\Cat(\mathbf{Cat})$ which has horizontal natural
transformations as its 2-cells.

The following theorem summarizes Section 7 of
\cite{everaertinternal} and Theorem 3.3 of \cite{lack2monads}
applied to $\mathcal{K}=\Cat(\mathbf{C})$ to conclude that the
$\tau_{triv}$-model structure coincides with the trivial model
structure on the underlying category of the 2-category $\Cat(\bfC)$.

\begin{thm} \label{trivialtopologyproducestrivialstructure}
Let $\mathbf{C}$ be a finitely complete category such that
$\Cat(\mathbf{C})$ is finitely complete and finitely cocomplete. Let
$\tau_{triv}$ denote the trivial topology\footnote{In the trivial
topology the only covering sieve on an object is the maximal sieve.}
on $\mathbf{C}$. Then the following hold.
\begin{enumerate}
\item
A morphism $p$ in $\bfC$ is $\tau_{triv}$-epi if and only if there
exists a morphism $q$ such that $pq=id$.
\item
An internal functor $F$ is an acyclic $\tau_{triv}$-fibration if and
only if there exists an internal functor $G$ such that $FG=id$ and
$GF \cong id$.
\item
Every object of $\bfC$ is $\tau_{triv}$-projective, and hence $\bfC$
has enough projectives.
\item
Assume for the rest of this theorem that the
$\tau_{triv}$-equivalences have the 2-out-of-3 property. Then we
have the $\tau_{triv}$-model structure on $\Cat(\bfC)$ of Theorem
\ref{topologymodel}.
\item
Every object of $\Cat(\bfC)$ is fibrant and cofibrant.
\item
The weak equivalences in the $\tau_{triv}$-model structure are
precisely the equivalences in the 2-category $\Cat(\bfC)$.
\item
The $\tau_{triv}$-model structure coincides with the trivial model
structure of \cite{lack2monads} on the underlying category of the
2-category $\mathcal{K}=\Cat(\bfC)$.
\end{enumerate}
\end{thm}
\begin{pf}
\begin{enumerate}
\item
This follows from Proposition \ref{epicharacterization}.
\item
This is Proposition 7.1 of \cite{everaertinternal}, and it follows
from (i) along with Proposition \ref{acyclicfibration}.
\item
This follows immediately from (i).
\item
The model structure follows from (iii) and the hypotheses we have
made.
\item
Every object is fibrant by Example \ref{everyobjectfibrant}, and
every object is cofibrant by (ii).
\item
Since every object is fibrant and cofibrant, the weak equivalences
are precisely the homotopy equivalences. In Section 3 of
\cite{everaertinternal}, a cocylinder on $\Cat(\bfC)$ is given such
that two internal functors are homotopic if and only if they are
naturally isomorphic. Hence, the weak equivalences are precisely the
equivalences.
\item
By Theorem 3.3 of \cite{lack2monads}, the weak equivalences in the
trivial model structure on $\Cat(\bfC)$ are the equivalences. By 3.4
of \cite{lack2monads}, the acyclic fibrations in the trivial model
structure are the morphisms in (ii). Since the classes of weak
equivalences and acyclic fibrations in the two model structures are
the same, we conclude that the two model structures coincide.
\end{enumerate}
\end{pf}

\begin{rmk}
The assumption that the $\tau_{triv}$-equivalences have the
2-out-of-3 property can be removed by proving directly that the
$\tau_{triv}$-equivalences are the equivalences in the 2-category
$\Cat(\bfC)$ using (i) and fully faithfulness.
\end{rmk}

The trivial model structure on a category of internal categories is
much like the Str{\o}m structure of \cite{strom}. This analogy is
made precise in Section 7 of \cite{everaertinternal}.

We finish this section by comparing the categorical model structure
on $\mathbf{Cat}$ with the trivial model structure on
$\mathbf{DblCat}$ that arises from $\mathbf{DblCat_h}$. Our
comparison uses the horizontal embedding $\bbH$ of $\mathbf{Cat}$
into $\mathbf{DblCat}$, which is the same as the embedding
$$\xymatrix@1{\Cat(\mathbf{Set})\ar[r] & \Cat(\mathbf{Cat})}$$ induced
by the embedding $\xymatrix@1{\mathbf{Set} \ar[r] & \mathbf{Cat}}$.

\begin{prop} \label{TrivialRightQuillenFunctor}
The functor $\xymatrix@1{(\bfH\text{-})_0 \co \mathbf{DblCat} \ar[r]
& \mathbf{Cat}}$ maps the weak equivalences and fibrations of the
trivial model structure on $\mathbf{DblCat}$ to weak equivalences
and fibrations in the categorical model structure on $\mathbf{Cat}$.
In particular, $(\bfH\text{-})_0$ is a right Quillen functor.
\end{prop}
\begin{pf}
The functor $(\bfH\text{-})_0$ preserves weak equivalences because
it is the underlying functor of the 2-functor $\xymatrix@1{\bfH\co
\mathbf{DblCat_h} \ar[r] & \mathbf{Cat}}$, which maps equivalences
to equivalences.

To prove that $(\bfH\text{-})_0$ preserves fibrations, we use the
characterization of isofibrations (in a 2-category in general, and
in $\mathbf{Cat}$ in particular) at the beginning of Section
\ref{trivialmodelstructure}. If $F$ is a double functor that is an
is an isofibration, and $\xymatrix@1{\beta \co b \ar@{=>}[r] & (\bfH
F)_0 \circ a}$ is a natural transformation in $\mathbf{Cat}$, then
we obtain the required $\alpha$ by applying $\bbH$ to the lifting
problem, and then applying $(\bfH\text{-})_0$ to the solution
$\alpha'$ of the new lifting problem in $\mathbf{DblCat}$.
\end{pf}

\begin{cor} \label{TrivialQuillenAdjunction}
The adjunction $$\xymatrix@C=4pc{\mathbf{Cat} \ar@{}[r]|{\perp}
\ar@/^1pc/[r]^-{\bbH} & \ar@/^1pc/[l]^-{(\bfH-)_0}
\mathbf{DblCat}}$$ is a Quillen adjunction.
\end{cor}
\begin{pf}
This follows immediately from Proposition
\ref{TrivialRightQuillenFunctor}.
\end{pf}

\begin{prop} \label{TrivialExtension}
The functor $\xymatrix@1{\bbH\co \mathbf{Cat} \ar[r] &
\mathbf{DblCat}}$ preserves and reflects weak equivalences,
fibrations, and cofibrations. In other words, a functor $F$ is a
weak equivalence (respectively fibration, respectively cofibration)
in the categorical model structure on $\mathbf{Cat}$ if and only if
$\bbH F$ is a weak equivalence (respectively fibration, respectively
cofibration) in the trivial model structure on $\mathbf{DblCat}$
that arises from the 2-category $\mathbf{DblCat_h}$. As a
consequence, the trivial model structure on $\mathbf{DblCat}$ from
$\mathbf{DblCat_h}$ extends the categorical model structure on
$\mathbf{Cat}$ as a horizontally embedded subcategory.
\end{prop}
\begin{pf}
The 2-functor $\xymatrix@1{\bbH\co\mathbf{Cat} \ar[r] &
\mathbf{DblCat}}$ is fully faithful in the 2-categorical sense, so
it preserves and reflects equivalences (=weak equivalences).

The functor $\bbH$ preserves and reflects isofibrations because of
the characterization of isofibrations at the beginning of Section
\ref{trivialmodelstructure} and the fact that horizontal natural
transformations between functors $\xymatrix@1{\bbE \ar[r] & \bbH
\bfC}$ are in bijective correspondence with natural transformations
between the underlying 1-functors.

Cofibrations are preserved by $\bbH$ by Corollary
\ref{TrivialQuillenAdjunction}. The functor $\bbH$ reflects
cofibrations as follows. If $F$ is a functor such that the double
functor $\bbH F$ is a cofibration, and $G$ is an acyclic fibration
in $\mathbf{Cat}$, then any diagram in $\mathbf{DblCat}$
$$\xymatrix@1{\bbH\bfC \ar[d]_{\bbH F} \ar[r]^{K} & \bbH\mathbf{C'} \ar[d]^{\bbH G} \\
\bbH\bfD \ar[r]_{L} & \bbH\mathbf{D'}}$$ admits a lift, as $\bbH G$
is an acyclic fibration by the above. Since $\bbH$ is fully
faithful, this lift gives us a lift in $\mathbf{Cat}$. Hence the
functor $F$ is a cofibration as well.
\end{pf}

After the discussion of model structures on {\bf DblCat} as a
category of internal categories in Section \ref{section:topologies},
we now turn to a model structure on {\bf DblCat} as a category of
algebras and show that this model structure is the same as the
categorically surjective model structure. We will make use of the
trivial model structure on the underlying category of the 2-category
$\Cat(\mathbf{Graph})$.

\section{A Model Structure for {\bf DblCat} as the 2-Category of
Algebras for a 2-Monad} \label{2monadstructure}

Every 2-category of strict algebras over a 2-monad $T$ with rank
(\ie which preserves $\alpha$-filtered colimits for some $\alpha$)
on a locally finitely presentable 2-category $\mathcal{K}$ admits a
canonical cofibrantly generated {\bf Cat}-enriched model structure
as in \cite{lack2monads}. It is obtained by transferring the trivial
model structure on the 2-category $\mathcal{K}$ described in Section
\ref{trivialmodelstructure}. A strict morphism of strict
$T$-algebras is a weak equivalence (respectively fibration) if and
only if its underlying morphism is an equivalence (respectively
isofibration). We prove that the model structure induced by the
categorically surjective topology $\tau'$ can be recovered in this
way.

The interest in having these two different descriptions of the $\tau'$-model
structure lies in the fact that they allow a characterization of the flexible
double categories (Corollary \ref{flexibledoublecategories} and Remark
\ref{alternativetocorollary}). We will see that the cofibrant replacement in
the $\tau'$-model structure of Remark \ref{tau'cofibrantreplacement} is left
2-adjoint to the inclusion of strict algebras, strict morphisms, and 2-cells
into strict algebras, pseudo morphisms, and 2-cells. In particular, the
cofibrant replacement in Remark \ref{tau'cofibrantreplacement} coincides with
the cofibrant replacement in the algebra structure of \cite{lack2monads}.
Another interesting aspect of the two descriptions of the $\tau'$-model
structure is that {\bf DblCat} provides a good setting for comparing the
categorical model structure on {\bf 2-Cat} in \cite{lack2Cat} and
\cite{lackBiCat} to a model structure induced by a 2-monad.

Recall that the adjunction $\mathbf{Graph} \dashv \mathbf{Cat}$
induces a Cartesian monad $M$ on the category $\mathbf{Graph}$ of
small directed graphs. These directed graphs are non-reflexive; a
choice of distinguished arrows called identities is {\it not} part
of the data. The algebras of $M$ are precisely the small categories.
A {\it Cartesian monad} is a monad whose underlying functor
preserves pullbacks and whose unit and multiplication are Cartesian
natural transformations, \ie each of their naturality squares is a
pullback square. Since $M$ is Cartesian, it induces a 2-monad
$\overline{M}$ on $\Cat(\mathbf{Graph})$, the 2-category of internal
categories in $\mathbf{Graph}$. The strict algebras for this 2-monad
$\overline{M}$ are pairs $(\bbD_0,\bbD_1)$ of $M$-algebras with
source, target, unit, and composition maps compatible with the
$M$-algebra structure. Thus, the strict $\overline{M}$-algebras are
precisely double categories. Similarly, strict morphisms and 2-cells
for strict $\overline{M}$-algebras are double functors and
horizontal natural transformations. The 2-category of strict
$\overline{M}$-algebras, strict morphisms, and 2-cells is
$$\overline{M}\text{-Alg}_\text{s}=\Cat(\mathbf{Cat})=\mathbf{DblCat_h}.$$
The term {\it algebra} will always mean {\it strict algebra}, so we
occasionally leave off the adjective {\it strict}.

Let $\xymatrix@1{U\co \mathbf{Cat} \ar[r] & \mathbf{Graph}}$ be the
forgetful functor. The functor $U$ induces a 2-functor
$$\xymatrix{\overline{U}\co\Cat(\mathbf{Cat}) \ar[r] & \Cat(\mathbf{Graph})}$$
which coincides with the forgetful 2-functor
$$\xymatrix{\overline{M}\text{-Alg}_\text{s} \ar[r] & \Cat(\mathbf{Graph})}.$$
An internal category in $\mathbf{Graph}$ is a (non-reflexive) double
graph $\bbE$ with a category structure on $(\Obj \bbE,\Hor \bbE)$
and on $(\Ver \bbE,\Sq \bbE)$, in other words horizontal
compositions are defined but vertical compositions are not. There
are no vertical identity 1-arrows in $\bbE$, and no vertical
identity squares on horizontal 1-morphisms.

\begin{thm} \label{2monadstructure=tau'structure}
The model structure induced by the 2-monad $\overline{M}$ is the
$\tau'$-model structure.
\end{thm}
\begin{pf}
First we prove that the weak equivalences are the same. Note that a
double functor $G$ is fully faithful if and only if $\overline{U}G$
is fully faithful as in Definition \ref{fullyfaithful}. A double
functor $G$ is a weak equivalence as a morphism of algebras if and
only if $\overline{U}G$ is a weak equivalence in the trivial model
structure on $\Cat(\mathbf{Graph})$, which is the case if and only
if $\overline{U}G$ is fully faithful and there exists a morphism $q$
of directed graphs such that $U(s \circ \overline{G}_0) \circ q =
\mbox{id}_{\bbB_0}$ by Theorem
\ref{trivialtopologyproducestrivialstructure}. That is equivalent to
$G$ being fully faithful and $s \circ \overline{G}_0$ being
$U$-split, which is precisely the definition of weak equivalence in
the $\tau'$-model structure using Proposition \ref{tau'epis} and
Lemma \ref{Usplitcharacterization}. Hence the weak equivalences
coincide.

Similarly, a double functor $G$ is a fibration as a morphism of
algebras if and only if $\overline{U}G$ is a fibration in the
trivial model structure on $\Cat(\mathbf{Graph})$, which is the case
if and only if there exists a morphism $q$ of directed graphs such
that $(U(r_G)_0) \circ q=id$ in Diagram (\ref{fibrationdefinition}),
which is the case if and only if $(r_G)_0$ is $U$-split. This is
equivalent to $G$ being a fibration in the $\tau'$-model structure.
Hence the fibrations coincide.
\end{pf}

\begin{cor}
The categorically surjective $\tau'$-model structure is cofibrantly
generated and admits an enrichment as a {\bf Cat}-enriched model
category.
\end{cor}

For a 2-monad $T$ on $\mathcal{K}$ as above, let
$T$-$\text{Alg}_\text{s}$ denote the 2-category of strict
$T$-algebras, strict morphisms, and 2-cells. As usual, we denote by
$T$-Alg the 2-category of strict $T$-algebras, pseudo morphisms, and
2-cells. As shown in \cite{blackwellkellypower}, the inclusion
$\xymatrix@1{T\text{-Alg}_\text{s} \ar[r] & T\text{-Alg}}$ admits a
left 2-adjoint denoted $A \mapsto A'$. The counit component $\xymatrix@1{q\co
A' \ar[r] & A}$ is a strict morphism, and if $q$ admits a section in
$T\text{-Alg}_\text{s}$, then $A$ is called {\it flexible}. The
flexible algebras are the closure under flexible colimits of the
free algebras. Strict morphisms from $A'$ to $B$ are in bijective
correspondence with pseudo morphisms from $A$ to $B$.

\begin{thm}[Theorem 4.12 in \cite{lack2monads}]
\label{flexiblealgebras} The cofibrant objects of $T\text{\rm
-Alg}_\text{\rm s}$ are precisely the flexible algebras; in
particular, any algebra of the form $A'$ is cofibrant, and is thus a
cofibrant replacement for $A$. Every free algebra is flexible.
\end{thm}

\begin{cor} \label{flexibledoublecategories}
The cofibrant objects in the $\tau'$-model structure are precisely
the flexible double categories. In particular, a double category
$\bbX$ is flexible if and only if $\bbX_0$ is a free category on a directed graph.
\end{cor}
\begin{pf}
This follows from Theorem \ref{2monadstructure=tau'structure}, Theorem
\ref{flexiblealgebras}, and Proposition \ref{tau'cofibrantobjects}.
\end{pf}

We next show that the cofibrant replacement in the $\tau'$-model
structure of Remark \ref{tau'cofibrantreplacement} is the same as
the left 2-adjoint $\bbA \mapsto \bbA'$ to the inclusion
$\xymatrix@1{\overline{M}\text{-Alg}_\text{s} \ar[r] &
\overline{M}\text{-Alg}}$. Our method is to verify that the
cofibrant replacement has the same 2-universal property as $\bbA'$.
For this we need to identify the pseudo morphisms of strict
$\overline{M}$-algebras.

\begin{prop}
The pseudo morphism between strict $\overline{M}$-algebras are the
pseudo double functors which are strict in the horizontal direction,
but weak in the vertical direction.
\end{prop}
\begin{pf}
Let $\bbA$ and $\bbB$ be strict algebras for the 2-monad
$\overline{M}$, in other words, $\bbA$ and $\bbB$ are internal
categories in {\bf Graph} equipped with internal functors
$\xymatrix@1{a\co\overline{M}(\bbA)\ar[r] & \bbA}$, and
$\xymatrix@1{b\co\overline{M}(\bbB)\ar[r] & \bbB}$ defining the
vertical composition. More specifically, $a_0$ and $b_0$ define the
vertical compositions of paths of vertical arrows, and $a_1$ and
$b_1$ define the vertical composition of vertical paths of squares.

A pseudo morphism $\xymatrix@1{(F,\varphi)\co \bbA\ar[r] & \bbB}$
consists of an internal functor $\xymatrix@1{F\co \bbA \ar[r]
&\bbB}$ in $\mathbf{Graph}$ and an invertible internal
transformation
$$
\xymatrix@C=4pc{ \overline{M}(\bbA)\ar[r]^{\overline{M}(F)} \ar[d]_a
\ar@{}[dr]|{\varphi}
    & \overline{M}(\bbB)\ar[d]^b
\\
\bbA\ar[r]_{F} & \bbB,
}
$$
given by a morphism of graphs $\xymatrix@1{\varphi\co
(\overline{M}(\bbA))_0\ar[r] & \bbB_1}$ which satisfies the
coherence conditions from \cite{blackwellkellypower}. It follows
from the identity coherence condition that $\varphi$ sends every
object of $\bbA$ to the corresponding horizontal identity arrow. To
every path $\langle u_1,\ldots,u_n\rangle$ of compatible vertical
arrows in $\bbA$, $\varphi$ assigns a horizontally invertible square
in $\bbB$ denoted by
$$
\xymatrix@C=4pc{
FA_0\ar[d]_{Fu_1}\ar@{=}[r]\ar@{}[3,1]|{\varphi_{u_1,\ldots,u_n}} &
FA_0\ar[ddd]^{F(u_n\circ\cdots\circ u_1)}
\\
FA_1\ar@{..>}[d] &
\\
FA_{n-1}\ar[d]_{Fu_n} &
\\
FA_n\ar@{=}[r] & FA_n \rlap{\,.}}
$$
It follows from the second coherence axiom that for any path of
paths $\langle\langle u_{11},\ldots,u_{1m_1}\rangle,\langle
u_{21},\ldots,u_{2m_2}\rangle,\ldots,\langle u_{n1},\ldots,
u_{nm_n}\rangle\rangle$, the pasting
$$
\xymatrix@C=12pc{
\ar[d]_{Fu_{11}}\ar@{=}[r]\ar@{}[dddr]|{\varphi_{u_{11},\ldots,u_{1m_1}}}
& \ar[ddd]|{\tb{F(u_{1m_1}\circ \cdots \circ u_{11})}}
    \ar@{=}[r] \ar@{}[7,1]|{\varphi_{(u_{1m_1} \circ \cdots \circ u_{11}),\ldots,(u_{nm_n}\circ\cdots\circ u_{n1})}}
    & \ar[7,0]|(.6){\tb{F((u_{nm_n}\circ\cdots\circ u_{n1})\circ \cdots \circ (u_{1m_1} \circ \cdots \circ u_{11}))}}
\\
\ar@{..}[d] &&
\\
\ar[d]_{Fu_{1m_1}} &&
\\
\ar@{=}[r] \ar@{..}[d] \ar@{}[dr]|\vdots  & \ar@{..}[d]&
\\
\ar[d]_{Fu_{n1}}\ar@{=}[r]\ar@{}[dddr]|{\varphi_{u_{n1},\ldots,u_{nm_n}}}
& \ar[ddd]|{\tb{F(u_{nm_n} \circ \cdots \circ u_{n1})}}&
\\
\ar@{..}[d] &&
\\
\ar[d]_{Fu_{nm_n}} &&
\\
\ar@{=}[r] &\ar@{=}[r]&
}
$$
is equal to $\varphi_{u_{11},\ldots,u_{1m_1}, \ldots, u_{n1},\ldots,
u_{nm_n}}$.

Note that the coherence conditions imply that $\varphi$ is
completely determined by its components $\varphi_{u_1,u_2}$ for
composable vertical arrows $u_1$ and $u_2$ in $\bbA$. These
$\varphi_{u_1,u_2}$ are the composition coherence isomorphisms for
the underlying vertical pseudo functor of $F$.

The internal natural transformation $\varphi$ associates to the
empty path on $A$ the unit coherence isomorphism of the underlying
vertical pseudo functor of $F$. The coherence conditions on
$\varphi$ contain the coherence conditions for the coherence
isomorphisms of a pseudo functor.

Conversely, given a pseudo double functor (weak in the vertical
direction) the natural isomorphism $\varphi$ is defined in terms of
the coherence isomorphisms.
\end{pf}

\begin{prop} \label{left2adjointcofibrantreplacement}
The cofibrant replacement in the $\tau'$-model structure of Remark
\ref{tau'cofibrantreplacement} is isomorphic to the left 2-adjoint
$\bbA \mapsto \bbA'$ to the inclusion
$\xymatrix@1{\overline{M}\text{\rm-Alg}_\text{\rm s} \ar[r] &
\overline{M}\text{\rm-Alg}}$.
\end{prop}
\begin{pf}
Let $Q$ denote the cofibrant replacement functor defined on objects
in Remark \ref{tau'cofibrantreplacement}. Our task is to present a
natural isomorphism of categories
\begin{equation} \label{left2adjointtoinclusion}
\overline{M}\text{-Alg}_{\text{s}}(Q\bbA,\bbB)\cong\overline{M}\text{-Alg}(\bbA,\bbB)
\end{equation}
for strict $\overline{M}$-algebras (double categories) $\bbA$ and
$\bbB$.

The category $(Q\bbA)_0$ is the free category on the underlying
graph of $\bbA_0$, and the category $(Q\bbA)_1$ is the pullback in
$\mathbf{Cat}$
$$
\xymatrix@C=3em{(Q\bbA)_1 \ar[r] \ar[d]_{(s,t)} & \bbA_1
\ar[d]^{(s,t)}
\\
(Q\bbA)_0 \times (Q\bbA)_0 \ar[r] & \bbA_0\times\bbA_0 \rlap{\,.}}
$$
In particular, squares of $Q\bbA$ (the morphisms of $(Q\bbA)_1$)
have the form $(\langle u_1,\ldots,u_m\rangle, \alpha, \langle
v_1,\ldots,v_n\rangle)$, where $\alpha$ is a square in $\bbA$ with
the vertical morphism $u=u_m\circ\cdots u_1$ as its horizontal
source and the vertical morphism $v=v_n\circ\cdots\circ v_1$ as its
horizontal target.

Given a double functor $\xymatrix@1{G\co Q\bbA \ar[r] & \bbB}$, we
may define a pseudo morphism $\xymatrix@1{(F,\phi)\co\bbA \ar[r] &
\bbB}$ as follows. On objects, $F(A)=G(A)$; on horizontal morphisms,
$F(f)=G(f)$; on vertical morphisms, $F(u)=G(\langle u \rangle)$; and
on squares, $F(\alpha)=G(\langle u\rangle,\alpha,\langle v\rangle)$,
where $u$ is the horizontal source of $\alpha$ and $v$ is the
horizontal target of $\alpha$. Further, $\phi$ has components
$\phi_{u_1,\ldots,u_m}=G(\langle u_1,\ldots,
u_m\rangle,i_{u_m\circ\cdots\circ u_1}^h,\langle u_m\circ\cdots\circ
u_1\rangle)$ (where $i_{u_m\circ\cdots\circ u_1}^h$ is the
horizontal identity square on the vertical morphism
$u_m\circ\cdots\circ u_1$).

Given a pseudo morphism $\xymatrix@1{(F,\varphi)\co\bbA \ar[r] &
\bbB}$, we may define a double functor $\xymatrix@1{G\co Q\bbA
\ar[r] & \bbB}$ as follows. On objects, $G(A)=F(A)$; on horizontal
morphisms, $G(f)=F(f)$; on vertical morphisms, $$G(\langle
u_1,\ldots, u_m\rangle)= F(u_m)\circ\cdots\circ F(u_1);$$ and on
squares $$G(\langle u_1,\ldots,u_m\rangle, \alpha, \langle
v_1,\ldots,v_n\rangle)=[\phi_{u_1,\ldots,u_m} \hspace{.4cm} F\alpha
\hspace{.4cm} (\varphi_{v_1, \ldots, v_m})^{-1} ].$$

It is straightforward to see that these two procedures are inverse
to each other, and define a bijection on the object sets of equation
(\ref{left2adjointtoinclusion}). We extend this bijection to an
isomorphism of categories. If $\xymatrix@1{\theta\co G_1 \ar@{=>}[r]
& G_2}$ is a horizontal natural transformation of double functors,
then we define $\xymatrix@1{\lambda \co F_1 \ar@{=>}[r] & F_2}$ by
$\lambda_A=\theta_A$ and $$\lambda_{\langle u_1, \ldots, u_m
\rangle}=[(\varphi_1)_{u_1, \ldots, u_m} \hspace{.4cm} \theta_{u_m
\circ \cdots \circ u_1} \hspace{.4cm} (\varphi_2)_{u_1, \ldots,
u_m}^{-1} ].$$ Conversely, if $\xymatrix@1{\lambda \co F_1
\ar@{=>}[r] & F_2}$ is a horizontal natural transformation of pseudo
morphisms, then we define $\xymatrix@1{\theta\co G_1 \ar@{=>}[r] &
G_2}$ by $\theta_A=\lambda_A$ and $\theta_u=\lambda_{\langle u
\rangle}$.

Naturality can be easily verified.
\end{pf}

\begin{rmk} \label{alternativetocorollary}
Proposition \ref{left2adjointcofibrantreplacement} allows us to give an alternative proof of Corollary
\ref{flexibledoublecategories}. In fact, if $\bbA$ is a cofibrant object in the
$\tau'$-model structure, the map $\xymatrix@1{\emptyset \ar[r] & \bbA}$ has the left
lifting property with respect to the counit component $\xymatrix@1{q\co Q\bbA
\ar[r] & \bbA}$, which is an acyclic fibration by Lemma 5.14 in
\cite{everaertinternal}. It follows easily that $q$ has a section,
and thus $\bbA$ is flexible. Conversely, if $\bbA$ is flexible, then
$q$ has a section $p$. In the diagram below,
let $f$ be an acyclic fibration in the $\tau'$-model structure and $h$ any map.
$$
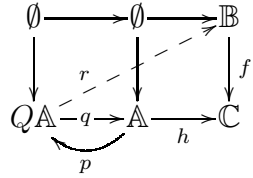
\xymatrix@1{\emptyset \ar[r] \ar[d] & \emptyset \ar[r] \ar[d] & \bbB \ar[d]^f \\ Q\bbA \ar[r]|{\lr{q}} \ar@{-->}[urr]^(.3)r
& \bbA \ar[r]_h \ar@/^1pc/[l]^p & \bbC}$$
Since $Q\bbA$ is cofibrant, there is a map
$\xymatrix@1{r\co Q\bbA \ar[r] & \bbB}$ with $fr=hq$. Hence the map
$\xymatrix@1{rp\co \bbA \ar[r] & \bbB}$ satisfies $frp=hq p=h$. This
shows that $\xymatrix@1{\emptyset \ar[r] & \bbA}$ is a cofibration; that is, $\bbA$ is
cofibrant.
\end{rmk}

The categorical model structure on {\bf 2-Cat} of \cite{lack2Cat}
and \cite{lackBiCat} has weak equivalences the strict 2-functors
that are biequivalences, fibrations those strict 2-functors with the
equivalence lifting property (as defined in \cite{lackBiCat}, not
\cite{lack2Cat}), and cofibrations those strict 2-functors whose
underlying functor has the left lifting property with respect to
functors that are surjective on objects and full. We can compare
this with the $\tau'$-model structure as follows.

\begin{prop} \label{embedding2-Catintau'}
Consider {\bf 2-Cat} vertically embedded in {\bf DblCat}. If a
2-functor is a cofibration in the $\tau'$-model structure on {\bf
DblCat}, then it is a cofibration in the categorical model structure
on {\bf 2-Cat}. A 2-category is cofibrant in the categorical model
structure on {\bf 2-Cat} if and only if it is cofibrant in the
$\tau'$-model structure on {\bf DblCat}. Thus a 2-category is
flexible as in \cite{lack2Cat} if and only if it is flexible as an
algebra over the 2-monad $\overline{M}$.
\end{prop}
\begin{pf}
Suppose $G$ is a 2-functor such that $\bbV G$ is a cofibration in
the $\tau'$-model structure on $\mathbf{DblCat}$. Then $(\bbV G)_0$
has the left lifting property with respect to $\tau'$-functors by
Proposition \ref{cofibration}. This implies that $(\bbV G)_0$ has
the left lifting property with respect to all functors that are
surjective on objects and full by Proposition \ref{tau'epis}. The
underlying functor of $G$ is $(\bbV G)_0$, so $G$ is a cofibration
in the categorical structure on {\bf 2-Cat}.

A 2-category is cofibrant in the categorical structure on {\bf
2-Cat} if and only if its underlying category is projective with
respect to all functors that are surjective on objects and full. But
this coincides with cofibrant 2-categories in the $\tau'$-model
structure on {\bf DblCat} by Corollary \ref{cofibrantobject} and Corollary \ref{sameprojectives}.
\end{pf}

\begin{rmk}
Our characterization of flexible double categories extends Lack's
characterization of flexible 2-categories as those
2-cat\-e\-go\-ries with underlying category a free category on a
graph (Theorem 4.8 (iv) of \cite{lack2Cat}). See Corollary
\ref{flexibledoublecategories} and Proposition
\ref{embedding2-Catintau'}.
\end{rmk}

The sets of weak equivalences with source and target 2-categories in
the two model structures have nontrivial intersection, but neither
set of weak equivalences is contained in the other. The $\bbV$-image
of a biequivalence is not necessarily essentially
$\tau'$-surjective, though it is fully faithful. For example,
consider the inclusion $\xymatrix@1{h\co \{0\} \ar[r] &
\{0\cong1\}}$. Then $h$ is a biequivalence, and the only object of
$(\bbP_{\bbV h})_0$ in Definitions \ref{mappingpathobjectdefinition}
and \ref{fullyfaithful} is $(0,\xymatrix@1{0 \ar[r]^{id_0} & 0})$.
Then $\delta_0 \circ \overline{\bbV h}_0$ cannot be surjective on
objects, as its target has two objects. By Proposition
\ref{tau'epis}, the functor $\delta_0 \circ \overline{\bbV h}_0$ is
not $\tau'$-epi. Thus $\bbV h$ is not essentially
$\tau'$-surjective.

For another reason why the left adjoint $\bbV$ in Proposition
\ref{embedding2-Catintau'} is not a left Quillen functor, consider
the 2-functor $j_1'$ given by inclusion of the terminal 2-category
\{1\} into the free-living adjoint equivalence $E'$. The free-living
adjoint equivalence $E'$ has objects 0 and 1. Morphisms are
$$\xymatrix@1{f\co 0 \ar[r] & 1}$$
$$\xymatrix@1{g\co 1 \ar[r] & 0},$$ as well as all concatenations of
$f$ and $g$. There is a unique 2-cell between every parallel pair of
morphisms. In particular, every 2-cell of $E'$ is invertible. The
2-functor $j_1'$ is a generating acyclic cofibration for the
categorical structure on {\bf 2-Cat} as described in
\cite{lackBiCat}. However $\bbV j_1'$ is not a cofibration in the
$\tau'$-model structure on $\mathbf{DblCat}$: its underlying functor
$(\bbV j_1')_0$ of object categories does not have the left lifting
property with respect to all $\tau'$-epis. For example, let $\bfC$
be the smallest category containing the underlying category
$(E')_0=(\bbV E')_0$ of $E'$ as well as an additional object $1'$
and an arrow $\xymatrix@1{1' \ar[r] & 0}$. The projection from
$\bfC$ to $E'$ takes $1'$ to $1$ and is the identity on $E'$, hence
it is $\tau'$-epi. Then the commutative diagram
$$\xymatrix{(\bbV\{1\})_0 \ar[r]^-{1 \mapsto 1'} \ar[d]_{(\bbV j_1')_0} & \bfC \ar[d]
\\ (\bbV E')_0 \ar[r]_{id} & (\bbV E')_0 }$$
does not admit a lift. Thus, $\bbV$ in Proposition
\ref{embedding2-Catintau'} preserves neither cofibrations nor weak
equivalences.

It is interesting to note that the {\bf Cat}-analogue of Theorem
\ref{2monadstructure=tau'structure} does not hold. In other words,
if we view {\bf Cat} as the category of algebras over the 2-monad
$M$ on {\bf Graph}, then the associated model structure on {\bf Cat}
is not the model structure associated to the topology of surjective
functions on {\bf Set}. A covering family in a basis for this
topology is a single surjective function, so that the epis for this
topology are the same as the epis for the trivial topology by
Proposition \ref{epicharacterization} and Remark
\ref{epicharacterizationremark}, namely the surjective maps
themselves. In fact, the trivial topology, simplicially surjective
topology, and categorically surjective topology on {\bf Set} all
give rise to the categorical model structure on {\bf Cat}, while the
2-monad structure on {\bf Cat} has weak equivalences the
isomorphisms of categories. When we pass to {\bf DblCat} on the
other hand, the three model structures associated to these three
topologies become distinct, and one of them agrees with the 2-monad
structure.

\section{Appendix: Horizontal Nerves and Pushouts}
\label{section:appendix} Though the horizontal nerve and
bisimplicial nerve preserve filtered colimits, they certainly do not
preserve general colimits, not even pushouts. The purpose of this
appendix is to explicitly describe the behavior of the horizontal
nerve on pushouts in {\bf DblCat} along $$\xymatrix@1{i\boxtimes
1_{\bfC}\co\bfA \boxtimes \bfC \ar[r] & \bfB \boxtimes \bfC}$$ where
$\xymatrix@1{i\co\bfA \ar[r] & \bfB}$ is either of the following
full inclusions from Section \ref{structuresonCat}.
$$\xymatrix{c\Sd^2\Lambda^k[m] \ar[r] & c\Sd^2 \Delta[m]}$$
$$\xymatrix@1{\{1\} \ar[r] & \bfI}$$
Theorem \ref{pushoutiso} is the main technical result needed for an
application of Kan's Lemma on Transfer \ref{Kan} to transfer model
structures across the adjunction $c_h \dashv N_h$ in Theorems
\ref{diagramtransfer} and \ref{categoricaltransfer}. In the
following, we use ``$\backslash$'' to denote set-theoretic
complement. We begin with some pushouts in {\bf Cat} which will aid
us in our description of the horizontal and vertical 1-categories of
the pushouts in Theorem \ref{pushoutinDblCat}. The squares will
require an induction argument.

\begin{lem} \label{pushoutinSet}
If $A \subseteq B$ and $D$ are sets, then the pushout in {\bf Set}
$$\xymatrix{A \ar@{^{(}->}[d] \ar[r] & D \ar[d] \\ B \ar[r] & P}$$ is
$P=D \coprod (B \backslash A).$
\end{lem}

\begin{lem} \label{pushoutinCat}
Suppose $\bfA$ is a full subcategory of $\bfB$ and
$$\xymatrix{\bfA \ar[r]^F \ar@{^{(}->}[d] & \bfD \ar[d] \\ \bfB \ar[r] & \bfP}$$
is a pushout in {\bf Cat}. Then the objects of $\bfP$ are
$$
\Obj \bfP=\Obj\bfD \coprod (\Obj\bfB \backslash \Obj \bfA)$$ and
morphisms of $\bfP$ have two forms:
\begin{enumerate}
\item
A morphism $\xymatrix@1{B_0\ar[r]^{f} & B_1}$ with $f \in (\Mor \bfB
\backslash \Mor \bfA)$.
\item
A path $\xymatrix@1@C=3pc{X_1 \ar[r]^{f_1} & D_1 \ar[r]^{d} &
\ar[r]^{f_2} D_2 & X_2 }$ where $d$ is a morphism in $\bfD$, and
$f_1,f_2 \in (\Mor \bfB \backslash \Mor \bfA) \cup
\{\text{identities on } \Obj\bfP\}.$ If $f_1$ is nontrivial, then
$D_1 \in \bfA$. If $f_2$ is nontrivial, then $D_2 \in\bfA$.
\end{enumerate}
\end{lem}
\begin{pf}
To calculate a pushout of categories, one takes the free category on
the pushout of the underlying graphs, and then mods out by the
relations necessary to make the natural maps from $\bfA,\bfB,\bfD$
to the free category into functors as in Theorem \ref{catcolimit}.
Thus the objects of $\bfP$ are
$$
\Obj\bfD \coprod (\Obj\bfB \backslash \Obj\bfA)
$$
by Lemma \ref{pushoutinSet}. The edges of the pushout graph are
$$
\Mor \bfD \coprod (\Mor \bfB \backslash \Mor \bfA),
$$
again by Lemma \ref{pushoutinSet}. The free category on this
consists of finite composable paths of these edges.

Suppose
$$
\xymatrix@C=3pc{P_0 \ar[r]^{f_1} & P_1 \ar[r]^-{f_2} & P_2
\ar@{.}[r] & P_{k-1} \ar[r]^-{f_k} & P_k}
$$ is a morphism in the pushout $\bfP$. Then we can reduce it to
the form (i) or (ii) using the relations induced by $\bfA,\bfB,$ and
$\bfD$ as follows. Suppose $f_{j-1}$ and $f_{j+1}$ are in $\Mor
\bfD$, while $f_j$ is in $(\Mor  \bfB \backslash \Mor \bfA)$. Then
$P_{j-1}$ and $P_j$ must be objects of $\bfA$. But by the fullness
of $\bfA$, $f_j$ must be in $\Mor  \bfA$, and we have arrived at a
contradiction. Thus no morphism of $(\Mor  \bfB \backslash \Mor
\bfA)$ can be surrounded by morphisms of $\bfD$: there exist $0\leq
m \leq n \leq k+1$ such that for all $0 \leq j \leq m$ and all  $n
\leq j \leq k$ we have $f_j \in (\Mor \bfB \backslash \Mor \bfA)$,
and for all $m<j<n$ we have $f_j \in \Mor \bfD$. Next we compose the
$f_j$ in each range, and we obtain a path of the form (i) or (ii).
\end{pf}

\begin{rem}
{\em A morphism $f$ of $\bfB$ is in $\Mor \bfB \backslash \Mor\bfA$
if and only if its source or target is in $\Obj\bfB \backslash \Obj
\bfA$ by the fullness of $\bfA$ in $\bfB$.}
\end{rem}

\begin{lem} \label{ABsetpushoutinCat}
If $A \subseteq B$ are sets and $\bfC$ and $\bfD$ are categories,
then the pushout in {\bf Cat}
$$\xymatrix{A_{\text{\rm disc}} \times \bfC \ar[r] \ar@{^{(}->}[d]_{i \times 1_{\bfC}}
& \bfD \ar[d] \\ B_{\text{\rm disc}} \times \bfC \ar[r] & \bfP}$$ is
$\bfP=\bfD \coprod ((B \backslash A)_{\text{\rm disc}} \times
\bfC)$. (The subscript {\rm `disc'} means discrete category on a
given set.)
\end{lem}
\begin{pf}
Since $B_{\text{disc}} \times \bfC=A_{\text{disc}} \times \bfC
\coprod ((B \backslash A)_{\text{disc}} \times \bfC)$, the pushout
of the underlying graphs is
\begin{equation} \label{ABsetpushoutinCatequation}
\bfD \coprod ((B \backslash A)_{\text{disc}} \times \bfC)
\end{equation}
by Lemma \ref{pushoutinSet}. The free category on this graph, modulo
the appropriate relations as in Theorem \ref{catcolimit}, is once
again (\ref{ABsetpushoutinCatequation}).

Alternatively, this Lemma also follows easily from Lemma
\ref{pushoutinCat}.
\end{pf}

\begin{lem} \label{XsetpushoutinCat}
Suppose $\bfA$ is a full subcategory of $\bfB$, $C$ is a set, and
$$\xymatrix{\bfA \times C_{\text{\rm disc}} \ar[r]^-F \ar@{^{(}->}[d] & \bfD \ar[d] \\ \bfB \times C_{\text{\rm disc}}
\ar[r] & \bfP}$$
is a pushout in {\bf Cat}. Then the objects of $\bfP$ are
$$\Obj \bfP=\Obj \bfD \coprod ((\Obj \bfB \backslash \Obj \bfA) \times C)$$ and
the morphisms of $\bfP$ have two forms:
\begin{enumerate}
\item
A morphism $\xymatrix@1{(B_0,c)\ar[r]^{f} & (B_1,c)}$ with $c \in C$
and $f=(f',c) \in (\Mor \bfB \backslash \Mor\bfA) \times C$.
\item
A path $\xymatrix@1@C=3pc{X_1 \ar[r]^{f_1} & D_1 \ar[r]^{d} &
\ar[r]^{f_2} D_2 & X_2 }$ where $d$ is a morphism in $\bfD$, and
each of $f_1$ and $f_2$ is either in $\Mor \bfB \backslash \Mor \bfA
\times C$ or an identity morphism.
\end{enumerate}
Moreover, if $f_1$ or $f_2$ is not an identity morphism in (ii),
then the path has one of the two respective forms
$$\xymatrix@C=3pc{(B_1,c_1) \ar[r]^{(f_1',c)} & (A_1,c_1) \ar[r]^-d & D_2 \ar[r]^{f_2} & X_2}$$
$$\xymatrix@C=3pc{X_1 \ar[r]^{f_1} & D_1 \ar[r]^-d & (A_2,c_2) \ar[r]^{(f_2',c_2)} & (B_2,c_2)}$$
where $c_1,c_2 \in C$, $B_1,B_2 \in \Obj \bfB\backslash\Obj\bfA$,
$A_1,A_2 \in\Obj\bfA$, $f_1',f_2'\in\Mor\bfB\backslash\Mor\bfA$, and
$d \in \Mor \bfD$.
\end{lem}
\begin{pf}
This follows from Lemma \ref{pushoutinCat}.
\end{pf}

Let us recall the two full inclusions $\xymatrix@1{i\co \bfA \ar[r]
& \bfB}$ under consideration. The first case in which we are
interested is the full inclusion of posets
$\xymatrix@1{c\Sd^2\Lambda^k[m] \ar[r] & c\Sd^2\Delta[m]}$. Here
$\xymatrix@1{c\co \mathbf{SSet} \ar[r] & \mathbf{Cat}}$ denotes the
fundamental category functor as described in Section
\ref{section:categorification} and $\xymatrix@1{\Sd\co \mathbf{SSet}
\ar[r] & \mathbf{SSet}}$ is the subdivision functor defined in
\cite{goerssjardine} and recalled on Page
\pageref{subdivisionobjects}.

The second full inclusion $\xymatrix@1{i\co \bfA \ar[r] & \bfB}$ of
interest is $\xymatrix@1{\{1\} \ar[r] & \bfI}$. The category $\bfI$
consists of two objects $0$ and $1$ and four morphisms: an
isomorphism between $0$ and $1$, and the identity maps. The discrete
subcategory $\{1\}$ is clearly full.

We can now give an explicit description of pushouts in {\bf DblCat}
along $\xymatrix@1{i\boxtimes1_\bfC\co \bfA \boxtimes \bfC \ar[r] &
\bfB \boxtimes \bfC}$ which we use immediately in Theorem
\ref{pushoutiso} for the transfer.

\begin{thm} \label{pushoutinDblCat}
Let $\xymatrix@1{i\co \bfA \ar[r] & \bfB}$ be either of the
following full inclusions.
$$\xymatrix@1{c\Sd^2\Lambda^k[m] \ar[r] & c\Sd^2
\Delta[m]}$$
$$\xymatrix@1{\{1\} \ar[r] & \bfI}$$ Let $\bfC$ be a
category (\eg the finite ordinal $[n]$), and $\bbD$ a double
category. Then the pushout
$$\xymatrix{\mathbf{A} \boxtimes \mathbf{C} \ar[r]^-F
\ar[d]_{i\boxtimes 1_{\mathbf{C}}} & \bbD \ar[d] \\
\mathbf{B} \boxtimes \mathbf{C} \ar[r] & \mathbb{P}}$$ in {\bf
DblCat} has the following explicit description:
\begin{equation} \label{objectsetpushoutinDblCat}
\Obj \bbP=\Obj \bbD \coprod ((\Obj \bfB \backslash \Obj \bfA )
\times \Obj \bfC)
\end{equation}
\begin{equation} \label{horizontal1catpushoutinDblCat}
(\bfH\bbP)_0=(\bfH\bbD)_0 \coprod ((\Obj \bfB \backslash \Obj\bfA
)_{\text{disc}} \times \bfC)
\end{equation}
\begin{equation} \label{vertical1catpushoutinDblCat}
\aligned \Mor (\bfV\bbP)_0 =\{&\text{paths of the form (i) and (ii)}
\\
& \text{in Lemma \ref{XsetpushoutinCat} with }C=\Obj\bfC
\\
& \text{and } \bfD=(\bfV\bbD)_0 \}.
\endaligned
\end{equation}
Squares of $\bbP$ have two forms:
\begin{enumerate}
\item \label{squarepushoutinDblCati}
A square $\begin{array}{c} \xymatrix{\ar[r] \ar[d]
\ar@{}[dr]|{\beta} & \ar[d] \\ \ar[r] & }
\end{array}$
in $\Sq (\bfB \boxtimes \bfC) \backslash \Sq (\bfA \boxtimes \bfC
)$.
\item \label{squarepushoutinDblCatii}
A vertical path of squares $\begin{array}{c} \xymatrix{\ar[r] \ar[d]
\ar@{}[dr]|{\beta_1} & \ar[d] \\ \ar[r] \ar[d] \ar@{}[dr]|{\delta} &
\ar[d] \\ \ar[r] \ar[d] \ar@{}[dr]|{\beta_2} & \ar[d] \\ \ar[r] & }
\end{array}$
where $\delta$ is a square in $\bbD$ and each of $\beta_1$ and
$\beta_2$ is either a vertical identity square (on a horizontal
morphism) in $\bbP$ or is in $\Sq (\bfB \boxtimes \bfC) \backslash
\Sq (\bfA \boxtimes \bfC )$. Moreover, in the case of
$\xymatrix@1{c\Sd^2\Lambda^k[m] \ar[r] & c\Sd^2 \Delta[m]}$, the
square $\beta_1$ is always a vertical identity square.
\end{enumerate}

Note that $\Sq (\bfB \boxtimes \bfC) \backslash \Sq(\bfA \boxtimes
\bfC )=$
$$\left\{
\begin{array}{c}
\xymatrix{(B,C) \ar[r]^{(1_B,g)} \ar[d]_{(f,1_C)} \ar@{}[dr]|{} &
(B,C') \ar[d]^{(f,1_{C'})} \\ (B',C) \ar[r]_{(1_{B'},g)} & (B',C') }
\end{array}\Biggl\lvert
\begin{array}{c}
g \in \Mor \bfC, \\ f \in \Mor \bfB \backslash \Mor  \bfA
\end{array}
\right\}.$$
\end{thm}
\begin{pf}
We use Theorem~\ref{colimitformulaDblCat}. First we calculate the
pushout $\bbS$ of the underlying double derivation schemes. The
object set $\Obj \bbS=\Obj \bbP$ is the pushout of the object sets,
so (\ref{objectsetpushoutinDblCat}) follows from
Lemma~\ref{pushoutinSet}. The horizontal and vertical 1-categories
of $\bbS$ (and $\bbP$) are the pushouts of the horizontal and
vertical 1-categories, so (\ref{horizontal1catpushoutinDblCat})
follows from Lemma \ref{ABsetpushoutinCat} and
(\ref{vertical1catpushoutinDblCat}) follows from
Lemma~\ref{XsetpushoutinCat}. By Lemma~\ref{pushoutinSet} again, the
pushout of the {\it sets} of squares is
\begin{equation} \label{squaresinS}
\Sq \bbS=\Sq \bbD \coprod (\Sq (\bfB \boxtimes \bfC) \backslash \Sq
(\bfA \boxtimes \bfC)).
\end{equation}
Thus we have calculated the pushout $\bbS$ of the underlying double
derivation schemes, its horizontal and vertical 1-categories
coincide with those of $\bbP$, and they have the form claimed in the
theorem. It only remains to show that the squares of $\bbP$ have the
form claimed in the theorem.

The double category $\bbP$ is the free double category on the double
derivation scheme $\bbS$ modulo the smallest congruence making the
natural morphisms of double derivation schemes from $\mathbf{A}
\boxtimes \mathbf{C},\mathbf{B} \boxtimes \mathbf{C},$ and $\bbD$ to
$\bbP$ into double functors. Squares of $\bbP$ are represented by
allowable compatible arrangements in $\bbS$. To prove that squares
of $\bbP$ have the form (\ref{squarepushoutinDblCati}) or
(\ref{squarepushoutinDblCatii}), it suffices to show that any
allowable compatible arrangement of squares in $\bbS$ can be
transformed into (\ref{squarepushoutinDblCati}) or
(\ref{squarepushoutinDblCatii}) using the relations of the
congruence and the double category associativity, identity, and
interchange axioms. The congruence allows us to compose squares
according to the relations in the double categories $\mathbf{A}
\boxtimes \mathbf{C},\mathbf{B} \boxtimes \mathbf{C},$ and $\bbD$.

We must treat the two inclusions $i$ separately.

Let $\xymatrix@1{i\co \bfA \ar[r] & \bfB}$ be the full inclusion
$\xymatrix@1{c\Sd^2\Lambda^k[m] \ar[r] & c\Sd^2 \Delta[m]}$. Recall
from Page \pageref{subdivisionobjects} that $c\Sd^2\Lambda^k[m]$ and
$c\Sd^2\Delta[m]$ are respectively the posets of nondegenerate
simplices of $\Sd\Lambda^k[m]$ and $\Sd\Delta[m]$, and that there is
a morphism $\xymatrix@1{(u_0, \dots, u_p) \ar[r] & (v_0, \dots,
v_q)}$ in $\bfB$ if and only if
$$\{u_0,\dots, u_p\} \subseteq \{v_0,\dots, v_q\}.$$
Also, an object $(v_0, \dots, v_q)$ of $\bfB$ is in $\bfA$ if and
only {\it all} $v_r$ are faces of $\Lambda^k[m]$. Thus, we see for
any path of composable morphisms in $\bfB$
$$\xymatrix@1{B_0 \ar[r]^{f_1} & B_1 \ar[r]^{f_2} & B_2 \ar@{.}[r] & B_{n-1} \ar[r]^{f_n} & B_n}$$
with $B_j$ not in $\bfA$, all $f_\ell$ and  $B_\ell$ with $\ell \geq
j$ are also not in $\bfA$. Thus, once a path leaves $\bfA$, it
cannot return to $\bfA$. In particular, if $\xymatrix@1{B \ar[r] &
B'}$ is a morphism in $\bfB$ and $B$ is not in $\bfA$, then $B'$ is
also not in $\bfA$. Another useful property of $\bfB$ is that every
morphism has a unique decomposition into irreducibles. These special
features of the posets $\bfA$ and $\bfB$ allow us to put the squares
of $\bbP$ into the desired form (\ref{squarepushoutinDblCati}) or
(\ref{squarepushoutinDblCatii}), as we do now.

Suppose $R$ is an allowable compatible arrangement of squares in
$\bbS$, \ie a representative of a square in $\bbP$. If $R$ consists
entirely of squares in $\bbD$, then it is equivalent to its
composition in $\bbD$, so it has the form
(\ref{squarepushoutinDblCatii}) and we are finished.

So suppose that $R$ contains at least one square in $\Sq \bfB
\boxtimes \bfC \backslash \Sq \bfA \boxtimes \bfC$. Then $R$ has at
least one vertex $(B,C)$ in $\bfB \boxtimes \bfC$ but not in $\bfA
\boxtimes \bfC$, \ie $B$ is in $\bfB$ but not in $\bfA$. Any
horizontal morphism in $R$ with source (respectively target) $(B,C)$
is in $\bfB \boxtimes \bfC$ but not in $\bfA \boxtimes \bfC$, as
$(B,C)$ is not in $\bfA \boxtimes \bfC$. Thus the target
(respectively source) of such a morphism has the form $(B,C')$ and
is also in $\bfB \boxtimes \bfC$ but not $\bfA \boxtimes \bfC$. Any
vertical morphism in $R$ with {\it source} $(B,C)$ is in $\bfB
\boxtimes \bfC$ but not in $\bfA \boxtimes \bfC$, as $(B,C)$ is not
in $\bfA \boxtimes \bfC$. Thus the target of such a vertical
morphism is of the form $(B',C)$ with $B'$ not in $\bfA$ by the
special feature of the posets $\bfA$ and $\bfB$ described in the
preceding paragraph. From the original vertex $(B,C)$ we traverse
down a vertical morphism with source $(B,C)$ if there is one,
otherwise we traverse to the right along a horizontal morphism with
source $(B,C)$. In either case, we arrive at another vertex
$(B_1,C_1)$ which is in $\bfB \boxtimes \bfC$ but not in $\bfA
\boxtimes \bfC$. From this vertex we repeat the procedure, moving
either to the right or down. We continue in this way until we reach
the bottom edge of the allowable compatible arrangement $R$. We
conclude that the entire bottom edge of the diagram consists of
objects and horizontal morphisms in $\bfB \boxtimes \bfC$ but not in
$\bfA \boxtimes \bfC$, and hence not in $\bbD$.

Each of these horizontal morphisms on the bottom edge is the bottom
edge of a square in $\bfB \boxtimes \bfC$ but not in $\bfA \boxtimes
\bfC$, since squares of $\bbD$ only have vertices in $\bbD$ (some
objects of $\bbD$ are identified with objects of $\bfA \boxtimes
\bfC$). Thus, the bottom portion of $R$ looks like Figure
\ref{bottomportion} with all squares in $\bfB \boxtimes \bfC$ but
not in $\bfA \boxtimes \bfC$.
\begin{figure}
\setlength{\unitlength}{.7mm}
\begin{picture}(70,35)
\put(0,0){\line(1,0){60}} \put(0,0){\line(0,1){30}}
\put(10,0){\line(0,1){30}} \put(20,0){\line(0,1){35}}
\put(30,0){\line(0,1){35}} \put(40,0){\line(0,1){30}}
\put(50,0){\line(0,1){32}} \put(60,0){\line(0,1){32}}
\put(0,30){\line(1,0){10}} \put(10,20){\line(1,0){10}}
\put(20,35){\line(1,0){10}} \put(30,30){\line(1,0){10}}
\put(40,25){\line(1,0){10}} \put(50,32){\line(1,0){10}}
\put(5,5){\makebox(0,0)[b]} \put(15,5){\makebox(0,0)[b]}
\put(25,5){\makebox(0,0)[b]} \put(35,5){\makebox(0,0)[b]}
\put(45,5){\makebox(0,0)[b]} \put(55,5){\makebox(0,0)[b]}
\put(70,5){\makebox(0,0)[b]}
\end{picture}
\caption{The bottom portion of $R$.} \label{bottomportion}
\end{figure}
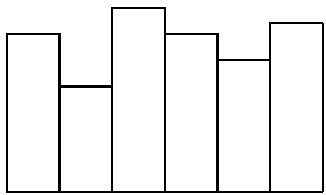
\begin{figure}
\setlength{\unitlength}{.7mm}
\begin{picture}(70,35)
\put(0,0){\line(1,0){60}} \put(0,0){\line(0,1){30}}
\put(10,0){\line(0,1){30}} \put(20,0){\line(0,1){35}}
\put(30,0){\line(0,1){35}} \put(40,0){\line(0,1){30}}
\put(50,0){\line(0,1){32}} \put(60,0){\line(0,1){32}}
\put(0,30){\line(1,0){10}} \put(10,20){\line(1,0){10}}
\put(20,35){\line(1,0){10}} \put(30,30){\line(1,0){10}}
\put(40,25){\line(1,0){10}} \put(50,32){\line(1,0){10}}
\put(5,5){\makebox(0,0)[b]} \put(15,5){\makebox(0,0)[b]}
\put(25,5){\makebox(0,0)[b]} \put(35,5){\makebox(0,0)[b]}
\put(45,5){\makebox(0,0)[b]} \put(55,5){\makebox(0,0)[b]}
\put(70,5){\makebox(0,0)[b]} \linethickness{.75mm}
\put(0,20){\line(1,0){60}}
\end{picture}
\caption{A factorization of the squares in Figure
\ref{bottomportion}.} \label{bottomportioncut}
\end{figure}
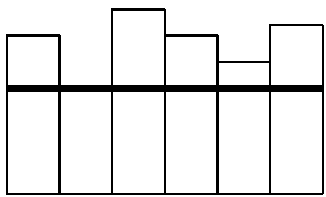

Next we factor the vertical morphisms of Figure \ref{bottomportion}
into irreducibles, which we can do since these vertical morphisms
are of the form $(f,C)$ where $f$ is a morphism in $\bfB$ and $C$ is
an object of $\bfC$. By the uniqueness of the factorization and the
form of squares in $\bfB \boxtimes \bfC$, we can factor these
squares at the height of the shortest one as illustrated in Figure
\ref{bottomportioncut}. We include these new horizontal morphisms
into the allowable compatible arrangement $R$, and obtain a new
compatible arrangement $R_1$. The compatible arrangement $R_1$ is
also allowable, since the same cuts that make $R$ allowable also
make $R_1$ allowable.

The bold horizontal line in Figure \ref{bottomportioncut} is a full
length cut on an allowable compatible arrangement $R_1$, hence it
divides $R_1$ into two allowable compatible arrangements by
Proposition \ref{anycut}. We denote the upper allowable compatible
arrangement by $R_{1,1}$ and the lower allowable compatible
arrangement by $R_{1,2}$. Then $R_{1,1}$ has at least one square
less than $R$, since we cut off at the height of the shortest square
whose bottom edge is on the bottom edge of $R$. If we argue by
induction on the number of squares in an allowable compatible
arrangement, we may assume that $R_{1,1}$ is equivalent to a square
of the form (\ref{squarepushoutinDblCati}) or
(\ref{squarepushoutinDblCatii}). The allowable compatible
arrangement $R_{1,2}$ is equivalent to a square of the form
(\ref{squarepushoutinDblCati}), as it can be composed horizontally.
Finally, we compose $R_{1,1}$ with $R_{1,2}$ to conclude that $R$ is
also equivalent to a compatible arrangement of the form (i) or (ii).

We only need an argument for the triviality of $\beta_1$ whenever a
compatible arrangement is equivalent to one of the form (ii).
Suppose $\beta_1$ is in $\Sq (\bfB \boxtimes \bfC) \backslash \Sq
(\bfA \boxtimes \bfC )$. Then its lower two vertices cannot be in
$\bfA \boxtimes \bfC$ (for if they were, the upper two vertices must
also be in $\bfA \boxtimes \bfC$, and the square $\beta_1$ would be
in $\bfA \boxtimes \bfC$). Thus, the upper two vertices of the
square $\delta$ are not in $\bfA \boxtimes \bfC$, a contradiction.
Thus $\beta_1$ must be a vertical identity. This completes the proof
of Theorem \ref{pushoutinDblCat} for the case
$\xymatrix@1{c\Sd^2\Lambda^k[m] \ar[r] & c\Sd^2 \Delta[m]}$.

Now we turn to the squares in the second case. Let $\xymatrix@1{i\co
\bfA \ar[r] & \bfB}$ be the full inclusion $\xymatrix@1{\{1\} \ar[r]
& \bfI}$ where $\bfI$ is the category with two objects $0$ and $1$
and an isomorphism between them. We will again argue by induction on
the number of squares in the allowable compatible arrangement, but
the special features of the inclusion $\xymatrix@1{\{1\} \ar[r] &
\bfI}$ are different from those of the previous case.
Note that $\bfB \boxtimes \bfC$ only has the four types of squares
listed in Figure \ref{fourtypes}.
\begin{figure}
$$\xymatrix{(0,C_1) \ar[r] \ar@{=}[d]  & (0,C_2) \ar@{=}[d]
& (1,C_1) \ar@{=}[d] \ar[r] & (1,C_2) \ar@{=}[d]
\\ (0,C_1) \ar[r] & (0,C_2) & (1,C_1) \ar[r] & (1,C_2)}$$
$$\xymatrix{(0,C_1) \ar[r] \ar[d] & (0,C_2) \ar[d] & (1,C_1) \ar[d]
\ar[r] & (1,C_2) \ar[d] \\ (1,C_1) \ar[r] & (1,C_2) & (0,C_1) \ar[r]
& (0,C_2)}$$ \caption{The Four Types of Squares in $\bfB \boxtimes
\bfC$.} \label{fourtypes}
\end{figure}
The only vertical morphisms in $\bfB \boxtimes \bfC$ that are
identified with a morphism in $\bbD$ are the trivial vertical
morphisms $\xymatrix@1{\mbox{id}_{(1,C)}^v\co (1,C) \ar[r] &
(1,C)}$.

Suppose that any allowable compatible arrangement of squares in
$\bbS$ with fewer than $n$ squares is equivalent in $\bbP$ to one of
the form (\ref{squarepushoutinDblCati}) or
(\ref{squarepushoutinDblCatii}). Let $R$ be an allowable compatible
arrangement of $n$ squares in $\bbS$. Since $R$ is allowable, it
admits a full length cut $\mathcal{C}$ which divides $R$ into two
allowable compatible arrangements each with fewer than $n$ squares.
We now recombine these two smaller allowable compatible arrangements
to show that $R$ is equivalent to a compatible arrangement of the
form (i) or (ii), but the argument is slightly different depending
on whether $\mathcal{C}$ is horizontal or vertical.

Suppose the full length cut $\mathcal{C}$ is horizontal. Let $R_1$
and $R_2$ be the allowable compatible arrangements above and below
$\mathcal{C}$ respectively. Since $R_1$ and $R_2$ have fewer than
$n$ squares, they must be equivalent to compatible arrangements of
the form (i) or (ii). If $R_1$ and $R_2$ both are equivalent to
compatible arrangements of the form (ii), then by the fullness of
$\bfA$ in $\bfB$ their vertical composite is also of the form (ii),
and hence $R$ is equivalent to a compatible arrangement of the form
(ii). If one or both of $R_1$ and $R_2$ has the form (i), then one
can similarly conclude that $R$ is equivalent to a compatible
arrangement of form (i) or (ii).

Suppose the full length cut $\mathcal{C}$ is vertical. Let $Q^\ell$
and $Q^r$ be the allowable compatible arrangements to the left and
to the right of $\mathcal{C}$ respectively. Since $Q^\ell$ and $Q^r$
have fewer than $n$ squares, they must be equivalent to compatible
arrangements of the form (i) or (ii). There are several cases to
consider.

If both $Q^\ell$ and $Q^r$ are equivalent to compatible arrangements
of the form (i), then their horizontal composite $R$ is clearly in
$\bfB \boxtimes \bfC$, and hence also equivalent to a compatible
arrangement of the form (i) or (ii).

If $Q^\ell$ is equivalent to a compatible arrangement of the form
(i) and $Q^r$ is equivalent to a compatible arrangement of the form
(ii), then $\beta_1^r$ and $\beta_2^r$ must be in $\bfB \boxtimes
\bfC$ as in Figure \ref{i-ii}.
\begin{figure}
$$\xymatrix{ & (B_1,C_1^\ell) \ar[ddd] \ar[r] \ar@{}[dddr]|\beta^\ell & (B_1,C)
\ar[ddd]^{p^\ell} & & (B_1,C) \ar[r] \ar[d]_{j^r}
\ar@{}[dr]|{\beta_1^r}
& \ar[d] & \\
Q^\ell & & & & (B_2^r,C) \ar[d]_{k^r} \ar[r] \ar@{}[dr]|{\delta^r} & \ar[d] & Q^r \\
& & & & (B_3^r,C) \ar[d]_{m^r} \ar[r] \ar@{}[dr]|{\beta_2^r} & \ar[d] & \\
& (B_4,C_1^\ell) \ar[r] & (B_4,C) & & (B_4,C) \ar[r] & & }$$
\caption{ \hspace{1mm}}\label{i-ii}
\end{figure}
Further, the vertical morphism $\xymatrix@1{k^r\co (B_2^r,C) \ar[r]
& (B_3^r,C)}$ must be the identity
$$\xymatrix{\mbox{id}_{(1,C)}^v\co (1,C) \ar[r] & (1,C),}$$ since $k^r=(m^r)^{-1}p^\ell(j^r)^{-1}$
lies in both $\bfB \boxtimes \bfC$ and $\bbD$, and the only vertical
morphisms in both $\bfB \boxtimes \bfC$ and $\bbD$ are such vertical
identities. Then we can subdivide $\beta^\ell$ in $\bfB \boxtimes
\bfC$ as in Figure \ref{i-iisubdivided}.
\begin{figure}
$$\xymatrix{ & (B_1,C_1^\ell) \ar[d] \ar[r] & (B_1,C) \ar[d]
& & (B_1,C) \ar[r] \ar[d] \ar@{}[dr]|{\beta_1^r}
& \ar[d] & \\
Q^\ell & (1,C_1^\ell) \ar@{=}[d] \ar[r] & (1,C) \ar@{=}[d]  & &
(1,C) \ar@{=}[d] \ar[r] \ar@{}[dr]|{\delta^r}
& \ar[d] & Q^r \\
& (1,C_1^\ell) \ar[d] \ar[r] & (1,C) \ar[d] & & (1,C) \ar[d] \ar[r] \ar@{}[dr]|{\beta_2^r} & \ar[d] & \\
 & (B_4,C_1^\ell) \ar[r] & (B_4,C) & & (B_4,C) \ar[r] & & }$$
\caption{\hspace{1mm}} \label{i-iisubdivided}
\end{figure}
The middle square of $Q^\ell$ is now an identity square on a
horizontal morphism in $\bbD$, and hence is also a square in $\bbD$.
Finally, we horizontally compose $Q^\ell$ and $Q^r$ and use the
interchange law to obtain a compatible arrangement of the form (i)
or (ii). Hence $R$ is equivalent to a compatible arrangement of the
form (i) or (ii).

Next we consider the case where $Q^\ell$ and $Q^r$ are both
equivalent to compatible arrangements of the form (ii) {\it and} the
squares $\beta_1^\ell, \beta_2^\ell, \beta_1^r, \beta_2^r$ are in
$\bfB \boxtimes \bfC$ as in Figure \ref{ii-ii}.
\begin{figure}
$$\xymatrix{ &  \ar[d] \ar[r] \ar@{}[dr]|{\beta_1^\ell} & (B_1,C_1)
\ar[d]^{j^\ell} & & (B_1,C_1) \ar[r] \ar[d]_{j^r}
\ar@{}[dr]|{\beta_1^r}
& \ar[d] & \\
Q^\ell &  \ar[d] \ar[r] \ar@{}[dr]|{\delta^\ell} & (B_2^\ell,C_1)
\ar[d]^{k^\ell} & & (B_2^r,C_1) \ar[d]_{k^r} \ar[r]
\ar@{}[dr]|{\delta^r} & \ar[d] & Q^r
\\ &  \ar[d] \ar[r] \ar@{}[dr]|{\beta_2^\ell} & (B_3^\ell,C_2) \ar[d]^{m^\ell} &
& (B_3^r,C_2) \ar[d]_{m^r} \ar[r] \ar@{}[dr]|{\beta_2^r} & \ar[d] & \\
 & \ar[r] & (B_4,C_2) & & (B_4,C_2) \ar[r] & & }$$
\caption{\hspace{1mm}} \label{ii-ii}
\end{figure}
Then $B_2^\ell=1=B_2^r$ and $B_3^\ell=1=B_3^r$, since the only
objects of $\bfB \boxtimes \bfC$ that are identified with an object
of $\bbD$ are of the form $(1,C)$. Thus $j^\ell=j^r$ and
$m^\ell=m^r$, as there is a unique vertical morphism from any object
of $\bfB \boxtimes \bfC$ to another. Since $j^\ell$ and $m^\ell$ are
invertible and $m^\ell k^\ell j^\ell = m^\ell k^r j^\ell$, we see
also that $k^\ell=k^r$. Hence $Q^\ell$ and $Q^r$ can be horizontally
composed to obtain a compatible arrangement equivalent to (i) or
(ii).

Next we consider the case where $Q^\ell$ and $Q^r$ are both
equivalent to compatible arrangements of the form (ii), but the
squares $\beta_1^\ell, \beta_2^\ell, \beta_1^r, \beta_2^r$ may be
vertical identity squares in $\bbP$, \ie not necessarily in $\bfB
\boxtimes \bfC$, as in Figure \ref{ii-ii-identity}.
\begin{figure}
$$\xymatrix{ &  \ar[d] \ar[r] \ar@{}[dr]|{\beta_1^\ell} & P_1
\ar[d]^{j^\ell} & & P_1 \ar[r]^{f^r} \ar[d]_{j^r}
\ar@{}[dr]|{\beta_1^r}
& \ar[d] & \\
Q^\ell &  \ar[d] \ar[r] \ar@{}[dr]|{\delta^\ell} & D_1^\ell
\ar[d]^{k^\ell} & & D_1^r \ar[d]_{k^r} \ar[r] \ar@{}[dr]|{\delta^r}
& \ar[d] & Q^r
\\ &  \ar[d] \ar[r] \ar@{}[dr]|{\beta_2^\ell} & D_2^\ell \ar[d]^{m^\ell} &
& D^r \ar[d]_{m^r} \ar[r] \ar@{}[dr]|{\beta_2^r} & \ar[d] & \\
 & \ar[r] & P_2 & & P_2 \ar[r] & & }$$
\caption{\hspace{1mm}} \label{ii-ii-identity}
\end{figure}
Suppose $\beta_1^\ell$ is a vertical identity square. Then
$P_1=D_1^\ell$. We claim that $\beta_1^r$ is also a vertical
identity square; there are two cases to prove. If $D_1^\ell$ is an
object of $\bbD$ that is not of the form $(1,C)$, then $f^r$ cannot
be in $\bfB \boxtimes \bfC$ (as its source is not in $\bfB \boxtimes
\bfC$). Hence $\beta_1^r$ is a vertical identity square. For the
second case, if $D_1^\ell$ is of the form $(1,C)$, then $j^r$ is a
vertical arrow in $\bfB \boxtimes \bfC$ with source and target
$(1,C)$. By the special form of squares in $\bfB \boxtimes \bfC$ in
Figure \ref{fourtypes}, we see that $\beta_1^r$ is also a vertical
identity square. Thus, we have proved, if $\beta_1^\ell$ is a
vertical identity square, then $\beta_1^r$ is also a vertical
identity square. One can similarly show that if any one of
$\beta_1^\ell, \beta_2^\ell, \beta_1^r, \beta_2^r$ is a vertical
identity square, then the square next to it is also.

Let us continue the case where $Q^\ell$ and $Q^r$ are both
equivalent to compatible arrangements of the form (ii) as in Figure
\ref{ii-ii-identity}, and suppose again that $\beta_1^\ell$ is a
vertical identity square. Then $\beta_1^r$ is also a vertical
identity square. If either of $\beta_2^\ell$ or $\beta_2^r$ is a
trivial identity square, then so is the other, in which case
$Q^\ell$ and $Q^r$ can be horizontally composed to give a compatible
arrangement equivalent to one of the form (i) or (ii). If neither
$\beta_2^\ell$ nor $\beta_2^r$ is a vertical identity square, then
they are both in $\bfB \boxtimes \bfC$, and we can argue as in
Figure \ref{ii-ii} to conclude $(B_4,C_2)=P_2=(B_4,C_2)$,
$m^\ell=m^r$, and $k^\ell=k^r$, in which case $Q^\ell$ and $Q^r$ can
be horizontally composed to give a compatible arrangement equivalent
to one of the form (i) or (ii).

The other cases of Figure \ref{ii-ii-identity}, where one or more of
$\beta_1^\ell, \beta_2^\ell, \beta_1^r, \beta_2^r$ is a vertical
identity square in $\bbP$, are similar.

Thus every every square of $\bbP$ is equivalent to a compatible
arrangement of the form (i) or (ii), for both inclusions $i$ under
consideration. This completes the proof of Theorem
\ref{pushoutinDblCat}.
\end{pf}

The two inclusions of Theorem \ref{pushoutinDblCat} have some
features in common, and the theorem holds for an entire class of
inclusions $\xymatrix@1{i\co \bfA \ar[r] & \bfB}$. We will return to
the this topic and its interaction with
\cite{dawsonparepronkextensions} in the future. Theorem
\ref{pushoutinDblCat} allows us to characterize the behavior of the
horizontal nerve on such pushouts in Theorem \ref{pushoutiso}, which
we need to transfer the model structures from
$\mathbf{Cat}^{\Delta^{op}}$ in Section
\ref{section:simplicialcategories}.

\begin{thm} \label{pushoutiso}
Let $\xymatrix@1{i\co\bfA \ar[r] & \bfB}$ be either of the following
full inclusions.
$$\xymatrix{c\Sd^2\Lambda^k[m] \ar[r] & c\Sd^2\Delta[m]}$$
$$\xymatrix@1{\{1\} \ar[r] & \bfI}$$  Let $\bfC$ be a
finite ordinal $[n]$ viewed as a category, $\bbD$ a double category,
and $\bbP$ the pushout
$$\xymatrix{\mathbf{A} \boxtimes \mathbf{C} \ar[r]
\ar[d]_{i\boxtimes 1_{\mathbf{C}}} & \bbD \ar[d] \\
\mathbf{B} \boxtimes \mathbf{C} \ar[r] & \mathbb{P}}$$ in {\bf
DblCat}.  Then the induced map
\begin{equation*}
N_h(\bbD) \coprod_{N_h(\mathbf{A} \boxtimes \mathbf{C})} N_h
(\mathbf{B}\boxtimes \mathbf{C})  \xymatrix{ \ar[r] &
N_h(\mathbb{P})}
\end{equation*}
is an isomorphism of simplicial objects in {\bf Cat}.
\end{thm}
\begin{pf}
We calculate the pushout
\begin{equation} \label{horizontalnervepushout}
N_h(\bbD) \coprod_{N_h(\mathbf{A} \boxtimes \mathbf{C})} N_h
(\mathbf{B}\boxtimes \mathbf{C})
\end{equation}
levelwise and compare it with $N_h(\bbP)$, which was described in
Theorem \ref{pushoutinDblCat}. The horizontal nerve of an external
product of categories is known from Proposition
\ref{horizontalnerveexternalproduct}.

In level 0, the pushout (\ref{horizontalnervepushout}) is
\begin{equation*}
 \bbD_0\coprod_{\bfA
\times (\Obj \bfC)_{disc}} \bfB \times ( \Obj \bfC)_{\text{disc}}
\end{equation*}
which is the same as the vertical 1-category of $\bbP$ and thus
$(N_h\bbP)_0$.


In level $k \geq 1$ the pushout (\ref{horizontalnervepushout}) is
\begin{equation*}
N_h(\bbD)_k \coprod_{\bfA \times N\bfC_k} \bfB \times N\bfC_k
\end{equation*}
by Proposition \ref{horizontalnerveexternalproduct}. An application
of Lemma 10.5 to level $k$ together with Theorem 10.6 give
immediately that
\begin{equation*}
N_h(\bbD)_k \coprod_{\bfA \times N \bfC_k} (\bfB\times N\bfC_k)
\xymatrix{\ar[r] &} (N_h\bbP)_k
\end{equation*}
is full (and is the identity on objects). To see that this functor
is also faithful, we only need to concern ourselves with the squares
of the form (ii) in Theorem \ref{pushoutinDblCat}. For the Thomason
structure these squares always have a representative where $\beta_2$
is of the form
$$\xymatrix{(A,C)\ar[d]\ar[r]&(A,C') \ar[d]\\ (B,C)\ar[r] & (B,C')}$$
such that $A$ is a maximal element of $c\Sd^2\Lambda^k[n]$. This
determines $\beta_2$ and $\delta$ uniquely. For the categorical
structure all non-identity $\beta_2$ squares
are of the form $$\xymatrix{(1,C)\ar[d]\ar[r] &(1,C')\ar[d]\\
(0,C)\ar[r] &(0,C').}$$ A similar argument for $\beta_1$ shows that
$\beta_1$, $\delta$, and $\beta_2$ are determined uniquely.

\end{pf}


\end{document}